%% file: globalrv-20210214v21-b.tex
\documentclass[a4paper,12pt]{article}
\usepackage{dsfont}
\usepackage{amssymb}
\usepackage{latexsym}
\usepackage{amsmath}
\usepackage{xcolor}
\usepackage{comment}
\usepackage{amsthm}
\usepackage[dvips]{graphicx}
\usepackage{layout} 
\usepackage{ulem}
\usepackage{enumerate}
\usepackage{bm}
\usepackage{bbm} 
\usepackage[bbgreekl]{mathbbol} 
\usepackage{authblk}
\usepackage{setspace} 

\usepackage{here} 

\newif\ifrs
\rstrue
\ifrs \usepackage{mathrsfs} \fi  
\newif\ifcol
\coltrue 

\newtheorem{theorem}{Theorem}[section]

\newtheorem{lemma}[theorem]{Lemma}

\newtheorem{remark}[theorem]{Remark}

\numberwithin{equation}{section}
\newtheorem{theorem*}{Theorem}
\newtheorem{ass*}[theorem*]{Assumption}
\newtheorem{note*}[theorem*]{Note}
\newtheorem{lemma*}[theorem*]{Lemma}
\newtheorem{definition*}[theorem*]{Definition}
\newtheorem{proposition*}[theorem*]{Proposition}
\newtheorem{corollary*}[theorem*]{Corollary}
\newtheorem{remark*}[theorem*]{Remark}
\newtheorem{example*}[theorem*]{Example}
\numberwithin{equation}{section}
\newif\ifcol
\coltrue
\ifcol
\newcommand{\colorr}{\color{black}}
\newcommand{\sred}{\color{black}}
\newcommand{\tred}{\color{black}}
\newcommand{\fred}{\color{black}}
\newcommand{\fblue}{\color{black}}
\newcommand{\vred}{\color{black}}
\newcommand{\xred}{\color{black}}
\newcommand{\colorg}{\color{black}}
\newcommand{\colorb}{\color{black}}

\newcommand{\colorn}{\color[rgb]{1,1,1}}

\newcommand{\coloro}{\color{black}}
\newcommand{\coloroy}{\color[rgb]{1,0.95,0}}

\else
\newcommand{\sred}{\color{black}}
\newcommand{\tred}{\color{black}}
\newcommand{\fred}{\color{black}}
\newcommand{\fblue}{\color{black}} 
\newcommand{\colorb}{\color{black}}

\newcommand{\colorr}{\color{black}}
\newcommand{\colorg}{\color{black}}

\newcommand{\colorn}{\color{black}}
\newcommand{\coloro}{\color{black}}

\newcommand{\coloroy}{\color{black}}
\fi
\excludecomment{en-text}
\includecomment{jp-text}
\includecomment{comment}
\input nakamacro291103+++.tex

\def\mfA{{\mathfrak A}}
\def\mfB{{\mathfrak B}}
\def\mfS{{\mathfrak S}}

\def\tk{{t_k}}
\def\tkm{{t_{k-1}}}
\def\ol{\overline}
\def\wh{\widehat}
\def\wt{\widetilde}
\def\hU{\widehat{U}}
\def\hR{\widehat{R}}

\begin{document}
	
	\title{{\colorb Global jump filters and realized volatility} 
		\footnote{
			This work was in part supported by 
			Japan Science and Technology Agency CREST JPMJCR14D7; 
			Japan Society for the Promotion of Science Grants-in-Aid for Scientific Research 
			No. 17H01702 (Scientific Research);  
			and by a Cooperative Research Program of the Institute of Statistical Mathematics. 
		}
	}
	
	{\colorb 
		\author[1]{Haruhiko Inatsugu}
		\author[1,2]{Nakahiro Yoshida}
		\affil[1]{Graduate School of Mathematical Sciences, University of Tokyo
			\footnote{Graduate School of Mathematical Sciences, University of Tokyo: 3-8-1 Komaba, Meguro-ku, Tokyo 153-8914, Japan. e-mail: nakahiro@ms.u-tokyo.ac.jp}
		}
		\affil[2]{Japan Science and Technology Agency CREST
		}
	}
	\maketitle
	\ \\
{\it Summary} 
\begin{en-text}
	{\colorb
		We propose a new estimation method for integrated volatility of a 
		semimartingale with jumps using a jump-detection filter. 
		We apply the ``global filter" proposed by Inatsugu and Yoshida \cite{InatsuguYoshida2020Accepted} which uses all observations and 
		obtain the order statistic of absolute increments as a threshold for 
		jump detection. 
		We prove the moment convergence and asymptotic mixed normality of 
		our ``global realized volatility" and its variant 
		``winsorized global volatility." 
		By numerical simulation, we show that our estimators 
		outperform previous realized volatility estimators that use 
		a few adjacent increments to mitigate the effects of jumps. 
	}
\end{en-text}
{\colorg
	For a semimartingale with jumps, 
	we propose a new estimation method for integrated volatility, 
	i.e., the quadratic variation of the continuous martingale part, 
	based on the global jump filter proposed by Inatsugu and Yoshida \cite{InatsuguYoshida2020Accepted}. 
	%
	To decide whether each increment of the process has jumps, 
	the global jump filter adopts 
	the upper $\alpha$-quantile of the absolute increments as the threshold. 
	This jump filter is called global since it uses all the observations to classify one increment. 
	We give a rate of convergence and prove asymptotic mixed normality of 
	the global realized volatility and its variant ``Winsorized global volatility''.
	By simulation studies, we show that our estimators 
	outperform previous realized volatility estimators that use 
	a few adjacent increments to mitigate the effects of jumps. 
}

\section{Introduction}
Let $(\Omega,\calf,P)$ be a probability space equipped with 
a filtration $\bbF=(\calf_t)_{t\in[0,T]}$. 
We consider a one-dimensional semimartingle $X=(X_t)_{t\in[0,T]}$ 
having a decomposition 
\bea\label{0301191155}
X_t &=& X_0+\int_0^t b_sds+\int_0^t{\colorb \sigma_s}dw_s+J_t\quad(t\in[0,T])
\eea
where 
$X_0$ is an $\calf_0$-measurable random variable, 
$b=(b_t)_{t\in[0,T]}$ and $\sigma=(\sigma_t)_{t\in[0,T]}$ are 
\cadlag $\bbF$-adapted processes, and 
$w=(w_t)_{t\in[0,T]}$ is an $\bbF$-standard Wiener process. 
$J=(J_t)_{t\in[0,T]}$ is the jump part of $X$. 
We will assumed that $J$ is finitely active, that is, 
$J_t=\sum_{s\in(0,t]}{\tred\Delta J_s}$ for ${\tred\Delta J_s=J_s-J_{s-}}$ 
and $\sum_{t\in[0,T]}1_{\{\Delta J_t\not=0\}}<\infty$ a.s. 
%
%
In this paper, we are interested in the estimation of the integrated volatility
\bea\label{0211281514}
\Theta
&=& 
\int_0^T\sigma_t^2dt
\eea
based on the data $(X_\tj)_{j=0,1,...n}$, where 
$\tj=\tj^n=jT/n$. 

The jump part $J$ can be endogenous or exogenous, as well as $b$ and $\sigma$, 
however, $J$ is a nuisance in any case. 
The simple realized volatility is heavily damaged when jumps exist. 
To avoid the effects of the jumps, various methods have been proposed so far. 
For example, 
the bipower variation (Barndorff-Nielsen and Shephard. \cite{Barndorff-NielsenShephard2004b}, Barndorff-Nielsen et al. \cite{Barndorff-NielsenShephardWinkel2006}) and the minimum 
realized volatility (Andersen et al. \cite{Andersen2012}) are shown to 
be consistent estimators of the integrated volatility even in the presence of jumps. 
The idea of these methods is that, to mitigate the effect of jumps, 
they employ adjacent increments in constructing the estimator.  
Another direction to handle jumps is to introduce a threshold to detect jumps.  

{\tred
	Parametric inference for sampled 
	{\coloro diffusion type processes} was studied by 
	Dohnal \cite{Dohnal1987}, 
	Prakasa Rao \cite{PrakasaRao1983,prakasa1988statistical}, 
	Yoshida \cite{Yoshida1992b, yoshida2011}, 
	Kessler \cite{Kessler1997}, 
	Genon-Catalot and Jacod \cite{Genon-CatalotJacod1993},
	Uchida and Yoshida \cite{UchidaYoshida2013, UchidaYoshida2012Adaptive,uchida2014adaptive}, 
	Ogihara and Yoshida \cite{ogihara2014quasi}, 
	Kamatani and Uchida \cite{KamataniUchida2014}
	and others. 
	{\colorg Limit theorems used to analyse the realized volatility appeared 
		in the studies of parametric inference.} 
	If no jump part exists, then the distribution of the increment $\Delta_jX=X_\tj-X_\tjm$ 
	admits Gaussian approximation in a short time interval, and 
	a quasi-likelihood function can be constructed with the conditional Gaussian density. 
	When the jump part exists, the local Gaussian approximation is no longer valid. 
	Then it is necessary to detect jumps and classify the increments to apply the local Gaussian quasi-likelihood function 
	for estimation of the parameters in the continuous part.}
Threshold method was investigated {\tred by} Shimizu and Yoshida \cite{ShimizuYoshida2006} 
{\sred and Ogihara and Yoshida \cite{OgiharaYoshida2011}} 
in the context of the parametric inference for a stochastic differential equation with jumps. 
The idea of thresholding {\colorg is rather old, going}
back to the studies of limit theorems for L\'evy processes 
as latest. Mancini {\colorb \cite{Mancini2001}} {\tred used this idea in a nonparametric situation. 
	{\colorg Koike \cite{koike2014estimator} applied the threshold method to covariance estimation for asynchronously observed semimartingales with jumps.}
	The classical jump filters compare 
	the size of increment with a threshold determined 
	by a {\colorg (conditionally/unconditionally deterministic)} function 
	of the length of {\colorg the} time interval. 
	If an increment 
	is so large that exceeds the threshold, it is regarded as having jumps. 
	Otherwise, the increment is regarded as having no jump. 
	Once classified, the increments are used to estimate the parameters 
	in continuous and jump parts, respectively.}

{\tred 
	Though the efficiency of the traditional thresholding parametric estimators 
	has been established theoretically,  
	it is known that their real performance strongly depends on a choice of 
	tuning parameters; see, e.g., 
	Iacus and Yoshida \cite{iacus2017simulation}.}
{\coloro Examining each individual increment without other data
	is not always effective in 
	finding jumps.} It sometimes overlooks relatively small jumps 
{\tred 
	due to a conservative level of threshold to try to incorporate all Brownian increments.
	To resolve this problem,} 
Inatsugu and Yoshida {\colorb\cite{InatsuguYoshida2020Accepted}} introduced 
{\tred the so-called} global filters that 
examine all increments simultaneously and regard an increment of high rank in order of absolute size as a jump. Using the information about the size of other increments helps us detect jumps more accurately than {\tred the previous methods that ignore 
	such information. 
	Moreover, Inatsugu and Yoshida {\colorb\cite{InatsuguYoshida2020Accepted}} also removed 
	the assumption of low intensity of small jumps, {\colorg that} was used in 
	Shimizu and Yoshida \cite{ShimizuYoshida2006} 
	and Ogihara and Yoshida \cite{OgiharaYoshida2011}. 
	This is a theoretical advantage of the global jump filters, 
	in addition to their outperformance 
	in practice. 
}

In this paper, we {\tred will} apply the global filtering method to nonparametric volatility estimation. 
Specifically, we {\tred will} construct the ``global realized volatility (GRV) estimator" of 
the integrated volatility 
for the semimartingale $X$ having the decomposition (\ref{0301191155}). 
{\tred 
	Though $J$ and the jump part of $\sigma$ are assumed to be finitely active for each $n$, 
	we permit the number of jumps to diverge as $n$ tends to infinity.}
We will investigate the theoretical properties of GRV and then 
conduct numerical simulations to study their performance compared with 
{\tred traditional} methods, that is, 
{\tred the deterministic threshold estimator,} 
the bipower variation and the minimum 
realized volatility. 

The organization of this paper is as follows. 
{\tred Section \ref{Sec2} introduces} the GRV and its variant, the winsorized GRV (WGRV). 
In Section \ref{Sec3}, 
we introduce the local-global realized volatility (LGRV) and 
prove its convergence to {\tred the spot volatility}.
{\tred The LGRV will be used 
	for normalizing the increments to compute the global filter. 
	Section \ref{0211181814} gives} 
the rate of convergence of the GRV and WGRV in the 
situation where the intensity of jumps is high. In this case, 
we need a high and  
fixed cut-off rate $\alpha$ {\tred to eliminate harmful jumps}. 
In Section \ref{0211181815}, 
we allow the cut-off rate to vary according to the sample size. 
This ``moving threshold" method is for the situation where the intensity of jumps is 
moderate and small cut-off rate is applicable. 
{\tred {\colorg Section \ref{0211281401}} 
	briefly discusses the situation where true volatility is constant. 
	In this case, normalizing increments is not necessry,} so the estimator gets
a little simpler. 
{\tred Section \ref{Sec7} presents some}
simulation results to compare the {\tred real} performance
of the GRV, WGRV, bipower varition, and the mininum realized volatility. 

{\tred 
	Concluding, let us mention some technical aspects. 
	The global jump filter causes theoretical difficulty. 
	By nature, it uses all the data to classify {\coloro each increment} $\Delta_jX$. 
	This completely destroys the martingale structure in the model, which makes it difficult 
	to use orthogonality between the selected increments to validate 
	the law of large numbers and the central limit theorem. 
	However, it is possible to asymptotically recover the orthogonality 
	by {\colorg the} glocal and global filtering lemmas presented in 
	Sections \ref{Sec3} and \ref{0211181814}. 
	Technically, the argument here is closed within the semimartingale theory,  
	although the global filter breaks adaptivity of {\colorg the functionals}, in other words, 
	a quadratic variation with anticipative weights is treated. 
	On the other hand, Yoshida \cite{yoshida2020asymptotic} 
	suggests a use of the Malliavin calculus to {\colorg analyse} 
	robustified volatility estimators with anticipative weights.  
}

\section{Realized volatilities with a global jump filter}\label{Sec2}
The global jump filter introduced by Inatsugu and Yoshida {\colorb \cite{InatsuguYoshida2020Accepted}} uses 
the order statistics of the transformed increments of the observations. 
Suppose that an estimator $S_{n,j-1}$ of the spot volatility $\sigma(X_\tjm)^2$ 
(up to a common scaling factor) 
is given for each $j\in I_n=\{1,...,n\}$. 
Denote $\Delta_jU=U_\tj-U_\tjm$ for a process $U=(U_t)_{t\in[0,T]}$. 
Then the distribution of 
the scaled increment $S_{n,j-1}^{-1/2}\Delta_jX$ 
is expected to be well approximated by the standard normal distribution $N(0,1)$. 
Therefore, if the value 
\bea\label{0211281531}
V_j &=& \big|(S_{n,j-1})^{-1/2}\Delta_jX\big|
\eea
is relatively very large among $\calv_n=\{V_k\}_{k\in I_n}$, then 
plausibly we can infer that the $V_j$ involves jumps with high probability. 
The idea of the global jump filter is to eliminate the increment $\Delta_jX$ from the data 
if the corresponding $V_j$ is ranked {\colorb within} the top $100\alpha$\% in $\calv_n$.
More precisely, let 
\beas 
\calj_n(\alpha)
&=& \big\{j\in I_n;\>V_j<V_{(s_n(\alpha))}\big\}
\eeas
where 
\beas 
s_n(\alpha) &=& \lfloor n(1-\alpha) \rfloor
\eeas
for $\alpha\in[0,1)$, 
and we denote by $r_n(U_j)$ the rank of $U_j$ among the variables $\{U_i\}_{i\in I_n}$. 
Let 
\bea\label{0302110914} 
q(\alpha)
&=&
\int_{\{|z|\leq c(\alpha)^{1/2}\}}z^2\phi(z;0,1)dz
\eea
where $\phi(z;0,1)$ is the density function of $N(0,1)$ and 
$c(\alpha)$ defined by
\beas 
P\big[\zeta^2\leq c(\alpha)\big]&=&1-\alpha
\eeas
for $\zeta\sim N(0,1)$ and $\alpha\in[0,1)$.
Then the {\bf global realized volatility} (globally truncated realized volatility, GRV) 
with cut-off ratio $\alpha$ is defined by 
\bea\label{0211260834}
\bbV_n(\alpha)
&=& 
\sum_{j\in\calj_n(\alpha)}q(\alpha)^{-1}|\Delta_jX|^2K_{n,j}
\eea
where 
$K_{n,j}=1_{\{|\Delta_jX|\leq n^{-1/4}\}}$. 
As remarked in Inatsugu and Yoshida {\colorb \cite{InatsuguYoshida2020Accepted}}, the indicator function $K_{n,j}$ is set 
just for relaxing the conditions for validation. 
Generalization by using like $1_{\{|\Delta_jX|\leq B_1n^{-\delta_1}\}}$ 
with constants $B_1>0$ and $\delta_1\in(0,1/4]$ is straightforward, but 
we prefer simplicity in presentation of this article. 
In practice, the probability {\colorb that} $K_{n,j}$ executes the task is exponentially small 
by the large deviation principle. 
However, the moments of $\Delta J_t$ are not controllable without assumption, 
and we can simply avoid it by the cut-off function $K_{n,j}$.

{\colorb Winsorization} is a popular technique in robust statistics. 
In the present context, 
the {\bf Winsorized global realized volatility} (WGRV) is given by 
\beas 
\bbW_n(\alpha) 
&=& 
\sum_{j=1}^n{\sf w}(\alpha)^{-1}
\big\{|\Delta_jX|\wedge\big(S_{n,j-1}^{1/2}V_{(s_n(\alpha))}\big)\big\}^2
K_{n,j}
\eeas
where 
\beas 
{\sf w}(\alpha)
&=&
\int_\bbR\big(z^2\wedge c(\alpha)\big)\phi(z;0,1)dz.
\eeas

The cut-off ratio $\alpha\in[0,1)$ is a tuning parameter in estimation procedures. 
The bigger $\alpha$ provides the more stable estimates even in high intensity of jumps. 
On the other hand, the smaller $\alpha$ gives the more precise estimates 
if the intensity of jumps is low. 
Making trade-off between stability and precision is necessary in practice. 
As a matter of fact, these cases require different theoretical treatments. 
We will consider fixed $\alpha$ in Section \ref{0211181814}, 
and shrinking $\alpha$ in Section \ref{0211181815}.

\section{Local-global filter}\label{Sec3}
{\fred 
An estimator $S_{n,j-1}$ for the spot volatility (up to a constant scaling) is necessary to construct a global realized volatility. 
Naturally, we use the data around time $t$ to estimate $\sigma_t$. 
Since these data are also contaminated with jumps, we need a jump filter to construct 
a temporally-local estimator $S_{n,j-1}$. The idea of the global jump filter with the order statistics 
of the data around $t$ serves to eliminate the effects of jumps, 
not only theoretically but also practically as demonstrated by the simulation studies 
of Section \ref{Sec7}. 
In this section, we propose a local-global realized volatility and validate it 
by establishing in Section \ref{0301241722} the rate of convergence of the estimator. 
Since the local-global filter involves the order statistics, that destroy the martingale structure, 
we try to recover it by somewhat sophisticated lemmas given in Section \ref{0301241725}. 
The minimum realized volatility (minRV) made of the  temporally-local data is also a candidate of 
an estimator for the spot volatility. 
A rate of convergence of the local minRV is mentioned in Section \ref{0301241726}. 
}
\subsection{Glocal filtering lemmas}\label{0301241725}
For each $j\in I_n$, let 
\beas 
\underline{j}_n
&=& 
\left\{\begin{array}{cl}
1&(j\leq \kappa_n)\y
j-\kappa_n&(\kappa_n+1\leq j\leq n-\kappa_n)\y
n-2\kappa_n&(j\geq n-\kappa_n+1)
\end{array}\right.
\eeas
for $\kappa_n\in\bbZ_+$ satisfying $2\kappa_n+1\leq n$. 
Let 
$I_{n,j}=\{\underline{j}_n,\underline{j}_n+1,...,\underline{j}_n+2\kappa_n\}$. 
Let 
\beas 
\hU_{j,k}\yeq h^{-1/2}
\sigma_{t_{\underline{j}_n-1}}^{-1}
\Delta_kX\quad \text{and}\quad
W_j\yeq h^{-1/2}\Delta_jw
\eeas
for $j,k\in I_n$. 
Both variables $\hU_{j,k}$ and $W_j$ depend on $n$. 
Let 
\beas 
\hR_{j,k}
&=& 
\hU_{j,k}-W_k-h^{-1/2}\sigma_{t_{\underline{j}_n-1}}^{-1}\Delta_kJ
\eeas
for $j,k\in I_n$. 
Denote $L^\inftym=\cap_{p>1}L^p$. 

{\sred
Let $N=\sum_{s\in(0,\cdot]}1_{\{\Delta J_s\not=0\}}$. 
Let $\wt{\sigma}=\sigma-J^\sigma$ for $J^\sigma=\sum_{s\in(0,\cdot]}\Delta\sigma_s$, 
and let $N^\sigma=\sum_{s\in(0,\cdot]}1_{\{\Delta J^\sigma_s\not=0\}}$. 
{\colorb We assume that 
$N^\sigma_T<\infty$ a.s. }
Moreover, let $\ol{N}=N+N^\sigma$. 
Let $\wt{X}=X-J$. 
A counting process will be identified with a random measure. 
Let ${\tt I}_{n,j}=\big(t_{\underline{j}_n-1},t_{\underline{j}_n+2\kappa_n}\big]$. 
}
{\sred 
\bd
\im[{\bf [G1]}] 
{\bf (i)} For every $p>1$, $\sup_{t\in[0,T]}\|\sigma_t\|_p<\infty$ and 
\beas 
\big\|\wt{\sigma}_t-\wt{\sigma}_s\big\|_p &\leq& C(p)|t-s|^{1/2}
\quad(t,s\in[0,T])
\eeas
for some constant $C(p)$ for every $p>1$. 
\im[{\bf (ii)}] $\sup_{t\in[0,T]}\|b_t\|_p<\infty$ for every $p>1$. 
\im[{\bf (iii)}] $\sigma_t\not=0$ a.s. for every $t\in[0,T]$, an 
$\sup_{t\in[0,T]}\big\|\sigma_t^{-1}\big\|_p<\infty$ for every $p>1$. 
%
\ed
}
\begin{en-text}
\begin{lemma}\label{0211190303}
Under $[G1]$, 
\beas 
\sup_{j\in I_n}\sup_{k\in I_{n,j}}
\bigg\|\hR_{j,k}1_{\big\{\Delta_kN=0,\>\sum_{k'\in I_{n,j}}\Delta_{k'}N^\sigma=0\big\}}
\bigg\|_p
&=& 
O\left(\left(\frac{\kappa_n}{n}\right)^{1/2}\right)
\eeas
as $n\to\infty$ for every $p>1$. 
\end{lemma}
\proof 
For $j\in I_n$, 
on the event $\big\{\>\sum_{k'\in I_{n,j}}\Delta_{k'}N^\sigma=0\big\}$, 
\beas 
\hR_{j,k}1_{\{\Delta_kN=0\}}
&=& 
\big(h^{-1/2}\sigma_\tjm^{-1}\Delta_k\wt{X}-h^{-1/2}\Delta_jw \big)1_{\{\Delta_kN=0\}}
\nn\\&=&
\big(h^{-1/2}\sigma_\tjm^{-1}\Delta_k\wt{X}-h^{-1/2}\Delta_jw \big)
-\big(h^{-1/2}\sigma_\tjm^{-1}\Delta_k\wt{X}-h^{-1/2}\Delta_jw \big)1_{\{\Delta_kN>0\}}
\eeas
Use orthogonality in 
$\sum_{k\in I_{n,j}}\int_{t_{k-1}}^{t_k}(\wt{\sigma}_s-\wt{\sigma}_{t_{k-1}})dw_s$. 
\qed
\begin{lemma}\label{0211190303}
Under $[G1]$, 
\bea\label{0211191046}
\sup_{j\in I_n}\sup_{k\in I_{n,j}}
\bigg\|\hR_{j,k}1_{\big\{\sum_{k'\in I_{n,j}}\Delta_{k'}N^\sigma=0\big\}}
\bigg\|_p
&=& 
O\left(\left(\frac{\kappa_n}{n}\right)^{1/2}\right)
\eea
as $n\to\infty$ for every $p>1$. 
\end{lemma}
\proof 
For $j\in I_n$, let 
$E(j)=\big\{\>\sum_{k'\in I_{n,j}}\Delta_{k'}N^\sigma=0\big\}$. 
Then, for $k\in I_{n,j}$, 
\beas 
\hR_{j,k}1_{E(j)}
&=& 
\big(h^{-1/2}\sigma_{t_{\underline{j}_n-1}}^{-1}\Delta_k\wt{X}-h^{-1/2}\Delta_kw \big)1_{E(j)}
\nn\\&=&
R^{(\ref{0211191011})}+R^{(\ref{0211191012})}
\eeas
where 
\bea\label{0211191011}
R^{(\ref{0211191011})}
&=& 
h^{-1/2}\int_\tkm^\tk
\sigma_{t_{\underline{j}_n-1}}^{-1}(\wt{\sigma}_t-\wt{\sigma}_{t_{\underline{j}_n-1}})dw_t
+
h^{-1/2}\sigma_{t_{\underline{j}_n-1}}^{-1}\int_\tkm^\tk b_tdt
\eea
and 
\bea\label{0211191012}
R^{(\ref{0211191012})}
&=&
-\big(h^{-1/2}\sigma_{t_{\underline{j}_n-1}}^{-1}\Delta_k\wt{X}-h^{-1/2}\Delta_jw \big)1_{E(j)^c}. 
\eea
Apply the Burkholder-Davis-Gundy inequality to the martingale part of (\ref{0211191011}), 
and the H\"older inequality gives 
\beas 
\big\|R^{(\ref{0211191012})}\big\|_p
&\simleq& 
\eeas  
Then we obtain the estimate (\ref{0211191046}). 
\qed
\end{en-text}
\begin{lemma}\label{0211190303}
Under $[G1]$, 
\bea\label{0211191046}
\sup_{j\in I_n}\sup_{k\in I_{n,j}}
\big\|\hR_{j,k}1_{\{N^\sigma({\tt I}_{n,j})=0\}}\big\|_p
&=& 
O\left(\left(\frac{\kappa_n}{n}\right)^{1/2}\right)
\eea
as $n\to\infty$ for every $p>1$. 
\end{lemma}
\proof 
For $j\in I_n$, let 
$E(j)=\{N^\sigma({\tt I}_{n,j})=0\}$. 
Then, for $k\in I_{n,j}$, 
\bea\label{0211191011} 
\hR_{j,k}1_{E(j)}
&=& 
\big(h^{-1/2}\sigma_{t_{\underline{j}_n-1}}^{-1}\Delta_k\wt{X}-h^{-1/2}\Delta_kw \big)1_{E(j)}
\nn\\&=&
h^{-1/2}\int_\tkm^\tk
\sigma_{t_{\underline{j}_n-1}}^{-1}(\wt{\sigma}_t-\wt{\sigma}_{t_{\underline{j}_n-1}})dw_t1_{E(j)}
\nn\\&&
+
h^{-1/2}\sigma_{t_{\underline{j}_n-1}}^{-1}\int_\tkm^\tk b_tdt1_{E(j)}
\eea
We obtain (\ref{0211191046}) 
by applying the Burkholder-Davis-Gundy inequality to the martingale part of (\ref{0211191011})
after the trivial estimate $1_{E(j)}\leq1$. 
\qed\halflineskip

For $j\in I_n$, denote by $r_{n,j}(U_k)$ the rank of the element $U_k$ 
among a collection of random variables $\{U_\ell\}_{k\in I_{n,{\sred j}}}$. 
Let 
\beas &&
0<\eta_2<\eta_1,\qquad
\ol{\kappa}_n \yeq 2\kappa_n+1,\quad 
\nn\\&&
{\sf a}_n\yeq\lfloor (1-\alpha_0)\ol{\kappa}_n-\ol{\kappa}_n^{\>1-\eta_2}\rfloor,\qquad
\wh{\sf a}_n\yeq\lfloor{\sf a}_n-\ol{\kappa}_n^{\>1-\eta_2}\rfloor
\eeas
for $\alpha_0\in[0,1)$. 
Let 
\bea\label{0211201930}
L_{n,j,k}
&=& 
\big\{r_{n,j}(|W_k|)\leq {\sf a}_n-\ol{\kappa}_n^{\>1-\eta_2}\big\}
\cap
\big\{|W|_{(j,{\sf a}_n)}-|W_k|<\ol{\kappa}_n^{\>-\eta_1}\big\}
\eea
where $\big(|W|_{(j,k)}\big)_{k\in I_{n,j}}$ are the ordered statistics made from $\{|W_k|\}_{k\in I_{n,j}}$. 
In the same way as Lemma 1 of Inatsugu and Yoshida {\colorb\cite{InatsuguYoshida2020Accepted}}, 
we obtain the following result. 
\begin{lemma}\label{0211191227}
Let $\alpha_0\in(0,1)$. 
Suppose that ${\sred \eta_1}<1/2$ 
and that 
$n^{-\ep}\kappa_n\to\infty$ as $n\to\infty$ 
for some $\ep\in(0,1)$. 
Then 
\beas 
\sup_{j\in I_n}P\bigg[\bigcup_{k\in I_{n,j}}L_{n,j,k}\bigg]\yeq O(n^{-L})
\eeas
as $n\to\infty$ for every $L>0$. 
\end{lemma}
\halflineskip

Define $\calk_{n,j}(\alpha_0)$ by 
\beas 
\calk_{n,j}(\alpha_0)
&=&
\big\{k\in I_{n,j};\>r_{n,j}
({\colorb |\Delta_kX|}) 
\leq (1-\alpha_0)\ol{\kappa}_n\big\},
\eeas
where 
{\colorb 
$r_{n,j}(|\Delta_kX|)$
}
is the rank of ${\colorb |\Delta_kX|}$ 
among 
$\{ {\colorb |\Delta_{k'}X| }\}_{k'\in I_{n,j}}$. 
Let 
\beas 
\wh{\calk}_{n,j}(\alpha_0)
&=&
\big\{k\in I_{n,j};\>r_{n,j}(|W_k|)\yleq\hat{\sf a}_n\big\}{\colorb .}
\eeas
Let 
\beas 
\Omega_{n,j}
&=& 
\bigcap_{k\in I_{n,j}}\bigg[
\bigg\{|\hR_{j,k}|1_{\{N^\sigma({\tt I}_{n,j})=0\}}<2^{-1}\ol{\kappa}_n^{\>-\eta_1}\bigg\}
\cap L_{n,j,k}^c\bigg]{\colorb .}
\eeas
%
Let {\sred
\bea\label{0301021550}
\call_n=\big\{j\in I_n;\>\ol{N}({\tt I}_{n,j})\not=0\big\}.
\eea
}
\begin{lemma}\label{0211200634}
\bd
\im[{\bf (a)}] 
$\wh{\calk}_{n,j}(\alpha_0) \subset \calk_{n,j}(\alpha_0)$ 
on $\Omega_{n,j}$ if $j\in\call_n^c$. 

\im[{\bf (b)}]
$\ds
1_{\Omega_{n,j}}1_{\{j\in\call_n^c\}}\>
\#\big(\calk_{n,j}(\alpha_0)\setminus\wh{\calk}_{n,j}(\alpha_0)\big)
\leq
{\sred4}
\>\ol{\kappa}_n^{\>1-\eta_2}
\quad(j\in I_{n.j},\>n\in\bbN).
$
\ed
\end{lemma}
\proof 
Let $n\in\bbN$ and suppose that ${\colorr j\in\call_n^c}$. 
We will work on $\Omega_{n,j}$. 
For a pair 
{\colorr$(k_1,k_2)\in I_{n,j}^2$}, 
suppose that 
\bea\label{0211201926}
r_{n,j}(|W_{k_1}|)\leq\wh{\sf a}_n\quad\text{and}\quad 
r_{n,j}(|W_{k_2}|)\geq{\sf a}_n.
\eea
Then 
$
|\hU_{j,k_1}|
<
|W_{k_1}|+2^{-1}\ol{\kappa}_n^{\>-\eta_1}
$, since 
$\Delta_{k_1}N=0$ and 
{\colorr$N^\sigma({\tt I}_{n,j})=0$} when 
{\colorr $j\in\call_n^c$}, 
and then 
$|\hR_{j,{k_1}}|<2^{-1}\ol{\kappa}_n^{\>-\eta_1}$ 
on $\Omega_{n,j}$. 
By the first inequality of (\ref{0211201926}), 
$r_{n,j}(|W_{k_1}|)\leq{\sf a}_n-\ol{\kappa}_n^{\>1-\eta_2}$, 
and hence on $\Omega_{n,j}\subset L_{n,j,{\sred k_1}}^c$, we have 
$|W|_{(j,{\sf a}_n)}-|W_{k_1}|\geq\ol{\kappa}_n^{\>-\eta_1}$ 
by the definition (\ref{0211201930}) of $L_{n,j,k}$. 
Therefore 
\bea\label{0211201012}
|\hU_{j,k_1}| &<& |W|_{(j,{\sf a}_n)}-2^{-1}\ol{\kappa}_n^{\>{\sred -\eta_1}}.
\eea
The assumption 
{\colorr$j\in\call_n^c$} 
entails 
$|\hR_{j,{k_2}}|<2^{-1}\ol{\kappa}_n^{\>-\eta_1}$ 
on $\Omega_{n,j}$, and hence 
$|W_{k_2}|-2^{-1}\ol{\kappa}_n^{\>-\eta_1}<|\hU_{j,k_2}|$ 
due to $\Delta_{k_2}J=0$. 
From (\ref{0211201012}), we have got 
\bea\label{0211202027}
|\hU_{j,k_1}|<|\hU_{j,k_2}|
\eea
on $\Omega_{n,j}$ 
if 
$j\in\call_n^c$ and if a pair $(k_1,k_2)\in I_{n,j}^2$ satisfies (\ref{0211201926}). 

We are working on $\Omega_{n,j}$ yet. 
Suppose that $j\in\call_n^c$ and $k_1\in\wh{\calk}_{n,j}(\alpha_0)$. 
Then the inequality (\ref{0211202027}) holds 
for any $k_2\in I_{n,j}$ satisfying $r_{n,j}(|W_{k_2}|)\geq{\sf a}_n$. 
So, there are at least 
$\lfloor\alpha_0\ol{\kappa}_n+1\rfloor\big(\leq
\alpha_0\ol{\kappa}_n{\colorg+\ol{\kappa}_n^{\>1-\eta_2}}+1\leq
\ol{\kappa}_n-{\sf a}_n+1\big)$ variables $\hU_{j,k_2}$ that satisfy (\ref{0211202027}). 
Then $r_{n,j}(|\hU_{j,k_1}|)\leq (1-{\sred\alpha_0})\ol{\kappa}_n$, and hence 
$k_1\in\calk_{n,j}(\alpha_0)$. 
Thus, we found 
\beas 
\wh{\calk}_{n,j}(\alpha_0) &\subset& \calk_{n,j}(\alpha_0)
\eeas
on $\Omega_{n,j}$ if $j\in\call_n^c$, that is, (a). 

We still work on $\Omega_{n,j}$. 
Suppose that $j\in\call_n^c$ and 
$k_2\in\calk_{n,j}(\alpha_0)\setminus\wh{\calk}_{n,j}(\alpha_0)$. 
When $r_{n,j}(|W_{k_2}|)<{\sf a}_n$, 
since $r_{n,j}(|W_{k_2}|)>\wh{\sf a}_n$ due to 
$k_2\in{\sred\wh{\calk}}_{n,j}(\alpha_0)^c$, 
we see 
\bea\label{0211211229} 
1_{\{j\in\call_n^c\}}\>
\#\big\{k_2\in \calk_{n,j}(\alpha_0)\setminus\wh{\calk}_{n,j}(\alpha_0);\>
r_{n,j}(|W_{k_2}|)<{\sf a}_n\big\}
&\leq&
\ol{\kappa}_n^{\>1-\eta_2}
\eea
on $\Omega_{n,j}$.
When $r_{n,j}(|W_{k_2}|)\geq{\sf a}_n$, for 
any $k_1$ satisfying $r_{n,j}(|W_{k_1}|)\leq\wh{\sf a}_n$, we have (\ref{0211202027}). 
Therefore 
\beas 
\#\big\{k_1\in I_{n,j};\>|\hU_{j,k_1}|<|\hU_{j,k_2}|\big\}
&\geq&
1_{\{j\in\call_n^c\}}1_{\{r_{n,j}(|W_{k_2}|)\geq{\sf a}_n\}}\>{\sred\wh{\sf a}}_n,
\eeas
in other words, 
\bea\label{0211211225}
r_{n,j}(|\hU_{j,k_2}|) &>& {\sred\wh{\sf a}}_n
\eea
on $\Omega_{n,j}$ if $j\in\call_n^c$ and $r_{n,j}(|W_{k_2}|)\geq{\sf a}_n$. 
Moreover, 
$r_{n,j}(|\hU_{j,k_2}|)\leq\lfloor(1-\alpha_0)\ol{\kappa}_n\rfloor$ 
since $k_2\in\calk_{n,j}(\alpha_0)$.
Combining this estimate with (\ref{0211211225}), we obtain 
\bea\label{0211211232}
1_{\{j\in\call_n^c\}}\>
\#\big\{
k_2\in\calk_{n,j}(\alpha_0)\setminus\wh{\calk}_{n,j}(\alpha_0);\>
r_{n,j}(|W_{k_2}|)\geq{\sf a}_n\big\}
&\leq&
(1-\alpha_0)\ol{\kappa}_n-{\sred\wh{\sf a}}_n
\nn\\&\leq&
{\sred2\>}\ol{\kappa}_n^{\>1-\eta_2}+1
\eea
From (\ref{0211211229}) and (\ref{0211211232}), we obtain (b). 
\qed\halflineskip

For $\eta_3\in\bbR$, $j\in I_n$ and a sequence of random variables $(V_j)_{j\in I_n}$, let  
\beas 
\cald_{n,j}
&=& 
{\sred\ol{\kappa}_n^{\eta_3}}\bigg|\frac{1}{{\sred\ol{\kappa}_n}}\sum_{k\in\calk_{n,j}(\alpha_0)}V_k
-
\frac{1}{{\sred\ol{\kappa}_n}}\sum_{k\in\wh{\calk}_{n,j}(\alpha_0)}V_k
\bigg|
\eeas
The following lemma follows from Lemma \ref{0211200634} immediately. 
\begin{lemma}\label{0211231513}
{\bf (i)} Let $p\geq1$. Then 
\beas 
\big\|\cald_{n,j}\big\|_p
&\leq&
{\sred4\>}
{\sred\ol{\kappa}_n^{\eta_3-\eta_2}}\bigg\|\max_{k\in I_{n,j}}|V_k|
{\colorr 1_{\Omega_{n,j}\cap\{j\in\call_n^c\}}}
\bigg\|_p
+{\sred\ol{\kappa}_n^{\eta_3}}\bigg\|\max_{k\in I_{n,j}}|V_k|1_{\Omega_{n,j}^c}\bigg\|_p
+{\sred\ol{\kappa}_n^{\eta_3}}\bigg\|\max_{k\in I_{n,j}}|V_k|1_{\{j\in\call_n\}}\bigg\|_p
\eeas
for $j\in I_n$, $n\in\bbN$. 
\bd\im[{\bf (ii)}]
Let $p\geq1$ and $\eta_4>0$. Then 
\beas 
\big\|\cald_{n,j}\big\|_p
&\leq&
{\sred4\>}
{\sred\ol{\kappa}_n^{\eta_3-\eta_2}}
\bigg(\kappa_n^{\eta_4}
+{\sred\ol{\kappa}_n}\max_{k\in I_{n,j}}\bigg\||V_k|1_{\{|V_k|>\kappa_n^{\eta_4}\}}
{\colorr 1_{\Omega_{n,j}\cap\{j\in\call_n^c\}}}
\bigg\|_p
\bigg)
\nn\\&&
+{\sred\ol{\kappa}_n^{\eta_3}}\bigg\|\max_{k\in I_{n,j}}|V_k|1_{\Omega_{n,j}^c}\bigg\|_p
+{\sred\ol{\kappa}_n^{\eta_3}}\bigg\|\max_{k\in I_{n,j}}|V_k|1_{\{j\in\call_n\}}\bigg\|_p
\eeas
for $j\in I_n$, $n\in\bbN$. 
\ed
\end{lemma}
\halflineskip

Let 
\beas
\wt{\calk}_{n,j}(\alpha_0)
&=& 
\big\{k\in I_{n,j};\>|W_k|\leq c(\alpha_0)^{1/2}\big\}
\eeas
For $\eta_3>0$, $j\in I_n$ and a sequence of random variables $(V_j)_{j\in I_n}$, let  
\beas 
\wt{\cald}_{n,j}
&=& 
{\sred\ol{\kappa}_n^{\eta_3}}\bigg|\frac{1}{{\sred\ol{\kappa}_n}}\sum_{k\in\wh{\calk}_{n,j}{\sred(\alpha_0)}}V_k
-
\frac{1}{{\sred\ol{\kappa}_n}}\sum_{k\in\wt{\calk}_{n,j}{\sred(\alpha_0)}}V_k
\bigg|
\eeas
Let 
\bea\label{0211241445}
\wt{\Omega}_{n,j}=\big\{\big||W|_{(j,\wh{\sf a}_n)}-c(\alpha_0)^{1/2}\big|
<\check{C}\kappa_n^{-\eta_2}\big\}
\eea
for $j\in I_n$, where $\check{C}$ is a positive constant. 
\begin{lemma}\label{0211231405}
Let $\eta_3\in\bbR$. 
Then
\bd
\im[(i)] For $p\geq1$ and $j\in I_n$, 
\beas 
\big\|\wt{\cald}_{n,j}\big\|_p
&\leq&
{\sred\ol{\kappa}_n^{\eta_3}}\bigg\|\max_{k'\in I_{n,j}}|V_{k'}|\>
\frac{1}{{\sred\ol{\kappa}_n}}\sum_{k\in I_{n,j}}
1_{\big\{\big||{\sred W_k}|-c(\alpha_0)^{1/2}\big|{\sred\><\>}\check{C}\kappa_n^{-\eta_2}\big\}}\bigg\|_p
+
{\sred\ol{\kappa}_n^{\eta_3}}\bigg\|1_{\wt{\Omega}_{n,j}^c}\max_{k'\in I_{n,j}}|V_{k'}|\bigg\|_p
\eeas

\im[(ii)] For $p_1>p\geq1$ and  $j\in I_n$, 
\beas 
\big\|\wt{\cald}_{n,j}\big\|_p
&\leq&
{\sred\ol{\kappa}_n^{\eta_3}}\bigg\|\max_{k\in I_{n,j}}|V_k|\bigg\|_p
P\bigg[\big||{\sred W_1}|-c(\alpha_0)^{1/2}\big|{\sred\><\>}\check{C}\kappa_n^{-\eta_2}\bigg]
\nn\\&&
+
{\sred\ol{\kappa}_n^{\eta_3}}\bigg\|\max_{k\in I_{n,j}}|V_k|\bigg\|_{pp_1(p_1-p)^{-1}}
\bigg\|\frac{1}{{\sred\ol{\kappa}_n}}\sum_{k\in I_{n,j}}\bigg(
1_{\big\{\big|{\colorb |W_k|}-c(\alpha_0)^{1/2}\big|{\sred\><\>}	\check{C}\kappa_n^{-\eta_2}\big\}}
\nn\\&&\hspace{180pt}
-P\bigg[\big||W_k|-c(\alpha_0)^{1/2}\big|{\sred\><\>}\check{C}\kappa_n^{-\eta_2}\bigg]
\bigg)\bigg\|_{p_1}
\nn\\&&
+
{\sred\ol{\kappa}_n^{\eta_3}}P\big[\wt{\Omega}_{n,j}^c\big]^{1/p_1}
\bigg\|\max_{k\in I_{n,j}}|V_k|\bigg\|_{pp_1(p_1-p)^{-1}}
\eeas
\ed
\end{lemma}
\proof
For $k\in I_{n,j}$, 
\beas&&
\wt{\Omega}_{n,j}\cap\big\{\>r_{n,j}(|W_k|)\yleq\wh{\sf a}_n\big\}^c\cap\big\{|W_k|\leq c(\alpha_0)^{1/2}\big\}
\nn\\&=&
\big\{\big||W|_{(j,\wh{\sf a}_n)}-c(\alpha_0)^{1/2}\big|
<\check{C}\kappa_n^{-\eta_2}\big\}\cap\big\{|W_k|>|W|_{(j,\wh{\sf a}_n)}\big\}\cap\big\{|W_k|\leq c(\alpha_0)^{1/2}\big\}
\nn\\&\subset&
\big\{\big||W_k|-c(\alpha_0)^{1/2}\big|<\check{C}\kappa_n^{-\eta_2}\big\}
\eeas
and 
\beas&&
\wt{\Omega}_{n,j}\cap\big\{\>r_{n,j}(|W_k|)\yleq\wh{\sf a}_n\big\}\cap\big\{|W_k|\leq c(\alpha_0)^{1/2}\big\}^c
\nn\\&=&
\big\{\big||W|_{(j,\wh{\sf a}_n)}-c(\alpha_0)^{1/2}\big|
<\check{C}\kappa_n^{-\eta_2}\big\}\cap\big\{|W_k|\leq|W|_{(j,\wh{\sf a}_n)}\big\}\cap\big\{|W_k|> c(\alpha_0)^{1/2}\big\}
\nn\\&\subset&
\big\{\big||W_k|-c(\alpha_0)^{1/2}\big|<\check{C}\kappa_n^{-\eta_2}\big\}. 
\eeas
Thus we obtain (i). 
Property (ii) follows from (i). 
\qed\halflineskip

\begin{lemma}\label{0211241300}
If the constant $\check{C}$ in (\ref{0211280310}) is sufficiently large, then 
\beas
\sup_{j\in I_n}P\big[\wt{\Omega}_{n,j}^c\big]&=&O(n^{-L})
\eeas
as $n\to\infty$ for any $L>0$. 
\end{lemma}
\proof
We have 
\bea\label{0211241430}
&&
P\big[|W|_{(j,\wh{\sf a}_n)}-c(\alpha_0)^{1/2}<-\check{C}\kappa_n^{-\eta_2}\big]
\nn\\&\leq&
P\bigg[|W|_{\big(j,\lfloor{\sf a}_n-\ol{\kappa}_n^{\>1-\eta_2}-1\rfloor\big)
}<c(\alpha_0)^{1/2}-\check{C}\kappa_n^{-\eta_2}\bigg]
\nn\\&\leq&
P\bigg[\sum_{k\in I_{n,j}}1_{A_{{\sred n,k}}}\geq \lfloor{\sf a}_n-\ol{\kappa}_n^{\>1-\eta_2}-1\rfloor\bigg]
\nn\\&=&
P\bigg[\ol{\kappa}_n^{-1/2}\sum_{k\in I_{n,j}}\big\{1_{A_{n,k}}-P[A_{n,k}]\big\}
\geq {\colorb C_{n}}
\bigg]
\eea
where
\beas 
A_{n,k}
&=& 
\big\{|W_k|< c(\alpha_0)^{1/2}-\check{C}\kappa_n^{-\eta_2}\big\},
\nn\\
{\colorb C_{n}}
&=&
\ol{\kappa}_n^{\>-1/2}\big({\sf a}_n-\ol{\kappa}_n^{\>1-\eta_2}-2-\ol{\kappa}_n
P[A_{n,1}]\big).
\eeas
{\sred By using the mean-value theorem, we obtain}
\beas 
{\colorb C_{n}}
&\sim&
\ol{\kappa}_n^{\>-1/2}\bigg[
(1-\alpha_0)\ol{\kappa}_n-2\ol{\kappa}_n^{\>1-\eta_2}
-\ol{\kappa}_n\big\{1-\alpha_0-
{\colorb 2 \phi\big(c(\alpha_0)^{1/2};0,1\big)\>}
\check{C}\kappa_n^{-\eta_2}\big\}\bigg]
\nn\\&\gtrsim&
\ol{\kappa}_n^{\>{\colorb \frac{1}{2}}-\eta_2}
\eeas
as $n\to\infty$ if we choose a sufficiently large $\check{C}$. 
Therefore, the $L^p$-boundedness of the random variables in (\ref{0211241430}) gives 
\bea\label{0211241424}
\sup_{j\in I_n}
P\big[|W|_{(j,\wh{\sf a}_n)}-c(\alpha_0)^{1/2}<-\check{C}\kappa_n^{-\eta_2}\big]
&=&
O(n^{-L})
\eea
as $n\to\infty$ for any $L>0$. 
In a similar way, we know 
\bea\label{0211241441}
P\big[|W|_{(j,\wh{\sf a}_n)}-c(\alpha_0)^{1/2}>\check{C}\kappa_n^{-\eta_2}\big]
&=& 
O(n^{-L})
\eea
as $n\to\infty$ for any $L>0$. 
Then we obtain the result from (\ref{0211241424}) and (\ref{0211241441}). 
\qed\halflineskip
\begin{en-text}
\beas 
\wt{\Omega}_{n,j}^c
&=& 
\big\{|W|_{(j,\wh{\sf a}_n)}-c(\alpha_0)^{1/2}>\check{C}\kappa_n^{-\eta_2}\big\}
\cup
\big\{|W|_{(j,\wh{\sf a}_n)}-c(\alpha_0)^{1/2}<-\check{C}\kappa_n^{-\eta_2}\big\}
\eeas
\end{en-text}

\subsection{{\sred Local-global realized volatility}}\label{0301241722}
We introduce the {local-global realized volatility} (LGRV) 
\bea\label{0211260215} 
\bbL_{n,j}(\alpha_0)
&=& 
\frac{n}{\ol{\kappa}_nT}
\sum_{k\in\calk_{n,j}(\alpha_0)}q(\alpha_0)^{-1}|\Delta_kX|^2K_{n,k}{\colorb .}
\eea
\begin{theorem}\label{0301021509}
Suppose that $[G1]$ is fulfilled.
\begin{en-text}
and that 
\beas 
&c_0\in(0,1), \quad \eta_1\in\bigg(0,\half\big(\frac{1}{c_0}-1\big)\bigg), \quad 
\eta_2\in(0,1/2),&
\nn\\
&\eta_3>0,\quad\eta_4>0,\quad \eta_3+\eta_4<\eta_2,
\quad B>0.&
\eeas
\end{en-text}
For $c_0\in(0,1)$ and $B>0$, 
suppose that $\kappa_n\sim Bn^{c_0}$ as $n\to\infty$. 
Then 
\bea\label{0211221359} 
\sup_{n\in\bbN}
\sup_{j\in I_n}\sup_{k\in I_{n,j}}n^{\>\gamma_*}
\big\|1_{\{j\in\call_n^c\}}\>\big(\bbL_{n,j}(\alpha_0)-\sigma_\tk^2\big)\big\|_p
&<&
\infty
\eea
as $n\to\infty$ for any constant $\gamma_*$ satisfying 
\beas 
\gamma_* &<& \min\bigg\{\half(1-c_0),\half c_0\bigg\}. 
\eeas
\end{theorem}
\proof
(I) 
We have $\kappa_n\sim n^{c_0}\sim h^{-c_0}$ and 
$n/\ol{\kappa}_n\sim n^{1-c_0}\sim h^{c_0-1}$. 
{\sred Let 
\beas 
\cald_{n,j}^*
&=&
\ol{\kappa}_n^{\eta_3}
\bigg\{
\frac{n}{\ol{\kappa}_n}
\sum_{k\in\calk_{n,j}(\alpha_0)}|\Delta_kX|^2K_{n,k}
-
\frac{n}{\ol{\kappa}_n}
\sum_{k\in\wh{\calk}_{n,j}(\alpha_0)}|\Delta_kX|^2K_{n,k}\bigg\}
\eeas
}
{\sred 
Applied to $V_k=n|\Delta_kX|^2K_{n,k}1_{\{j\in\call_n^c\}}$, 
Lemma \ref{0211231513} (ii)
}
gives 
\bea\label{0211241150}
\big\|\cald_{n,j}^*1_{\{j\in\call_n^c\}}\big\|_p
&\leq&
\Phi^{(\ref{0211240511})}_{n,j}+\Phi^{(\ref{0211240512})}_{n,j}
\eea
{\sred for every $p>1$,}
where
\bea\label{0211240511}
\Phi^{(\ref{0211240511})}_{n,j}
&=& 
{\sred4}{\sred\ol{\kappa}_n^{\eta_3-\eta_2}}
\bigg(\kappa_n^{\eta_4}
+{\sred\ol{\kappa}_n}\max_{k\in I_{n,j}}\bigg\|n|\Delta_kX|^2
1_{\{n|\Delta_kX|^2
>\kappa_n^{\eta_4}\}}
{\colorr 1_{\{j\in\call_n^c\}}}
\bigg\|_p
\bigg)
\eea
and 
\bea\label{0211240512}
\Phi^{(\ref{0211240512})}_{n,j}
&=&
{\sred\ol{\kappa}_n^{\eta_3}}\bigg\|\max_{k\in I_{n,j}}n|\Delta_kX|^2K_{n,k}1_{\Omega_{n,j}^c}\bigg\|_p{\sred.}
\eea
\begin{en-text}
\bea\label{0211240513}
\Phi^{(\ref{0211240513})}_{n,j}
&=&
\kappa_n^{\eta_3}\bigg\|\max_{k\in I_{n,j}}n|\Delta_kX|^2K_{n,k}|1_{\{j\in\call_n\}}\bigg\|_p
\eea
for $j\in I_n$, $n\in\bbN$. 
\end{en-text}
%
Since there is no jump of $J$ on $\{j\in\call_n^c\}$, we see 
\bea\label{0211240530}
\sup_{j\in I_n}\sup_{k\in I_{n,j}}\big\|n|\Delta_kX|^21_{\{j\in\call_n^c\}}\big\|_p=O(1)
\eea
for every $p>1$, 
as a result, the $L^p$-norm on the right-hand side of (\ref{0211240511}) 
is of $O(n^{-L})$ for arbitrary $L>0$, and hence 
\bea\label{0211240541}
\Phi^{(\ref{0211240511})}_{n,j}
&=& 
O\big(\kappa_n^{\eta_3-\eta_2+\eta_4}\big)
\eea
as $n\to\infty$. 
Similarly to (\ref{0211240530}), we obtain 
\bea\label{0211240531}
\sup_{j\in I_n}P\big[\Omega_{n,j}^{{\sred c}}\big] &=& O(n^{-L})
\eea
as $n\to\infty$ for every $L>0$, from 
{\sred Lemma \ref{0211191227} as well as} 
Lemma \ref{0211190303} 
because $(n/\kappa_n)^{1/2}\kappa_n^{-\eta_1}\gg1$ when 
$2^{-1}(c_0^{-1}-1)>\eta_1$. 
{\sred Then}
\bea\label{0211240533}
\Phi^{(\ref{0211240512})}_{n,j}
&\leq&
\kappa_n^{\eta_3}n^{1/2}
P\big[\Omega_{n,j}^c\big]^{1/p}
\yeq  O(n^{-L})
\eea
for every $L>0$ and $p>1$. 
\begin{en-text}
\bea\label{0211240536}\koko
\Phi^{(\ref{0211240513})}_{n,j}
&=&
\kappa_n^{\eta_3}\bigg\|\max_{k\in I_{n,j}}n|\Delta_kX|^2K_{n,k}|1_{\{j\in\call_n\}}\bigg\|_p
\eea
\end{en-text}
\begin{en-text}
On the event $\{j\in\call_n^c\}\cap\Omega_{n,j}$, we have 
\bea\label{0211230448} 
\sup_{j\in I_n}\big\|1_{\{j\in\call_n^c\}}1_{\Omega_{n,j}}
\cald_{n,j}^*\big\|_p
&\simleq&
\kappa_n^{\eta_3}\times
\frac{n}{\kappa_n}\times\kappa_N^{1-\eta_2}\times h
\>\sim\>
h^{c_0(\eta_2-\eta_3)}
\eea
by Lemma \ref{0211200634}. 
\end{en-text}
From (\ref{0211241150}), (\ref{0211240541}) and (\ref{0211240533}), 
\bea\label{0211241151}
\big\|\cald_{n,j}^*1_{\{j\in\call_n^c\}}\big\|_p
&=&
O\big(\kappa_n^{\eta_3-\eta_2+\eta_4}\big)
\yeq 
O\big(n^{-c_0(\eta_2-\eta_3-\eta_4)}\big)
\eea
as $n\to\infty$ for every $p>1$. 
We recall that the parameters should satisfy 
\beas\label{0211241203}&&
0<\eta_2<\eta_1<\min\left\{\half,\half\left(\frac{1}{c_0}-1\right)\right\},
\quad
\eta_3+\eta_4<\eta_2{\colorb .}
\eeas
{\sred[
In particular, if $c_0=1/2$, then $0<\eta_2<\eta_1<1/2$. 
The positive parameters $\eta_3$ and $\eta_4$ can be sufficiently small 
at this stage. 
Remark that $c_0\eta_2<1/4$ when $c_0\leq1/2$. ]
}
\halflineskip

\noindent
(II) 
Let 
{\sred
\beas 
\wt{\cald}_{n,j}^*
&=& 
{\sred\ol{\kappa}_n^{\eta_3}}\bigg\{
\frac{n}{{\sred\ol{\kappa}_n}}
\sum_{k\in\wh{\calk}_{n,j}(\alpha_0)}|\Delta_kX|^2K_{n,k}
-
\frac{n}{{\sred\ol{\kappa}_n}}
\sum_{k\in\wt{\calk}_{n,j}(\alpha_0)}|\Delta_kX|^2K_{n,k}
\bigg\}{\colorb .}
\eeas
}
{\sred Applying 
Lemma \ref{0211231405} (ii) to 
$V_k=n|\Delta_kX|^2K_{n,k}1_{\{j\in\call_n^c\}}$}, 
we have 
\bea\label{0211241220}
\big\|\wt{\cald}_{n,j}^*1_{\{j\in\call_n^c\}}\big\|_p
&\leq&
\Phi^{(\ref{0211241221})}_{n,j}
+{\colorb \Phi^{(\ref{0211241222})}_{n,j}}
+\Phi^{(\ref{0211241223})}_{n,j},
\eea
where 
\bea\label{0211241221}
\Phi^{(\ref{0211241221})}_{n,j}
&=&
{\sred\ol{\kappa}_n^{\eta_3}}\bigg\|\max_{k\in I_{n,j}}n|\Delta_kX|^2K_{n,k}1_{\{j\in\call_n^c\}}\bigg\|_p
P\bigg[\big||W_{{\sred1}}|-c(\alpha_0)^{1/2}\big|{\colorb <}
\check{C}\kappa_n^{-\eta_2}\bigg], 
\eea
\bea\label{0211241222}
\Phi^{(\ref{0211241222})}_{n,j}
&=&
{\sred\ol{\kappa}_n^{\eta_3}}\bigg\|\max_{k\in I_{n,j}}n|\Delta_kX|^2K_{n,k}1_{\{j\in\call_n^c\}}\bigg\|_{pp_1(p_1-p)^{-1}}
\nn\\&&\times
\bigg\|\frac{1}{{\sred\ol{\kappa}_n}}\sum_{k\in I_{n,j}}\bigg(
1_{\big\{\big|{\colorb |W_k|}-c(\alpha_0)^{1/2}\big|{\colorb <}
	\check{C}\kappa_n^{-\eta_2}\big\}}
-P\bigg[\big||W_k|-c(\alpha_0)^{1/2}\big|{\colorb <}
\check{C}\kappa_n^{-\eta_2}\bigg]
\bigg)\bigg\|_{p_1}
\nn\\&&
\eea
and 
\bea\label{0211241223}
\Phi^{(\ref{0211241223})}_{n,j}
&=&
{\sred\ol{\kappa}_n^{\eta_3}}P\big[\wt{\Omega}_{n,j}^c\big]^{1/p_1}
\bigg\|\max_{k\in I_{n,j}}n|\Delta_kX|^2K_{n,k}1_{\{j\in\call_n^c\}}\bigg\|_{pp_1(p_1-p)^{-1}}
\eea
for $j\in I_n$, $n\in\bbN$. 
Then, 
{\sred paying $\kappa_n^{\eta_4}$ for the maximum,} 
we have the following estimates for any $p_1>p\geq1$: 
\bea\label{0211241241}
\sup_{j\in I_n}\Phi^{(\ref{0211241221})}_{n,j}
&=&
O\big(\kappa_n^{\eta_3{\sred+\eta_4}-\eta_2}\big)
=O\big(n^{-c_0(\eta_2-\eta_3{\sred-\eta_4})}\big),
\eea
\bea\label{0211241242}
\sup_{j\in I_n}\Phi^{(\ref{0211241222})}_{n,j}
&=&
O\big(\kappa_n^{\eta_3}{\sred \times\kappa_n^{\eta_4}}\times\kappa_n^{-(1+\eta_2)/2}\big)
\yeq 
O\big(n^{-c_0\big(\frac{1+\eta_2}{2}-\eta_3{-\sred \eta_4}\big)}\big), 
\eea
and 
\bea\label{0211241243}
\sup_{j\in I_n}\Phi^{(\ref{0211241223})}_{n,j}
&=&
O(n^{-L})
\eea
as $n\to\infty$ for any $L>0$ 
for a sufficiently large $\check{C}$; 
the estimate (\ref{0211241243}) follows from Lemma \ref{0211241300}. 
In this way, 
\bea\label{0211241453}
\big\|\wt{\cald}_{n,j}^*1_{\{j\in\call_n^c\}}\big\|_p
&=&
O\big(n^{-c_0(\eta_2-\eta_3{\sred-{\eta_4}})}\big)
+O\big(n^{-c_0\big(\frac{1+\eta_2}{2}-\eta_3{\sred-{\eta_4}}\big)}\big)
\eea
as $n\to\infty$ for every $p\geq1$. 
\halflineskip

\noindent
(III) 
On the event $\{j\in\call_n^c\}$, we have 
\bea\label{0301021407} 
\sum_{k\in\wt{\calk}_{n,j}(\alpha_0)}|\Delta_kX|^2K_{n,k}
&=&
\sum_{k\in\wt{\calk}_{n,j}(\alpha_0)}
\left(\int_\tkm^\tk \sigma_tdw_t+\int_\tkm^\tk b_tdt\right)^2K_{n,k}
\nn\\&=&
\Phi^{(\ref{0211221131})}_{n,j}+\Phi^{(\ref{0211221132})}_{n,j}+\Phi^{(\ref{0211221133})}_{n,j}+
\Phi^{(\ref{0211221134})}_{n,j}+\Phi^{(\ref{0211221135})}_{n,j}
\eea
where 
\bea\label{0211221131}
\Phi^{(\ref{0211221131})}_{n,j}
&=& 
\sum_{k\in I_{n,j}}\big(\sigma_{t_{\underline{j}_n}}\big)^2
hW_k^21_{\{|W_k|\leq c(\alpha_0)^{1/2}\}},
\eea
\bea\label{0211221132}
\Phi^{(\ref{0211221132})}_{n,j}
&=& 
\sum_{k\in I_{n,j}}\big(\sigma_{t_{\underline{j}_n}}\big)^2
hW_k^21_{\{|W_k|\leq c(\alpha_0)^{1/2}\}}\big(K_{n,k}-1\big)
\nn\\&&
+2\sum_{k\in\wt{\calk}_{n,j}(\alpha_0)}
\int_\tkm^\tk \int_\tkm^t\big(\wt{\sigma}_s-\wt{\sigma}_{t_{\underline{j}_n}}\big)dw_s\sigma_tdw_tK_{n,k}
\nn\\&&
+2\sum_{k\in\wt{\calk}_{n,j}(\alpha_0)}
\int_\tkm^\tk \int_\tkm^t\sigma_{t_{\underline{j}_n}} dw_s\big(\wt{\sigma}_t-\wt{\sigma}_{t_{\underline{j}_n}}\big) dw_tK_{n,k}
\nn\\&&
{\colorr
+2\sum_{k\in\wt{\calk}_{n,j}(\alpha_0)}
\int_\tkm^\tk \wt{\sigma}_{t_{\underline{j}_n}} 
\big(\wt{\sigma}_t-\wt{\sigma}_{t_{\underline{j}_n}}\big)dtK_{n,k}
}
\nn\\&&
{\colorr
+\sum_{k\in\wt{\calk}_{n,j}(\alpha_0)}
\bigg(\int_\tkm^\tk 
\big(\wt{\sigma}_t-\wt{\sigma}_{t_{\underline{j}_n}}\big)dw_t\bigg)^2K_{n,k}
},
\eea
\bea\label{0211221133}
\Phi^{(\ref{0211221133})}_{n,j}
&=& 
2\sum_{k\in\wt{\calk}_{n,j}(\alpha_0)}
\int_\tkm^\tk \int_\tkm^tb_sds\sigma_tdw_tK_{n,k}{\sred,}
\eea
\bea\label{0211221134}
\Phi^{(\ref{0211221134})}_{n,j}
&=& 
2\sum_{k\in\wt{\calk}_{n,j}(\alpha_0)}
\int_\tkm^\tk \int_\tkm^t\sigma_sdw_sb_tdtK_{n,k},
\eea
and 
\bea\label{0211221135}
\Phi^{(\ref{0211221135})}_{n,j}
&=& 
2\sum_{k\in\wt{\calk}_{n,j}(\alpha_0)}
\int_\tkm^\tk \int_\tkm^tb_sdsb_tdtK_{n,k}.
\eea
%

\begin{en-text}
{\colorr 
On $\{j\in\call_n^c\}$, 
\beas
\Phi^{(\ref{0211221132})}_{n,j}
&=& 
\sum_{k\in I_{n,j}}\big(\sigma_{t_{\underline{j}_n}}\big)^2
W_k^21_{\{|W_k|\leq c(\alpha_0)^{1/2}\}}\big(K_{n,k}-1\big)
\nn\\&&
+2\sum_{k\in\wt{\calk}_{n,j}(\alpha_0)}
\int_\tkm^\tk \int_\tkm^t\big(\wt{\sigma}_s-\wt{\sigma}_{t_{\underline{j}_n}}\big)dw_s\sigma_tdw_tK_{n,k}
\nn\\&&
+2\int_\tkm^\tk \int_\tkm^t\sigma_{t_{\underline{j}_n}} dw_s\big(\wt{\sigma}_t-\wt{\sigma}_{t_{\underline{j}_n}}\big) dw_tK_{n,k}
\eeas
}
\end{en-text}
By assumption, 
\bea\label{0211260123}&&
\sup_{j\in I_n}{{\sred \sup_{s\in[t_{\underline{j}_n-1},t_{\underline{j}_n+\ol{\kappa}_n}]}}}
\big\|1_{\{j\in\call_n^c\}}\big(\sigma_{{\sred s}}^2-\sigma_{t_{\underline{j}_n-1}}^2\big)\big\|_p
\nn\\&\yleq &
\sup_{j\in I_n}{{\sred \sup_{s\in[t_{\underline{j}_n-1},t_{\underline{j}_n+\ol{\kappa}_n}]}}}
\big\|\wt{\sigma}_{{\sred s}}^2-\wt{\sigma}_{t_{\underline{j}_n-1}}^2\big\|_p
\nn\\&\simleq&
(\kappa_nh)^{1/2}
\>\simleq\>
h^{\half(1-c_0)}
\eea
for every $p>1$. 
%
First, {\sred a primitive estimate gives}
\begin{en-text}
\bea\label{0211221143}
\frac{n}{\kappa_n}
\big\|\Phi^{(\ref{0211221132})}_{n,j}1_{\{j\in\call_n^c\}}\big\|_p
&\leq& 
\frac{n}{\kappa_n}
\bigg\|\sum_{k\in\wt{\calk}_{n,j}(\alpha_0)}
\int_\tkm^\tk \int_\tkm^t\sigma_sdw_s\sigma_tdw_t\bigg\|_p
\nn\\&&
+
\frac{n}{\kappa_n}
\bigg\|\sum_{k\in\wt{\calk}_{n,j}(\alpha_0)}
\int_\tkm^\tk \int_\tkm^t\sigma_sdw_s\sigma_tdw_t(1-K_{n,k})1_{\{j\in\call_n^c\}}\bigg\|_p
\nn\\&\simleq&
\frac{n}{\kappa_n}\big(\kappa_nh^2)^{1/2}
+
P\bigg[\max_{k\in\wt{\calk}_{n,j}(\alpha_0)}(1-K_{n,k})1_{\{j\in\call_n^c\}}\bigg]^{1/(2p)}
\nn\\&\simleq&
h^{c_0/2}+\kappa_n h^L
\nn\\&\simleq&
h^{c_0/2}
\eea
as $n\to\infty$, where $L$ is any number bigger than $3c_0/2$. 
\end{en-text}
\bea\label{0211221143}
\sup_{j\in I_n}
\frac{n}{\kappa_n}
\big\|\Phi^{(\ref{0211221132})}_{n,j}1_{\{j\in\call_n^c\}}\big\|_p
&\simleq& 
\frac{n}{\kappa_n}\times\frac{\kappa_n^{3/2}}{n^{3/2}}
\>\simleq \>
h^{\half (1-c_0)}
\eea
as $n\to\infty$; 
we note that the orthogonality cannot apply due to $\wt{\calk}_{n,j}(\alpha_0)$ 
even after $K_{n,k}$ is decoupled. 
We also have 
\bea\label{0211221212}
\sup_{j\in I_n}
\frac{n}{\kappa_n}
\big\|\Phi^{(\ref{0211221133})}_{n,j}1_{\{j\in\call_n^c\}}\big\|_p
&\simleq&
{\colorr h^{1/2}}{\colorb .}
\eea
For {\sred$\Phi^{(\ref{0211221134})}_{n,j}$ and $\Phi^{(\ref{0211221135})}_{n,j}$}, 
by the same way, we can get 
\bea\label{0211221214}
\sup_{j\in I_n}
\frac{n}{\kappa_n}
\big\|\Phi^{(\ref{0211221134})}_{n,j}1_{\{j\in\call_n^c\}}\big\|_p
&\simleq&
h^{1/2},
\eea
and 
\bea\label{0211221215}
\sup_{j\in I_n}
\frac{n}{\kappa_n}
\big\|\Phi^{(\ref{0211221135})}_{n,j}1_{\{j\in\call_n^c\}}\big\|_p
&\simleq&
{\colorr h}
\eea
as $n\to\infty$. 
Furthermore, we have
\bea\label{0211260110}&&
\sup_{j\in I_n}
\bigg\|
\bigg\{
\frac{1}{\ol{\kappa}_nh}
\Phi^{(\ref{0211221131})}_{n,j}
-\big(\sigma_{t_{\underline{j}_n}}\big)^2q(\alpha_0)\bigg\}
1_{\{j\in\call_n^c\}}
\bigg\|_p
\nn\\&\leq& 
\sup_{j\in I_n}
\bigg\|
\frac{1}{\ol{\kappa}_nh}\bigg\{
\Phi^{(\ref{0211221131})}_{n,j}
-\sum_{k\in I_{n,j}}\big(\sigma_{t_{\underline{j}_n}}\big)^2q(\alpha_0)h\bigg\}
1_{\{j\in\call_n^c\}}
\bigg\|_p
\nn\\&\leq& 
\sup_{j\in I_n}
\bigg\|
\frac{1}{\ol{\kappa}_n}
\sum_{k\in I_{n,j}}\big(\sigma_{t_{\underline{j}_n}}\big)^2
\big(W_k^21_{\{|W_k|\leq c(\alpha_0)^{1/2}\}}-q(\alpha_0)\big)
\bigg\|_p
\nn\\&=&
O(\kappa_n^{-1/2})
\yeq 
O(h^{c_0/2})
\eea
for every $p>1${\sred.}
\begin{en-text}
\bea\label{0211260119}&&
\sup_{j\in I_n}
\bigg\|
\frac{1}{\ol{\kappa}_nh}
\sum_{k\in I_{n,j}}\bigg\{
\big(\sigma_{t_{\underline{j}_n}}\big)^2q(\alpha_0)h
-\int_\tkm^\tk\sigma_t^2dt
\bigg\}
1_{\{j\in\call_n^c\}}
\bigg\|_p
\nn\\&=&
O(h^{\half(1-c_0)})
\eea
by using (\ref{0211260123}). 
We have \koko
\beas &&
\bigg\|1_{\{j\in\call_n^c\}}\bigg(
\frac{n}{\ol{\kappa}_nT}\sum_{k\in\wt{\calk}_{n,j}(\alpha_0)}\int_\tkm^\tk \sigma_t^2dtK_{n,k}
-q(\alpha_0)\sigma_{t_{\underline{j}_n-1}}^2\bigg)\bigg\|_p
\nn\\&\leq&
\bigg\|1_{\{j\in\call_n^c\}}
\frac{n}{\ol{\kappa}_nT}\sum_{k\in\wt{\calk}_{n,j}(\alpha_0)}\int_\tkm^\tk 
\big(\wt{\sigma}_t^2-\wt{\sigma}_{t_{\underline{j}_n-1}}^2\big)dtK_{n,k}\bigg\|_p
\nn\\&&
+
\bigg\|
\frac{1}{\ol{\kappa}_n}\sum_{k\in I_{n,j}}
\sigma_{t_{\underline{j}_n-1}}^2
\big(1_{\{|W_k|<c(\alpha_0)^{1/2}\}}-q(\alpha_0)\big)\bigg\|_p
\nn\\&&
+
\bigg\|1_{\{j\in\call_n^c\}}
\frac{1}{\ol{\kappa}_n}\sum_{k\in I_{n,j}}
\sigma_{t_{\underline{j}_n-1}}^2
1_{\{|W_k|<c(\alpha_0)^{1/2}\}}
(K_{n,k}-1)\bigg\|_p
\nn\\&\simleq&
(\kappa_nh)^{1/2}+\kappa_n^{-1/2}+O(n^{-L})
\nn\\&\simleq&
(\kappa_nh)^{1/2}+\kappa_n^{-1/2}
\eeas
as well as 
\beas
\big\|1_{\{j\in\call_n^c\}}\big(\sigma_\tk^2-\sigma_{t_{\underline{j}_n-1}}^2\big)\big\|_p
\yleq 
\big\|\wt{\sigma}_\tk^2-\wt{\sigma}_{t_{\underline{j}_n-1}}^2\big\|_p
\>\simleq\>
(\kappa_nh)^{1/2}
\eeas
for every $p>1$. 
\end{en-text}
Combining {\sred (\ref{0301021407}) and} 
(\ref{0211260123})-(\ref{0211260110}), we obtain 
\bea\label{0211221357}&&
\sup_{j\in I_n}
\sup_{k'\in I_{n,j}}
\bigg\|1_{\{j\in\call_n^c\}}\bigg(
\frac{n}{\ol{\kappa}_nT}\sum_{k\in\wt{\calk}_{n,j}(\alpha_0)}|\Delta_kX|^2K_{n,k}
-\sigma_{t_{k'}}^2q(\alpha_0)\bigg)\bigg\|_p
\nn\\&=&
O\big(n^{-(1-c_0)/2}\big)+O(n^{-c_0/2})
\eea
as $n\to\infty$ 
for every $p>1$. 
\begin{en-text}
The estimates 
(\ref{0211221143}), (\ref{0211221212}), 
(\ref{0211221214}), (\ref{0211221215}), 
and (\ref{0211221357}) hold uniformly in $j\in I_n$. 
\end{en-text}

\noindent
(IV) 
From (\ref{0211241151}), (\ref{0211241453}) and (\ref{0211221357}), 
we obtain the estimate 
\bea\label{0301021434}&&
\sup_{j\in I_n}
\sup_{k'\in I_{n,j}}{\sred\kappa_n^{\eta_3}}
\bigg\|1_{\{j\in\call_n^c\}}\bigg(
\frac{n}{\ol{\kappa}_nT}\sum_{k\in{{\sred\calk}}_{n,j}(\alpha_0)}|\Delta_kX|^2K_{n,k}
-\sigma_{t_{k'}}^2q(\alpha_0)\bigg)\bigg\|_p
\nn\\&=&
O\big(n^{-c_0(\eta_2-\eta_3-\eta_4)}\big)
+\bigg\{
O\big(n^{-c_0(\eta_2-\eta_3{\sred-{\eta_4}})}\big)
+O\big(n^{-c_0\big(\frac{1+\eta_2}{2}-\eta_3{\sred-{\eta_4}}\big)}\big)\bigg\}
\nn\\&&
+
{\sred\kappa_n^{\eta_3}}\bigg\{O\big(n^{-(1-c_0)/2}\big)+O(n^{-c_0/2})\bigg\}
\nn\\&=&
O(n^{-c_0(\eta_2-\eta_3-\eta_4)})+O\big(n^{c_0{\sred(\eta_3+\eta_4)}-(1-c_0)/2}\big)
\>=:\> \bbO_n
\eea
as $n\to\infty$ for every $p>1$. 
Here we are assuming the parameters satisfy
\bea\label{0211260224} 
&c_0\in(0,1), \quad B>0, \quad \eta_1\in\bigg(0,\min\bigg\{\half\big(\frac{1}{c_0}-1\big),\half\bigg\}\bigg), &
\nn\\
&\eta_2\in(0,\eta_1),\quad\eta_3>0,\quad\eta_4>0,\quad \eta_3+\eta_4<\eta_2.&
\eea
{\sred 
To obtain the last error bound in (\ref{0301021434}), we used the inequalities
\beas
-c_0\bigg(\frac{1+\eta_2}{2}-\eta_3-\eta_4\bigg)
<
-c_0\big(\eta_2-\eta_3-\eta_4\big)
\eeas
and
\beas 
c_0\eta_3-\frac{c_0}{2}
<
c_0\eta_3-c_0\eta_2
<
-c_0(\eta_2-\eta_3-\eta_4).
\eeas
}
The LGRV $\bbL_{n,j}(\alpha_0)$ of (\ref{0211260215}) does not depend on 
$\eta_i$ ($i=1,2,3,4$) within the ranges (\ref{0211260224}). 
When ${\sred c_0}>1/2$, we make 
\beas 
\half>\half\left(\frac{1}{c_0}-1\right)>\eta_1>\eta_2>\eta_3
\up\half\left(\frac{1}{c_0}-1\right),\quad \eta_4\down0
\eeas
to obtain {\sred$\bbO_n=O(1)$}.
When ${\sred c_0}\leq1/2$, we make 
\beas 
\half>\eta_1>\eta_2>\eta_3\up\half,\quad\eta_4\down0
\eeas
to obtain {\sred$\bbO_n=O(1)$}.
Thus, the proof {\sred of Theorem \ref{0301021509}} is concluded. 
\qed\halflineskip

According to the error bound (\ref{0211221359}), 
we should in general take 
$c_0=1/2$, i.e., 
$\kappa_n\sim  B n^{1/2}$ to obtain an optimal error estimate. 
However, this is not always true. If the process $\sigma$ is (unknown) constant for example, then 
we do not need any spot volatility estimator to construct 
a global jump filter, and the convergence of the resulting estimator for $\Theta$  
becomes much faster than that in the non-constant $\sigma$ case. 

\subsection{Local minimum RV}\label{0301241726}
Estimation of spot volatilities can be done by the minimum realized volatility 
{\tred (minRV)} method 
of Andersen et al. \cite{Andersen2012}. 
{\tred This method is localized to define 
	the local minRV by }
\bea \label{loc_minRV}
\bbM_{n, j} = \frac{\pi}{\pi-2} \frac{n}{\bar{\kappa}_n T} 
\sum_{k \in I_{n,j}} \big\{ |\Delta_k X| \wedge |\Delta_{k+1} X| \}^2. 
\eea
\begin{theorem}
	Suppose that $[G1]$ is fulfilled.
	For $c_0\in(0,1)$ and $B>0$, 
	suppose that $\kappa_n\sim Bn^{c_0}$ as $n\to\infty$. 
	Then 
	\beas
	\sup_{n\in\bbN}
	\sup_{j\in I_n}\sup_{k\in I_{n,j}}n^{{\vred \gamma_{**}}}
	\big\|
	{\fred1_{\{j\in\call_n^c\}}\>\big(
	\bbM_{n,j} - \sigma_\tk^2 \big)}
	\big\|_p
	&<&
	\infty
	\eeas
	as $n\to\infty$ 
{\fred for any $p>1$ and any constant ${\vred \gamma_{**}}$ satisfying 
\beas 
{\vred \gamma_{**}} &{\vred=}& \min\bigg\{\half(1-c_0),\half c_0\bigg\}. 
\eeas
}
\end{theorem}
%
{\fblue
\begin{proof}
Consider $k \in {\vred{ I}}_{n,j}$ for $ j \in {\vred\call}_n^c$. Then we can decompose 
$\Delta_k X$ as 
\beas
\Delta_k X = \sigma_{t_{\underline{j}_n}} \Delta_k w
+  \int_{t_{k-1}}^{t_k} (\sigma_t - \sigma_{t_{\underline{j}_n}}) dw_t 
+ \int_{t_{k-1}}^{t_k} b_t dt.
\eeas
By (\ref{0211260123}), 
\beas&&
\sup_{j \in I_n} \sup_{k \in{\vred { I}}_{n,j}} \Bigg\| 
\int_{t_{k-1}}^{t_k} (\sigma_t - \sigma_{t_{\underline{j}_n}}) dw_t \> 1_{\{ j \in {\vred\call}_n^c \}} \Bigg\|_{{\vred p}}
\nn\\&{\vred\leq} &
{\vred
\sup_{j \in I_n} \sup_{k \in { I}_{n,j}} \Bigg\|  
\int_{t_{k-1}}^{t_k} (\wt{\sigma}_t - \wt{\sigma}_{t_{\underline{j}_n}}) dw_t \Bigg\|_{2p}
\> \big\|1_{\{ j \in \calj_n^c \}} \big\|_{2p}
}\\
&{\vred\simleq}&
\sup_{j \in I_n} \sup_{k \in {\vred I}_{n,j}} 
\sqrt{ 
\int_{t_{k-1}}^{t_k} \Big\|  
{\vred 
\big(\wt{\sigma}_t - \wt{\sigma}_{t_{\underline{j}_n}} \big)^2
} 
\Big\|_{p} dt } \\
&{\vred=}&
\sqrt{O\big( h \times \kappa_n h \big)} = O\big(h^{1-\frac{1}{2}c_0} \big)
\eeas
{\vred for every $p>1$.}
Hence, we obtain 
\beas
|\Delta_k X|^2 = \sigma_{t_{\underline{j}_n}}^2 h W_k^2 + \mathcal{X}_k, 
\eeas
where $\mathcal{X}_k$ is a random variable satisfying 
$\sup_{j \in I_n} \sup_{k \in {\vred I}_{n,j}} \| \mathcal{X}_k\|_p = O(h^{\frac{1}{2}(3-c_0)})$. By using this approximation (and the equality $a \wedge b = \frac{1}{2}(a+b - |a-b|)$ for $a,b > 0$), we have 
\begin{align*}
|\Delta_k X|^2 \wedge |\Delta_{k+1} X|^2 = 
\sigma_{t_{\underline{j}_n}}^2 h \> 
{\vred\big(}W_k^2 \wedge W_{k+1}^2{\vred\big)} +  \mathcal{X}_k',
\end{align*}
where $\mathcal{X}_k'$ is a random variable satisfying 
$\sup_{j \in I_n} \sup_{k \in  {\vred I}_{n,j}} \| \mathcal{X}_k'\|_p = O(h^{\frac{1}{2}(3-c_0)})$. 
Hence, we obtain 
\bea \label{locminRV-1}
\bbM_{n, j} -  \sigma_{t_k}^2
&=&
\frac{\pi}{\pi-2} \frac{n}{\ol{\kappa}_n T} 
\sum_{k \in  {\vred I}_{n,j}} 
{\vred\big\{}|\Delta_k X|^2 \wedge |\Delta_{k+1} X|^2{\vred\big\}}
 -  \sigma_{t_k}^2 \nonumber \\
&=&
\frac{1}{\ol{\kappa}_n} 
\sum_{k \in  {\vred I}_{n,j}} \sigma_{t_{\underline{j}_n}}^2 
\Bigg( \frac{\pi}{\pi-2} {\vred\big\{}W_k^2 \wedge W_{k+1}^2{\vred\big\}} - 1 \Bigg) 
+ \big(\sigma_{t_{\underline{j}_n}}^2  - \sigma_{t_k}^2 \big) 
\nn\\&&
+ \frac{\pi}{\pi-2}  \frac{n}{\ol{\kappa}_n T} \sum_{k \in  {\vred I}_{n,j}}\mathcal{X}_k'
\eea
{\vred The} first term on the right-hand side of (\ref{locminRV-1}) is $O(\kappa_n^{-1/2}) = O(h^{c_0/2})$. 
As for the second term, (\ref{0211260123}) gives 
$\| \sigma_{t_{\underline{j}_n}}^2  - \sigma_{t_k}^2 \|_p = O(h^{\frac{1}{2}(1-c_0)})$. 
Finally, as for the third term, we can estimate as 
\beas
\Bigg\| \frac{n}{\ol{\kappa}_n T}  
\sum_{k \in  {\vred I}_{n,j}}\mathcal{X}_k' \Bigg\|_p 
\lesssim 
n \times O\Big(h^{\frac{1}{2}(3-c_0)}\Big) = 
O\Big(h^{\frac{1}{2}(1-c_0)}\Big). 
\eeas
With these estimates, we obtain the desired result. 
\end{proof}
}

\section{Rate of convergence of the global realized volatilities 
in high intensity of jumps}\label{0211181814}
{\fred 
In this section, we present a rate of convergence of the GRV and WGRV, 
both defined in Section \ref{Sec2}.}
When the frequency of the jumps is high, 
it is recommend that one should choose a value of $\alpha$ that is not extremely small 
in order to cover the jumps by the index set $\calj_n(\alpha)^c$. 
{\fred 
We will assume the properties of $S_{n,j-1}$ below, 
that we already proved in Section \ref{Sec3} for the LGRV and the local minRV. 
Thus, GRV and WGRV with a LGRV or 
the local minRV are global realized volatilities. 
}

\bd
\im[{\bf [G2]}] 
{\bf (i)} $S_{n,j-1}$ is positive a.s. and 
\beas 
{\sred\sup_{n\in\bbN}\sup_{j\in I_n}}\>
\big\|S_{n,j-1}^{-1}\big\|_p&<&\infty
\eeas
for every $p>1$. 

\bd
\im[(ii)] 
There exist positive constants $\gamma_0$ and $c$ such that 
\beas 
\sup_{n\in\bbN}\sup_{j\in I_n}
n^{\gamma_0}\big\|
1_{\{j\in\call_n^c\}}\big(
\sigma_\tjm^2-c\>S_{n,j-1}\big)\big\|_p
&<& \infty
\eeas
for every $p>1$. 
\ed
\ed

In $[G2]$, we do not assume that the value of constant $c$ is known. 
We note that 
\beas 
{\sred\sup_{n\in\bbN}\sup_{j\in I_n}}\>
\big\|1_{\{j\in\call_n^c\}}S_{n,j-1}\big\|_p
&<& \infty
\eeas
for every $p>1$ under $[G1]$ and $[G2]$. 
{\sred As shown in Theorem \ref{0301021509}, 
the LGRV in (\ref{0211260215}) can serve as $S_{n,j-1}$. 
}

If $\sigma_t$ is equal to a (possibly unknown) constant, then 
$\gamma_0$ can be arbitrarily large since we can let $S_{n,j-1}=1$. 
In other words, we do not need any pre-estimate of $\sigma_\tjm^2$. 
So, the constant volatility case is very special and it will 
be discussed briefly in Section \ref{0211281401} separately. 
This section logically includes the constant volatility case 
(hence a less efficient way for it) 
but we will consider a general non-constant volatility 
and assume a given local estimator attains a limited rate of convergence.

\begin{remark}\rm
When ${\sf v}=2^{-1}\inf_{\omega\in\Omega,t\in[0,T]}\sigma_t^2>0$ 
for a priori known constant ${\sf v}$, 
given a local estimator 
$\bbL^{loc}_{n,j-1}$ of $\sigma_\tjm^2$, 
we can use $S_{n,j-1}({\sf v})=\bbL^{loc}_{n,j}\vee{\sf v}$ for $S_{n,j-1}$. 
For example, it is the case when $X$ satisfies a stochastic differential equation 
with jumps and its diffusion coefficient is uniformly elliptic. 
When ${\sf v}=0$, an appropriate modification of $\bbL^{loc}_{n,j}$ is necessary 
and possible. We only give an idea without going into details here. 
Preset a positive constant ${\sf v}$. 
Using $S_{n,j-1}({\sf v})$ for $S_{n,j-1}$, we obtain an estimator $\wt{\bbV}_n[{\sf v}]$ of 
$\Theta({\sf v})=\int_0^T\sigma_t^21_{\{\sigma_t^2\geq{\sf v}\}}dt$, and indeed, 
the rate of convergence $\wt{\bbV}_n[{\sf v}]$ is established in this paper. 
Then it is natural to use $\wt{\bbV}_n[{\sf v}_n]$ to estimate 
$\Theta=\int_0^T\sigma_t^2dt$ with a sequence of numbers ${\sf v}_n$ 
tending to $0$ as $n\to\infty$. 
Consistency does not matter because the mappting ${\sf v}\mapsto \Theta({\sf v})$ is 
continuous and the operation ${\sf v}_n\down0$ is stable. 
Some work is necessary to give an explicit rate of convergence since 
the constant of the error bound for each ${\sf v}_n$ depends on ${\sf v}_n$. 
However, the cause of the error by the truncation at level ${\sf v}_n$ is 
the difference 
$\int_0^T\sigma_t^21_{\{\sigma_t^2<{\sf v}_n\}}dt$, and it is rather easy to control 
for small ${\sf v}_n$. 
\end{remark}

\subsection{Rate of convergence of the GRV with a fixed $\alpha$}
We consider the GRV given by (\ref{0211260834}):
\beas 
\bbV_n(\alpha)
&=& 
\sum_{j\in\calj_n(\alpha)}q(\alpha)^{-1}|\Delta_jX|^2K_{n,j}{\colorb .}
\eeas

Denote by $r_n({\tt U}_j)$ the rank of ${\tt U}_j$ among the variables $\{{\tt U}_i\}_{i\in I_n}$ 
as before, and 
$|{\tt U}|_{(r)}$ denotes the $r$-th ordered statistic of $\{|{\tt U}_i|\}_{i\in I_n}$. 
Let $0<\gamma_2<\gamma_1<\gamma_0$, and 
define numbers $a_n$ and $\wh{a}_n$ by 
\beas 
a_n\yeq\lfloor (1-\alpha)n-n^{1-\gamma_2}\rfloor\quad\text{and}\quad
\wh{a}_n\yeq\lfloor a_n-n^{1-\gamma_2}\rfloor,
\eeas
respectively. 
Define the event $N_{n,j}$ by 
\beas 
N_{n,j} 
&=& 
\big\{r_n(|W_j|)\leq a_n-n^{1-\gamma_2}\big\}\cap
\big\{|W|_{(a_n)}-|W_j|< n^{-\gamma_1}\big\}
\eeas
The following lemma is 
Lemma 2.6 of {\sred Inatsugu} and Yoshida {\colorb\cite{InatsuguYoshida2020Accepted}}{\colorb .}
\begin{lemma}\label{0211260910}
\beas 
P  {\colorb \Bigg[} \bigcup_{j=1,..,,n}N_{n,j} {\colorb \Bigg]}  \yeq O(n^{-L})
\eeas 
as $n\to\infty$ for every $L>0$. 
\end{lemma}

We need some notation: 
{\sred 
\beas 
\wh{\calj}_n(\alpha)
&=& 
\big\{j\in I_n;\> r_n(|W_j|)\leq \wh{a}_n\big\},\quad
\nn\\
U_j &=& {\sred c}^{-1/2}h^{-1/2}(S_{n,j-1})^{-1/2}\Delta_jX
\nn\\
R_j &=& U_j-W_j- {\sred c}^{-1/2}h^{-1/2}(S_{n,j-1})^{-1/2}\Delta_jJ,
\eeas
}
as well 
\begin{en-text}
\beas
\Omega_n &=& \bigg\{\ol{N}_T< n^{1-\gamma_2}\bigg\} {\sred \bigcap}
\bigg(\bigcap_{j=1,...,n}\bigg[\big\{|R_j|1_{\{\Delta_jN^\sigma=0\}}<2^{-1}n^{-\gamma_1}\big\}\cap(N_{n,j})^c\bigg]
\bigg). 
\eeas
\end{en-text}
\beas
\Omega_n &=& \bigg\{
{\xred\#\>\call_n}
< n^{1-\gamma_2}\bigg\} {\sred \bigcap}
\bigg(\bigcap_{j=1,...,n}\bigg[\big\{|R_j|1_{\{{\xred j\in\call_n^c}\}}<2^{-1}n^{-\gamma_1}\big\}\cap(N_{n,j})^c\bigg]
\bigg). 
\eeas

{\xred 
Let 
\bea\label{0211280352}
{\mathfrak L}_n &=& \big\{j\in I_n;\> \Delta_j\ol{N}\not=0\big\}.
\eea
The definition of $\call_n^{(k)}$ 
in Inatsugu and Yoshida \cite{InatsuguYoshida2020Accepted} 
of the extended version arXiv:1806.10706v3 is essentially 
the same as ${\mathfrak L}_n$, and different from 
$\call_n$ defined by (\ref{0301021550}). 
The random set ${\sf L}_n^{(k)}$ therein corresponds to $\call_n$. 
%
}

{\xred 
We assume that the distribution of the variable $\ol{N}_T$ depends on $n$, 
and consider the case where $\ol{N}_T$ may diverge as $n\to\infty$. 
More precisely, we will assume the following situation.
\footnote{We slightly relaxed Condition $[G3]$ of arXiv:2102.05307v1.}
}
{\sred 
\bd
\im[{\bf [G3]}]
There exists a constant $\xi\geq0$ 
such that 
$\|{\xred\#\>{\mathfrak L}_n}\|_p=O(n^{\xi})$ 
as $n\to\infty$ 
for every $p>1$. 
\ed
}

\begin{lemma}\label{0211261336}
Suppose that $[G1]$ and $[G2]$ are satisfied. 
Suppose that $0<\gamma_1<\gamma_0<1/2$. Then 
\bea\label{02112q61347} 
\sup_{j\in I_n} P\big[
{\sred |R_j|}1_{\{{\xred j\in\call_n^c}\}}\geq2^{-1}n^{-\gamma_1}\big]
&=& 
O(n^{-L})
\eea
as $n\to\infty$ for every $L>0$. 
In particular, 
{\sred if 
{\xred the conditions $\kappa_n=O(n^{1/2})$, 
$\gamma_2<\half-\xi$ and $[G3]$ are additionally satisfied}, then}
\bea\label{0211261348}
P[\Omega_n^c] &=& O(n^{-L})
\eea
as $n\to\infty$ for every $L>0$. 
\end{lemma}
\proof 
We have 
\beas
\sup_{{\colorb j} \in I_n}\big\|R_j1_{\{{\xred j\in\call_n^c}\}}\big\|_p
&=& 
O(n^{-\gamma_0})
\eeas
for every $p>1$. 
The Markov inequality implies (\ref{02112q61347}). 
This estimate and Lemma \ref{0211260910} give (\ref{0211261348}) 
{\xred 
if the Markov inequality is used with the estimate $\|\#\>\call_n\|_p=O(n^{\xi+1/2})$ from $[G3]$. 
}
\qed\halflineskip

{\xred 
Lemma 2.7 of Inatsugu and Yoshida {\colorb\cite{InatsuguYoshida2020Accepted}} 
(or see an extended version arXiv:1806.10706v3)
is rephrased as follows. 
Recall that $\call_n$ is given by (\ref{0301021550})}.
\begin{lemma} 
\bea\label{300217-1} 
\wh{\calj}_n(\alpha)\cap{\xred \call_n^c}
&\subset& 
\calj_n(\alpha) 
\eea
on $\Omega_n$. 
In particular 
\bea\label{300217-2} 
\# \big[\calj_n(\alpha)\ominus\wh{\calj}_n(\alpha)\big] &\leq& 
c_*n^{1-\gamma_2}
+{\xred \#\>\call_n
}
\eea
on $\Omega_n$, where $c_*$ is a positive constant. 
{\sred Here $\ominus$ denotes the symmetric difference operator of sets.  }
\end{lemma}
\halflineskip

For $\gamma_3>0$ and random variables $({\tt U}_j)_{j=1,...,n}$, let 
\beas 
\cald_n &=& n^{\gamma_3}
\Bigg| \frac{1}{n}\sum_{j\in\calj_n(\alpha)}{\tt U}_j-\frac{1}{n}\sum_{j\in\wh{\calj}_n(\alpha)}{\tt U}_j \Bigg|.
\eeas

{\xred 
We refer the reader to Lemmas 2.8 and 2.9 of Inatsugu and Yoshida \cite{InatsuguYoshida2020Accepted} (or see arXiv:1806.10706v3)
for proof of the following two lemmas. 
}
\begin{lemma}\label{300217-5}
{\bf (i)} Let $p_1>1$. Then
\beas 
\|\cald_n\|_p&\leq& 
\big(c_*n^{\gamma_3-\gamma_2}+n^{-1+\gamma_3}\|{\xred \#\>\call_n}\|_{p_1}\big)
\bigg\|\max_{j=1,...,n}\big|{\tt U}_j\big|\bigg\|_{pp_1(p_1-p)^{-1}}
\\&&
+n^{\gamma_3}\bigg\|\max_{j=1,...,n}\big|{\tt U}_j\big|1_{\Omega_n^c}\bigg\|_{p}
\eeas
for $p\in(1,p_1)$. 
\bd
\im[(ii)] Let $\gamma_4>0$ and $p_1>1$. Then
\beas 
\|\cald_n\|_p&\leq& 
\big(c_*n^{\gamma_3-\gamma_2}+n^{-1+\gamma_3}\|{\xred \#\>\call_n}\|_{p_1}\big)
\nn\\&&\times
\bigg(n^{\gamma_4}
+n\max_{j=1,...,n}\bigg\|\big|{\tt U}_j\big|1_{\{|{\tt U}_j|>n^{\gamma_4}\}}\bigg\|_{pp_1(p_1-p)^{-1}}
\bigg)
\\&&
+n^{\gamma_3}\bigg\|\max_{j=1,...,n}\big|{\tt U}_j\big|1_{\Omega_n^c}\bigg\|_{p}
\eeas
for $p\in(1,p_1)$. 
\ed
\end{lemma}
\halflineskip

Let 
\beas 
\wt{\cald}_n &=& n^{\gamma_3}
\bigg|\frac{1}{n}\sum_{j\in\wh{\calj}_n(\alpha)}{\tt U}_j-\frac{1}{n}\sum_{j\in\widetilde{\calj}_n(\alpha)}{\tt U}_j\bigg|.
\eeas
for a collection of random variables $\{{\tt U}_j\}_{j\in I_n}$ and 
\bea\label{0211280343} 
\wt{\calj}_n(\alpha)
&=& 
\big\{j\in I_n;\> |W_j|\leq c(\alpha)^{1/2}\big\}. 
\eea
Let 
\bea\label{0211280310}
\wt{\Omega}_n
&=&
\big\{\big| |W|_{({\sred\widehat{a}}_n)}-{\sred c(\alpha)^{1/2}}\big|
< \check{C}\>n^{-\gamma_2}\big\},
\eea
where $\check{C}$ is a positive constant. 
See Lemma 4 of Inatsugu and Yoshida {\colorb\cite{InatsuguYoshida2020Accepted}} for a proof of the following lemma. 
\begin{lemma}\label{0211261213}
Let $\check{C}>0$ and $\gamma_3>0$. Then 
\bd
\im[(i)] For $p\geq1$, 
\beas 
\|\wt{\cald}_n\|_p
&\leq& 
n^{\gamma_3}\bigg\|\max_{j'=1,...,n}|{\tt U}_{j'}|\>
\frac{1}{n}\sum_{j=1}^n
1_{\big\{ \big||W_j|-{\sred c(\alpha)^{1/2}}\big|{\colorb <}
	\check{C}\>n^{-\gamma_2}\big\}}
\bigg\|_p
+n^{\gamma_3}\bigg\|
1_{\wt{\Omega}_n^c}
\max_{j'=1,...,n}|{\tt U}_{j'}|\bigg\|_p.
\eeas

\im[(ii)] For $p_1>p\geq1$, 
\beas 
\|\tilde{\cald}_n^{(k)}\|_p
&\leq& 
n^{\gamma_3}\bigg\|\max_{j=1,...,n}|{\tt U}_{j}|\>
\bigg\|_p
P\bigg[\big||W_1|-{\sred c(\alpha)^{1/2}}\big|{\colorb <}
\check{C}\>n^{-\gamma_2}\bigg]
\\&&
+n^{\gamma_3}\bigg\|\max_{j=1,...,n}|{\tt U}_{j}|
\bigg\|_{pp_1(p_1-p)^{-1}}
\nn\\&&\hspace{30pt}\times
\bigg\|
\frac{1}{n}\sum_{j=1}^n\bigg(
1_{\big\{ \big|{\sred|W_j|-}{\sred c(\alpha)^{1/2}}\big| {\colorb <}
	\check{C}\>n^{-\gamma_2}\big\}}
-
P\bigg[\big||W_1|-{\sred c(\alpha)^{1/2}}\big|{\colorb <}
\check{C}\>n^{-\gamma_2}\bigg]
\bigg)
\bigg\|_{p_1}
\\&&
+n^{\gamma_3}
P[\wt{\Omega}_n^c]^{1/p_1}
\bigg\|\max_{j=1,...,n}|{\tt U}_{j}|\bigg\|_{pp_1(p_1-p)^{-1}}.
\eeas
\ed
\end{lemma}
\halflineskip

%
%
\begin{lemma}\label{021127f1234}
If the constant $\check{C}$ in (\ref{0211241445}) is sufficiently large, then 
\beas 
P\big[\wt{\Omega}_n^c\big]
&=& 
O(n^{-L})
\eeas
as $n\to\infty$ for any $L>0$. 
\end{lemma}

{\sred 
Now we shall investigate the rate of convergence of $\bbV_n(\alpha)$ 
for a constant $\alpha\in(0,1)$. 
}
We note that, under $[G1]$ and $[G3]$, 
\bea\label{0211261352}
\bigg\|
\sum_{j\in
{\mathfrak L}_n}|\Delta_jX|^2K_{n,j}\bigg\|_p
&\leq&
n^{-1/2}\big\|{\xred\#\>{\mathfrak L}_n}\big\|_p
\yeq 
O(n^{-1/2+\xi}){\sred.}
\eea

Let 
\beas 
\wh{\bbV}_n(\alpha)
&=& 
\sum_{j\in\wh{\calj}_n(\alpha)}q(\alpha)^{-1}|\Delta_jX|^2K_{n,j}.
\eeas
\begin{lemma}\label{0211271153}
Suppose that $[G1]$ $[G2]$ and $[G3]$ are fulfilled. 
Suppose that $\xi<\half$. 
Let $\gamma_5<\min\big\{\gamma_0,\half-\xi\big\}$ 
{\xred and $\kappa_n=O(n^{1/2})$}. Then 
\beas 
\sup_{n\in\bbN}n^{\gamma_5}\big\|\bbV_n(\alpha)-\wh{\bbV}_n(\alpha)\big\|_p
&<& 
\infty. 
\eeas
\end{lemma}
\proof
By  (\ref{0211261352}), we obtain 
\beas 
\big\|\bbV_n(\alpha)-\wh{\bbV}_n(\alpha)\big\|_p
&=& 
\bigg\|\sum_{j\in\calj_n(\alpha)}
q(\alpha)^{-1}|\Delta_jX|^21_{\{\Delta_j\ol{N}=0\}}K_{n,j}
\nn\\&&\hspace{30pt}
-\sum_{j\in\wh{\calj}_n(\alpha)}
q(\alpha)^{-1}|\Delta_jX|^21_{\{\Delta_j\ol{N}=0\}}K_{n,j}
\bigg\|_p
+O(n^{-1/2+\xi}){\sred.}
\eeas
By Lemmas \ref{300217-5} and \ref{0211261336}, 
\beas &&
n^{\gamma_3}\big\|\bbV_n(\alpha)-\wh{\bbV}_n(\alpha)\big\|_p
\nn\\&\simleq& 
\big(c_*n^{\gamma_3-\gamma_2}+n^{-1+\gamma_3}\|{\xred \#\>\call_n}\|_{p_1}\big)
\nn\\&&\times
\bigg(n^{\gamma_4}
+n\max_{j=1,...,n}\bigg\|n|\Delta_jX|^21_{\{\Delta_j\ol{N}=0\}}K_{n,j}1_{\{n|\Delta_jX|^21_{\{\Delta_j\ol{N}=0\}}K_{n,j}>n^{\gamma_4}\}}\bigg\|_{pp_1(p_1-p)^{-1}}
\bigg)
\\&&
+n^{\gamma_3}\bigg\|\max_{j=1,...,n}\big(n|\Delta_jX|^21_{\{\Delta_j\ol{N}=0\}}K_{n,j}\big)1_{\Omega_n^c}\bigg\|_{p} {\colorb +O(n^{-1/2+\gamma_3+\xi})}
\nn\\&\simleq& 
c_*n^{\gamma_3+\gamma_4-\gamma_2}+n^{{\xred-1/2}+\gamma_3+\gamma_4+\xi}
+n^{-1/2+\gamma_3+\xi}
,
\eeas
where $1<p<p_1$. 
The parameters should satisfy
\beas 
0<\gamma_3<\gamma_2<\gamma_1<\gamma_0<\half,
\quad
\gamma_2<\half-\xi,
\quad 
\gamma_4>0.
\eeas
We make 
\beas 
\gamma_4\down0,\qquad
{\colorb \gamma_5 < }\gamma_3<\up\gamma_2<\up\gamma_1<\up\min\bigg\{\gamma_0,\half-\xi\bigg\}
\eeas
to obtain the desired exponent.
\qed\halflineskip

For $\wt{\calj}_n(\alpha)$ defined in (\ref{0211280343}), let 
\beas 
\wt{\bbV}_n(\alpha)
&=& 
\sum_{j\in\wt{\calj}_n(\alpha)}q(\alpha)^{-1}|\Delta_jX|^2K_{n,j}.
\eeas
\begin{lemma}\label{0211271219}
Suppose that $[G1]$ and $[G3]$ are fulfilled. 
Suppose that $\xi<\half$. 
Let 
$\gamma_6<\half-\xi$. 
Then 
\beas 
\sup_{n\in\bbN}
n^{\gamma_6}\big\|\wh{\bbV}_n(\alpha)-\wt{\bbV}_n(\alpha)\big\|_p
&<&
\infty{\colorb .}
\eeas
\end{lemma}
\proof
By  (\ref{0211261352}), we obtain 
\beas 
\big\|\wh{\bbV}_n(\alpha)-\wt{\bbV}_n(\alpha)\big\|_p
&=& 
\bigg\|\sum_{j\in\wh{\calj}_n(\alpha)}
q(\alpha)^{-1}|\Delta_jX|^21_{\{\Delta_j\ol{N}=0\}}K_{n,j}
\nn\\&&\hspace{30pt}
-\sum_{j\in\wt{\calj}_n(\alpha)}
q(\alpha)^{-1}|\Delta_jX|^21_{\{\Delta_j\ol{N}=0\}}K_{n,j}
\bigg\|_p
+O(n^{-1/2+\xi}){\colorb .}
\eeas
By Lemma \ref{0211261213}, we obtain 
\beas &&
n^{\gamma_3}\big\|\wh{\bbV}_n(\alpha)-\wt{\bbV}_n(\alpha)\big\|_p
\nn\\&\simleq& 
n^{\gamma_3}\bigg\|\max_{j=1,...,n}\big(n|\Delta_jX|^21_{\{\Delta_j\ol{N}=0\}}K_{n,j}\big)\>
\bigg\|_p
P\bigg[\big||W_1|-{\sred c(\alpha)^{1/2}}\big|{\colorb <}
\check{C}\>n^{-\gamma_2}\bigg]
\\&&
+n^{\gamma_3}\bigg\|\max_{j=1,...,n}\big(n|\Delta_jX|^21_{\{\Delta_j\ol{N}=0\}}K_{n,j}\big)
\bigg\|_{pp_1(p_1-p)^{-1}}
\nn\\&&\hspace{30pt}\times
\bigg\|
\frac{1}{n}\sum_{j=1}^n\bigg(
1_{\big\{ \big|{\sred |W_j|-c(\alpha)^{1/2}}\big|{\colorb <}
	\check{C}\>n^{-\gamma_2}\big\}}
-
P\bigg[\big||W_1|-{\sred c(\alpha)^{1/2}}\big|{\colorb <}
\check{C}\>n^{-\gamma_2}\bigg]
\bigg)
\bigg\|_{p_1}
\\&&
+n^{\gamma_3}
P[\wt{\Omega}_n^c]^{1/p_1}
\bigg\|\max_{j=1,...,n}\big(n|\Delta_jX|^21_{\{\Delta_j\ol{N}=0\}}K_{n,j}\big)\bigg\|_{pp_1(p_1-p)^{-1}}
+O(n^{\gamma_3-1/2+\xi})
\nn\\&\simleq&
n^{\gamma_3-\gamma_2}+n^{\gamma_3-\half-\frac{\gamma_2}{2}}
{\sred+n^{\gamma_3-1/2+\xi}}
\>\simleq\>
n^{\gamma_3{\sred+\gamma_4}-\gamma_2}+n^{\gamma_3{\sred+\gamma_4}-1/2+\xi},
\eeas
where $1\leq p<p_1$ and {\sred$\gamma_4$ is {\colorb an} arbitrary positive number}. 
{\sred Lemma \ref{021127f1234} was used in the above derivation.}
Making 
\beas 
{\sred 
\gamma_4\down0\quad\text{and}\quad}
\gamma_6<
\gamma_3<\up\gamma_2<\up
\half-\xi,
\eeas
we conclude the proof. 
\qed\halflineskip

\begin{lemma}\label{0211280809}
Suppose that ${\sred[G1]}$ and $[G3]$ are satisfied. Suppose that $\xi<1/2$. 
Then 
\beas
\bigg\|\wt{\bbV}_n(\alpha)
-\int_0^T \sigma_t^2dt\bigg\|_p
&=&
O\big(n^{-\half+\xi}\big)
\eeas
as $n\to\infty$ for every $p>1$. 
\end{lemma}
\proof 
Recall that ${\sred\mathfrak L_n}$ is defined by (\ref{0211280352}). 
We have 
\bea\label{0211280728} 
\sum_{j\in\wt{\calj}_n(\alpha)}|\Delta_jX|^2K_{n,j}
1_{\{j\in{\mathfrak L}_n^c\}}
&=&
\sum_{j\in\wt{\calj}_n(\alpha)}
\left(\int_\tjm^\tj \sigma_tdw_t+\int_\tjm^\tj b_tdt\right)^2K_{n,j}1_{\{j\in{\mathfrak L}_n^c\}}
\nn\\&=&
\Phi^{(\ref{0211280401})}_n+\Phi^{(\ref{0211280402})}_n+
\Phi^{(\ref{0211280403})}_n
\eea
where 
\bea\label{0211280401}
\Phi^{(\ref{0211280401})}_n
&=& 
\sum_{j\in I_n	}\sigma_\tjm^2
hW_j^21_{\{|W_j|\leq c(\alpha)^{1/2}\}},
\eea
\bea\label{0211280402}
\Phi^{(\ref{0211280402})}_n
&=& 
\sum_{j\in I_n}\sigma_\tjm^2
hW_j^21_{\{|W_j|\leq c(\alpha_0)^{1/2}\}}
\big(K_{n,j}1_{\{j\in{\mathfrak L}_n^c\}}
-1\big)
\nn\\&&
+2\sum_{j\in\wt{\calj}_n(\alpha)}
\int_\tjm^\tj \int_\tjm^t\big(\wt{\sigma}_s-\wt{\sigma}_\tjm\big)dw_s\sigma_tdw_tK_{n,j}
1_{\{j\in{\mathfrak L}_n^c\}}
\nn\\&&
+2\sum_{j\in\wt{\calj}_n(\alpha)}
\int_\tjm^\tj \int_\tjm^t\sigma_\tjm dw_s\big(\wt{\sigma}_t-\wt{\sigma}_\tjm\big) dw_tK_{n,j}
1_{\{j\in{\mathfrak L}_n^c\}}
\nn\\&&
{\colorr
+2\sum_{j\in\wt{\calj}_n(\alpha)}
\int_\tjm^\tj \wt{\sigma}_\tjm
\big(\wt{\sigma}_t-\wt{\sigma}_\tjm\big)dtK_{n,j}1_{\{j\in{\mathfrak L}_n^c\}}
}
\nn\\&&
{\colorr
+\sum_{j\in\wt{\calj}_n(\alpha)}
\bigg(\int_\tjm^\tj 
\big(\wt{\sigma}_t-\wt{\sigma}_\tjm\big)dw_t\bigg)^2K_{n,j}1_{\{j\in{\mathfrak L}_n^c\}}
}
\eea
and 
\bea\label{0211280403}
\Phi^{(\ref{0211280403})}_n
&=& 
2\sum_{j\in\wt{\calj}_n(\alpha)}
\int_\tjm^\tj \int_\tjm^tb_sds\sigma_tdw_tK_{n,j}1_{\{j\in{\mathfrak L}_n^c\}}
\nn\\&&+
2\sum_{j\in\wt{\calj}_n(\alpha)}
\int_\tjm^\tj \int_\tjm^t\sigma_sdw_sb_tdtK_{n,j}1_{\{j\in{\mathfrak L}_n^c\}}
\nn\\&&+
2\sum_{j\in\wt{\calj}_n(\alpha)}
\int_\tjm^\tj \int_\tjm^tb_sdsb_tdtK_{n,j}1_{\{j\in{\mathfrak L}_n^c\}}.
\eea
%
%

Let $\ep>\xi$. 
For $p>1$ and $\ep'>0$, 
\bea\label{0211280645} &&
\bigg\|
\sum_{j\in I_n}\sigma_\tjm^2
hW_j^21_{\{|W_j|\leq c({\sred\alpha})^{1/2}\}}
\big(K_{n,j}1_{\{j\in{\mathfrak L}_n^c\}}
-1\big)\bigg\|_p
\nn\\&\leq&
\bigg\|
\sum_{j\in I_n}\sigma_\tjm^2
hW_j^21_{\{|W_j|\leq c({\sred\alpha})^{1/2}\}}
{\xred1_{\{j\in{\mathfrak L}_n\}}}
\bigg\|_p
\\&&
+\bigg\|
\sum_{j\in I_n}\sigma_\tjm^2
hW_j^21_{\{|W_j|\leq c({\sred\alpha})^{1/2}\}}
\big(K_{n,j}
-1\big){\xred1_{\{j\in{\mathfrak L}_n^c\}}}\bigg\|_p
\nn\\&\leq&
\bigg\|
\max_{j\in I_n}\bigg(\sigma_\tjm^2hW_j^21_{\{|W_j|\leq c({\sred\alpha})^{1/2}\}}K_{n,j}\bigg)
{\>\xred\#{\mathfrak L}_n}
\bigg\|_p
+O(n^{-L})
\nn\\&\leq&
\big\|{\xred\#{\mathfrak L}_n}1_{\{{\xred\#{\mathfrak L}_n}>n^{\ep}\}}\big\|_{2p}
\bigg\|
\max_{j\in I_n}\bigg(\sigma_\tjm^2hW_j^21_{\{|W_j|\leq c({\sred\alpha})^{1/2}\}}K_{n,j}\bigg)
\bigg\|_{2p}
\nn\\&&
+
n^{\ep}
\bigg\|
\max_{j\in I_n}\bigg(\sigma_\tjm^2hW_j^21_{\{|W_j|\leq c({\sred\alpha})^{1/2}\}}K_{n,j}\bigg)
\bigg\|_{p}
+O(n^{-L})
\nn\\&\simleq&
n^{-\frac{L\ep}{2p}}\big\|{\xred\#{\mathfrak L}_n}\big\|_{2p+L}^{\frac{2p+L}{2p}}\times n^{-1+\ep'}
+n^{-1+\ep+\ep'}
+O(n^{-L})
\nn\\&\simleq&
n^{-\frac{L(\ep-\xi)}{2p}+\xi-1+\ep'}+n^{-1+\ep+\ep'}+O(n^{-L})
\nn\\&\simleq&
n^{-1+\ep+\ep'}
\nn
\eea
since $\ep>\xi$, where $L$ is a sufficiently large number 
chosen suitably depending on $(\ep,\xi,p,\ep')$. 

From the estimate (\ref{0211280645}), we have 
\bea\label{0211280407}
\big\|\Phi^{(\ref{0211280402})}_n
\big\|_p
&\simleq& 
h^{1/2}+h^{1-\ep-\ep'}
\>\simleq\>
h^{1/2}
\eea
if letting $\ep\down\xi<1/2$ and $\ep'\down0$. 
%

{\colorb
By the Burkholder-Davis-Gundy inequality, we have 
\begin{align*}
	\Bigg\|  \sum_{j \in \widetilde{\calj}_n(\alpha)} 
	\int_{\tjm}^{t_j} \int _{\tjm}^{t} b_s ds \sigma_t dw_t K_{n, j} 1_{\{j \in \mathfrak{L}_n^c \}} \Bigg\|_p 
	&\leq
	\sum_{j \in I_n} \Bigg\|  
	\int_{\tjm}^{t_j} \int _{\tjm}^{t} b_s ds \sigma_t dw_t  \Bigg\|_p \\
	&\lesssim
	\sum_{j \in I_n} \sqrt{
	\Bigg\|  
	\int_{\tjm}^{t_j} \Bigg( \int _{\tjm}^{t} b_s ds \sigma_t \Bigg)^2 dt  \Bigg\|_{p/2}} \\
	&\leq
	\sum_{j \in I_n} \sqrt{ 
	\int_{\tjm}^{t_j} 	\Bigg\|  
	\int _{\tjm}^{t} b_s ds \sigma_t \Bigg\|_p^2 dt }\\
	&\lesssim
	h^{1/2}. 
\end{align*}
From this and similar estimates, we have 
}
\bea\label{0211280408}
\big\|\Phi^{(\ref{0211280403})}_n
\big\|_p
&\simleq&
h^{1/2}
\eea
as $n\to\infty$ for every $p>1$. 
Moreover, 
\bea\label{0211280411}
\bigg\|
\Phi^{(\ref{0211280401})}_j
-\sum_{j\in I_n}\sigma_\tjm^2 q({\sred\alpha})h
\bigg\|_p
&\leq& 
\bigg\|h
\sum_{j\in I_n}\sigma_\tjm^2
\big(W_j^21_{\{|W_j|\leq c({\sred\alpha})^{1/2}\}}-q({\sred\alpha})\big)
\bigg\|_p
\nn\\&=&
O(h^{1/2})
\eea
for every $p>1$.

Obviously, 
\bea\label{0211280406}
\sup_{j\in I_n}
\big\|1_{\{j\in{\mathfrak L}_n^c\}}\big(\sigma_\tk^2-\sigma_\tjm^2\big)
\big\|_p
&\yleq &\sup_{j\in I_n}\big\|\wt{\sigma}_\tk^2-\wt{\sigma}_\tjm^2\big\|_p
\>\simleq\>
h^{1/2}
\eea
for every $p>1$. 
In view of (\ref{0211280406}), we deduce that 
\bea\label{0211280745}&&
\bigg\|
\sum_{j\in I_n}\sigma_\tjm^2 h
-\int_0^T \sigma_t^2dt
\bigg\|_p
\nn\\&\leq& 
\bigg\|
\sum_{j\in I_n}\int_\tjm^\tj \big|\wt{\sigma}_t^2
-\wt{\sigma}_\tjm^2 \big| dt
\bigg\|_p
+
\bigg\|
\sum_{j\in I_n}\int_\tjm^\tj \big(\sigma_t^2
-\sigma_\tjm^2 \big) dt1_{\{j\in{\mathfrak L}_n\}}
\bigg\|_p
\nn\\&\leq& 
O(h^{1/2})
+
\bigg\|
\max_{j\in I_n}\bigg\{\int_\tjm^\tj \big(|\sigma_t^2|
+|\sigma_\tjm^2| \big) dt\bigg\}\>
{\xred\#{\mathfrak L}_n}
\bigg\|_p
\nn\\&=& 
O(h^{1/2}),
\eea
following the passage from (\ref{0211280645}) to (\ref{0211280407}). 

Easily, 
\bea\label{0211280722}
\bigg\|
\sum_{j\in\wt{\calj}_n(\alpha)}|\Delta_jX|^2K_{n,j}
1_{\{j\in{\mathfrak L}_n\}}
\bigg\|_p
&\leq&
\big\|n^{-1/2}{\xred\#{\mathfrak L}_n}
\big\|_p
\>\simleq\>
n^{-\half+\xi}. 
\eea
Combining (\ref{0211280722}), (\ref{0211280728}), 
(\ref{0211280407}), (\ref{0211280408})
(\ref{0211280411}) and (\ref{0211280745}), 
we obtain 
\beas
\bigg\|\wt{\bbV}_n(\alpha)
-\int_0^T \sigma_t^2dt\bigg\|_p
&=&
O\big(n^{-\half+\xi}\big)
\eeas
as $n\to\infty$ for every $p>1$. 
\qed\halflineskip

\begin{theorem}\label{0211280814}
Suppose that $[G1]$ $[G2]$ and $[G3]$ are fulfilled. 
Suppose that $\xi<\half$ {\xred and $\kappa_n=O(n^{1/2})$}. 
Let $\alpha\in(0,1)$ and $\beta_0<\min\big\{\gamma_0,\half-\xi\big\}$. Then 
\beas 
\big\|  \bbV_n(\alpha)-\Theta \big\|_p
&=& 
O(n^{-\beta_0})
\eeas
as $n\to\infty$ for every $p>1$. 
\end{theorem}
\proof 
Use Lemmas \ref{0211271153}, \ref{0211271219} and \ref{0211280809}. 
\qed\halflineskip

\subsection{Rate of convergence of the WGRV with a fixed $\alpha$}
%
Next, we discuss the convergence of the WGRV with a fixed $\alpha$. 
Recall that the WGRV is defined as 

$$
\bbW_n(\alpha) = \sum_{j \in I_n} {\sf w}(\alpha)^{-1}  
\big\{ |\Delta_j X| \wedge (S_{n,j-1}^{1/2} V_{(s_n(\alpha))}) \big\}^2 K_{n,j}. 
$$

The WGRV has entirely the same rate of convergence as the GRV. 
\begin{theorem}
	Suppose that {\colorb $[G1]$},  {\colorb $[G2]$}, and {\colorb $[G3]$} are fulfilled. Suppose that 
	$\xi < \frac{1}{2}$. Let $\alpha \in (0, 1)$ and 
	$\beta_0 {\xred<} \min \big\{ \gamma_0, \frac{1}{2} - \xi \big\}$. Moreover, assume that $\kappa_n{\xred =O( n^{1/2})}$. Then 
	$$
	\left\| \bbW_n(\alpha) - \Theta \right\|_p = O(n^{-\beta_0})
	$$
	as $n \to \infty$ for every $p > 1$. 
\end{theorem}

\proof Decompose {\colorb $\bbW_n(\alpha)$} as 
\begin{align*}
	{\colorb \bbW_n(\alpha)}
	&=
	\sum_{j \in \calj_n(\alpha)} {\sf w}(\alpha)^{-1}  |\Delta_j X|^2 K_{n,j} 
	+  \sum_{j \in \calj_n(\alpha)^c} {\sf w}(\alpha)^{-1}  S_{n,j-1} V_{(s_n(\alpha))}^2 K_{n,j}\\ 
	&=
	\frac{q(\alpha)}{{\sf w}(\alpha)} {\bbV}_n(\alpha) 
	+  \sum_{j \in \calj_n(\alpha)^c} {\sf w}(\alpha)^{-1}  S_{n,j-1} V_{(s_n(\alpha))}^2 K_{n,j}. 
\end{align*}
Note that ${\sf w}(\alpha) = q(\alpha) + \alpha c(\alpha)$. 
Hence, it suffices to show that 
$$
\left \| \sum_{j \in \calj_n(\alpha)^c} S_{n,j-1} V_{(s_n(\alpha))}^2 K_{n,j} 
- \alpha c(\alpha) \Theta \right\|_p = O(n^{-\beta_0})
$$
{\colorb as $n \to \infty$ for every $p>1$.} 
Decompose the left-hand side as 
\begin{align*}
	\sum_{j \in \calj_n(\alpha)^c} S_{n,j-1} V_{(s_n(\alpha))}^2 K_{n,j} 
	- \alpha c(\alpha) \Theta 
	&= 
	\sum_{j \in \calj_n(\alpha)^c} S_{n,j-1} V_{(s_n(\alpha))}^2 K_{n,j}
	1_{\{j \in \call_n^c \}} - \alpha c(\alpha) \Theta  \\
	&\qquad + \sum_{j \in \calj_n(\alpha)^c} S_{n,j-1} V_{(s_n(\alpha))}^2 K_{n,j} 1_{\{j \in \mathfrak{L}_n \cap \call_n\}} \\
	&\qquad + \sum_{j \in \calj_n(\alpha)^c} S_{n,j-1} V_{(s_n(\alpha))}^2 K_{n,j} 1_{\{j \in \mathfrak{L}_n^c \cap \call_n\}} \\
	&{\colorb=:} A_1 + A_2 + A_3. 
\end{align*}
Since $S_{n,j-1} V_{(s_n(\alpha))}^2 K_{n,j} \leq |\Delta_j X|^2 K_{n,j} \leq 
n^{-1/2}$ for $j \in \calj_n(\alpha)^c$, we have 
$\big\|  A_2 \big\|_p \lesssim n^{-1/2+\xi}$. 
As for $A_3$, note that $\# \call_n \lesssim n^{\xi} \times \bar{\kappa}_n 
{\xred=O( n^{\xi + 1/2})}$ 
and that $\Delta_j X = \Delta_j \tilde{X}$ for $j \in \mathfrak{L}_n^c$. 
Hence we have  
$$
\| A_3 \|_p \leq 
\left \| \max_{j \in I_n} | \Delta_j \tilde{X}|^2 \# (\mathfrak{L}_n^c \cap \call_n) \right \|_p 
\lesssim n^{-1/2 + \xi + \ep},
$$
where $\ep$ is an arbitrarily small positive number. 

As for $A_1$, we can set $c = 1$ in the condition  {\colorb $[G2]$}(ii) without loss of generality. 
\begin{align*}
	\left\| A_1 - \alpha c(\alpha) \Theta \right\|_p
	&\leq 
	\left\| \Big( {\xred h^{-1}} V_{(s_n(\alpha))}^2 - c(\alpha) \Big) 
	{\xred h} \sum_{j \in \calj_n(\alpha)^c} S_{n,j-1} K_{n,j} 1_{\{ j \in \call_n^c \}} \right\|_p \\
	&\qquad +
	c(\alpha) \left\| {\xred h} \sum_{j \in \calj_n(\alpha)^c} \big( S_{n,j-1} - \sigma_{t_{j-1}}^2 \big) 1_{\{ j \in \call_n^c \}} \right \|_p\\
	&\qquad + 
	c(\alpha) \left\| {\xred h} \sum_{j \in \calj_n(\alpha)^c} S_{n,j-1} \big( 1 - K_{n,j} \big) 1_{\{ j \in \call_n^c \}} \right \|_p\\
	&\qquad + 
	c(\alpha) \left\| {\xred h}\sum_{j \in \calj_n(\alpha)^c} \sigma_{t_{j-1}}^2 1_{\{ j \in \call_n^c \}} 
	- \alpha \Theta  \right\|_p \\
	&{\colorb=:}  
	B_1 + B_2 + B_3 + B_4. 
\end{align*}

By condition  {\colorb $[G2]$}, $B_2 = O(n^{-\gamma_0})$. 
As for $B_3$, with the estimate $\big\| 1_{\{j \in \call_n^c\}}(1 - K_{n,j}) \big\|_p \leq P[|\Delta_j \tilde{X}| > n^{-1/4}]^{1/p} = O(n^{-L})$ (for all $p > 1$ and $L >0$) and the Cauchy-Schwarz inequality, we have 
$$
\| B_3 \|_p \leq 
{\xred h} \sum_{j \in I_n} \big\| 1_{\{j \in \call_n^c\}} S_{n,j-1} \big\|_{2p} 
\big\| 1_{\{j \in \call_n^c\}} (1 - K_{n,j}) \big\|_{2p} = O(n^{-L}). 
$$

For $B_4$, we use the following decomposition: 
\beas&&
	{\xred h}\sum_{j \in \calj_n(\alpha)^c} \sigma_{\tjm}^2 {\xred 1_{\{ j \in \call_n^c \}}}
	- \alpha \Theta 
	\\&= &
	\left( {\xred h} \sum_{j \in I_n} \sigma_{\tjm}^2 - \Theta \right) 
	- {\xred h} \sum_{j \in I_n} \sigma_{\tjm}^2 1_{\{ j \in \call_n\}} 
	+ {\xred h} \sum_{j \in \calj_n(\alpha)} \sigma_{\tjm}^2 1_{\{ j \in \call_n\}} \\ 
	\\&&\qquad +
	\left( (1-\alpha) \Theta - {\xred h} \sum_{j \in \widetilde{\calj}_n(\alpha)} \sigma_{\tjm}^2 \right) + 
	\left( {\xred h} \sum_{j \in \widetilde{\calj}_n(\alpha)} \sigma_{\tjm}^2 
	 - {\xred h} \sum_{j \in \calj_n(\alpha)} \sigma_{\tjm}^2 \right). 
\eeas
Hence, with the aid of {\colorb Lemmas} \ref{300217-5}, \ref{0211261213} 
and the estimate 
{\xred $\|\call_n\|_p\simleq n^{\xi+1/2}$,}
we have 
\bea\label{Ineq001}&&
	\Bigg\| {\xred h} \sum_{j \in \calj_n(\alpha)^c} \sigma_{\tjm}^2 
	{\xred 1_{\{ j \in \call_n^c \}}}
	- \alpha \Theta \Bigg\|_p
	\nn\\&\lesssim& 
	\Bigg\| {\xred h} \sum_{j \in I_n} \sigma_{\tjm}^2 - \Theta \Bigg\|_p
	 +
	 \Bigg\|  (1-\alpha) \Theta - {\xred h} \sum_{j \in \widetilde{\calj}_n(\alpha)} \sigma_{\tjm}^2  \Bigg\|_p 	+  {\xred O(n^{-\beta_0})}
\eea
{\xred since $\beta_0<\half-\xi$.}
{\colorb 
The first term of the right-hand side of the above inequality is $O(n^{-1/2})$ 
by (\ref{0211280745}). }
As for the second term {\xred on the right-hand side} of (\ref{Ineq001}), 
\begin{align*}
	 \Bigg\|  (1-\alpha) \Theta - {\xred h}\sum_{j \in \widetilde{\calj}_n(\alpha)} \sigma_{\tjm}^2  \Bigg\|_p 
	 &= 
	 \Bigg\|  (1-\alpha) \Theta - {\xred h} \sum_{j \in I_n} \sigma_{\tjm}^2 
	 1_{\{ |W_j| \leq c(\alpha)^{1/2} \}}  \Bigg\|_p \\ 
	 &\leq
	 \Bigg\| {\xred h} \sum_{j \in I_n} {\xred\sigma_{\tjm}^2 }
	 \Big( 1_{\{ |W_j| \leq c(\alpha)^{1/2}\}} - 
	 P\Big[ |W_j| \leq c(\alpha)^{1/2} \Big] \Big) \Bigg\|_p \\
	 &\qquad
	 + (1-\alpha) \Bigg\| {\xred h} \sum_{j \in I_n} \sigma_{t_{j-1}}^2  - \Theta \Bigg\|_p \\ 
	 &= {\colorb  O(n^{-1/2}) }.
\end{align*}
since 
Hence we have ${\colorb B_4} = O(n^{-\frac{1}{2} + \xi})$. 
%

%
{\xred 
Finally, for $B_1$, it suffices to show that
\bea\label{0302141014}
P \left[ \big| {\xred h^{-1}} V_{(s_n(\alpha))} - c(\alpha)^{1/2} \big| >  n^{-\beta_0} \right] 
&=& 
O(n^{-L})
\eea
as $n\to\infty$ for every $L>0$ and 
for every $\beta_0 < \min \{ \gamma_0,  \frac{1}{2} - \xi \}$. 
Let 
\beas 
A_{n,j} &=& \big\{ h^{-1/2} V_j < c(\alpha)^{1/2} - n^{-\beta_0}  \big\}
\\
\calv_{n,j}
&=& 
1_{\big\{|W_j|\leq c(\alpha)^{1/2}-n^{-\beta_0}
+2^{-1}n^{-\gamma_1}
\big\}}
\eeas
and 
\beas 
\mu_n 
&=& 
(1-\alpha)n - 1-n^{\half+\xi+\ep} -nE[\calv_{n,j}]
\eeas
for $\ep>0$. 
Then
\beas&&
P \left[ h^{-1/2} V_{(s_n(\alpha))} - c(\alpha)^{1/2} <  - n^{-\beta_0} \right] 
\\&\leq &
P \left[ \sum_{j \in I_n} 1_{A_{n,j}} \geq (1-\alpha)n - 1 \right] 
\\&=&
P \left[ \sum_{j \in I_n} 1_{A_{n,j}\cap\{j\in\call_n^c\}}
+\sum_{j \in I_n} 1_{A_{n,j}\cap\{j\in\call_n\}} \geq (1-\alpha)n - 1 \right] 
\\&\leq&
P \left[ \sum_{j \in I_n} 1_{A_{n,j}\cap\{j\in\call_n^c\}}
+\#\call_n \geq (1-\alpha)n - 1 \right] 
\\&\leq&
P \left[ \sum_{j \in I_n} \calv_{n,j}\geq (1-\alpha)n - 1-n^{\half+\xi+\ep}  \right]
+P[\#\call_n >n^{\half+\xi+\ep}]
 +P[\Omega_n^c]
 \\&\leq&
P \left[ \sum_{j \in I_n} \big(\calv_{n,j}-E[\calv_{n,j}]\big)
\geq \mu_n \right]
+P[\#\call_n >n^{\half+\xi+\ep}]
 +P[\Omega_n^c].
\eeas
We see
\beas 
\mu_n 
&\sim& 
(1-\alpha)n - 1-n^{\half+\xi+\ep} 
-n\big\{(1-\alpha)-c^*(n^{-\beta_0}-2^{-1}n^{-\gamma_1}\big)\big\}
\>\geq\>
\half c^*n^{1-\beta_0}
\eeas
for large $n$, where $c^*$ is some positive constant, 
if we take a sufficiently small $\ep$ and $\gamma_1\in(\beta_0,\gamma_0)$ thanks to 
$\beta_0 {\xred<} \min \big\{ \gamma_0, \frac{1}{2} - \xi \big\}$. 
%
Since $n^{-1/2}\mu_n\geq 2^{-1}c^*n^{\half-\beta_0}$, from $\beta_0<1/2$, we obtain 
\beas 
P \left[ n^{-1/2}\sum_{j \in I_n} \big(\calv_{n,j}-E[\calv_{n,j}]\big)
\geq  n^{-1/2}\mu_n \right]
&=&
O(n^{-L})
\eeas
for every $L>0$. Therefore,  
\beas 
P \left[ h^{-1/2} V_{(s_n(\alpha))} - c(\alpha)^{1/2} <  - n^{-\beta_0} \right] 
&=& 
O(n^{-L})
\eeas
as $n\to\infty$ for every $L>0$. 
Similarly, we can obtain the estimate 
$P \left[ h^{-1/2} V_{(s_n(\alpha))} - c(\alpha)^{1/2} > n^{-\beta_0} \right] =O(n^{-L})$ 
to show (\ref{0302141014}), 
which concludes the proof. 
}
\qed\halflineskip

\begin{en-text}
\noindent
Let $A_j = \Big\{  {\colorb n^{1/2}} V_j < c(\alpha)^{1/2} - n^{-\beta_0}  \Big\}$. 
Then we have 
\begin{align*}
	P \left[ n^{1/2} V_{(s_n(\alpha))} - c(\alpha)^{1/2} < {\colorb -} n^{-\beta_0} \right] 
	&\leq 
	P \left[ \sum_{j \in I_n} 1_{A_j} \geq (1-\alpha)n - 1 \right] \\
	&=
	P \left[ n^{-1/2} \sum_{j \in I_n} \Big( 1_{A_j} - P[A_j] \Big) \geq C_n \right],
\end{align*}
where 
$$
C_n = n^{-1/2} \left( (1-\alpha)n - 1 - \sum_{j \in I_n} P[A_j] \right).
$$

\noindent
We can estimate $\sum_{j \in I_n} P[A_j] $ as follows: 
\begin{align*}
	\sum_{j \in I_n} P[A_j] 
	&\leq 
	\sum_{j \in I_n} P \Big[ n^{1/2} \tilde{V}_j < c(\alpha)^{1/2} - n^{-\beta_0}, \ j \in  \call_n^c \Big] +E[ \# \call_n]\\
	&\lesssim
	\sum_{j \in I_n} P \Big[ n^{1/2} \tilde{V}_j < c(\alpha)^{1/2} - n^{-\beta_0}, \ j \in  \call_n^c \Big] + n^{\xi + \frac{1}{2}},
\end{align*}
where $\tilde{V}_j = S_{n,j-1}^{-1/2} |\Delta_j \tilde{X}|$. 
The probability on the right-hand side can be estimated as follows:
\begin{align*}
	&P \Big[ n^{1/2} \tilde{V}_j < c(\alpha)^{1/2} - n^{-\beta_0}, \ j \in  \call_n^c \Big] \\
	&\leq
	P \Bigg[ {\colorb |W_j|} < c(\alpha)^{1/2} - n^{-\beta_0}  
	+ n^{1/2} S_{n, j-1}^{-1/2} \bigg| \int_{t_{j-1}}^{t_j} b_t dt \bigg| 
	+ n^{1/2} S_{n, j-1}^{-1/2} \bigg| \int_{t_{j-1}}^{t_j} (\tilde{\sigma}_t - \tilde{\sigma}_{t_{j-1}}) dw_t \bigg| \\
	&\qquad \qquad \qquad 
	+ {\colorb |W_j|} {\colorb S_{n, j-1}^{-1/2}} \Big| {\colorb S_{n, j-1}^{1/2}} - \sigma_{t_{j-1}} \Big|,   
	\ j \in  \call_n^c  \Bigg] \\ 
	&\leq
	P \Bigg[ {\colorb |W_j|}< c(\alpha)^{1/2} - \frac{n^{-\beta_0}}{2}  \Bigg] + 
	P \Bigg[ n^{1/2} S_{n, j-1}^{-1/2} \bigg| \int_{t_{j-1}}^{t_j} b_t dt \bigg|  
		> \frac{n^{-\beta_0}}{6}, \ j \in \call_n^c \Bigg] \\
	&\qquad \qquad
	+ P \Bigg[  n^{1/2} S_{n, j-1}^{-1/2} \bigg| \int_{t_{j-1}}^{t_j} (\tilde{\sigma}_t - \tilde{\sigma}_{t_{j-1}}) dw_t \bigg| 
		> \frac{n^{-\beta_0}}{6}, \ j \in \call_n^c \Bigg] \\
	&\qquad \qquad
	+ P \Bigg[ {\colorb |W_j|}{\colorb S_{n, j-1}^{-1/2}} \Big| {\colorb S_{n, j-1}^{1/2}}- \sigma_{t_{j-1}} \Big|
		> \frac{n^{-\beta_0}}{6}, \ j \in \call_n^c \Bigg]\\ 
	&\lesssim
	P \Bigg[ {\colorb |W_j|} < c(\alpha)^{1/2} - \frac{1}{2} n^{-\beta_0} \Bigg] 
	+ O(n^{-p(\gamma_0 - \beta_0)})\\
	&\sim (1-\alpha) - \phi(c(\alpha)^{1/2}; 0, 1) n^{-\beta_0}
	+ O(n^{-p(\gamma_0 - \beta_0)}), \text{ for all } p > 1. 
\end{align*}
Here we used the following inequality{\colorb :} $C>0$,
\begin{align*}
	&P \Bigg[ {\colorb |W_j|}{\colorb S_{n, j-1}^{-1/2}} \Big| {\colorb S_{n, j-1}^{1/2}} - \sigma_{t_{j-1}} \Big|  
	> Cn^{-\beta_0}, \ j \in  \call_n^c  \Bigg] \\
&\!\!\leq 
	P \Bigg[ {\colorb |W_j|}{\colorb S_{n, j-1}^{-1/2}} n^{\gamma_0} \Big| {\colorb S_{n, j-1}^{1/2}} - \sigma_{t_{j-1}} \Big| 1_{\{  j \in  \call_n^c \}} 
	> Cn^{\gamma_0 - \beta_0}  \Bigg] \\
	&{\colorb \lesssim
	n^{-p(\gamma_0 - \beta_0)} 
	\Bigg\| W_j S_{n,j-1}^{-1/2} \Big( S_{n,j-1}^{1/2} + \sigma_{t_{j-1}} \Big)^{-1}  \times n^{\gamma_0} \Big( S_{n, j-1}- \sigma_{t_{j-1}}^2 \Big) 1_{\{  j \in  \call_n^c \}} \Bigg\|_p^p} \\
	&{\colorb \leq
	n^{-p(\gamma_0 - \beta_0)} 
	\Bigg\| W_j S_{n,j-1}^{-1/2} \Big( S_{n,j-1}^{1/2} + \sigma_{t_{j-1}} \Big)^{-1}  \Bigg\|_{2p}^p 
	\Bigg\| n^{\gamma_0} \Big( S_{n, j-1}- \sigma_{t_{j-1}}^2 \Big) 1_{\{  j \in  \call_n^c \}} \Bigg\|_{2p}^p} \\
	&{\colorb \lesssim
	n^{-p(\gamma_0 - \beta_0)} 
	\Bigg\| W_j S_{n,j-1}^{-1/2} \Big( S_{n,j-1}^{-1/2} + \sigma_{t_{j-1}}^{-1} \Big)  \Bigg\|_{2p}^p 
	\Bigg\| n^{\gamma_0} \Big( S_{n, j-1}- \sigma_{t_{j-1}}^2 \Big) 1_{\{  j \in  \call_n^c \}} \Bigg\|_{2p}^p} \\
	&{\colorb \lesssim
	n^{-p(\gamma_0 - \beta_0)} 
	\Big\| n^{\gamma_0} \Big( S_{n, j-1}- \sigma_{t_{j-1}}^2 \Big) 1_{\{  j \in  \call_n^c \}} \Big\|_{2p}^p = O(n^{-p(\gamma_0 - \beta_0)}) }, 
\end{align*}
{\colorb with the aid of the assumptions $[G1]$(iii), $[G2]$, and the inequality $(x + y)^{-1} \leq \frac{1}{4}(x^{-1} + y^{-1})$ for 
$x, y > 0$ by the convexity of $f(x)=x^{-1}$. }

\noindent
For $\beta_0 < \min \Big\{ \gamma_0, \ \frac{1}{2} - \xi \Big\}$, 
the order of the lower bound $C_n$ is given by 
{\colorb  
\begin{align*}
C_n &\gtrsim 
n^{-1/2} \Bigg\{ (1-\alpha)n - 1 - 
\sum_{j \in I_n} \Bigg( (1-\alpha) - \phi(c(\alpha)^{1/2};0, 1) n^{-\beta_0} + O(n^{-p(\gamma_0 - \beta_0)}) \Bigg) + O(n^{\xi + \frac{1}{2} })   \Bigg\} \\
&= n^{-1/2} \Bigg\{ - 1  
+ \phi(c(\alpha)^{1/2};0, 1) n^{1-\beta_0} + O(n^{1-p(\gamma_0 - \beta_0)}) + 
O(n^{\xi + \frac{1}{2} }) \Bigg\} \\
&\gtrsim
n^{\frac{1}{2} - \beta_0} + n^{\xi} \gtrsim n^{\frac{1}{2} - \beta_0}. 
%
%
%
%
\end{align*}
}
Since the random variable $n^{-1/2} \sum_{j \in I_n} (1_{A_j} - P[A_j])$ is 
bounded in $L^p$ for all $p > 1$, we have
\begin{align*}
	P \left[ n^{1/2} V_{(s_n(\alpha))} - c(\alpha)^{1/2} < {\colorb -} n^{-\beta_0} \right] 
	&\leq 
	P \left[ n^{-1/2} \sum_{j \in I_n} \Big( 1_{A_j} - P[A_j] \Big) \gtrsim 
	n^{\frac{1}{2}-\beta_0} \right] = O(n^{-L}).
\end{align*}

\noindent
As the other inequality $P \left[ n^{1/2} V_{(s_n(\alpha))} - c(\alpha)^{1/2} > n^{-\beta_0} \right] = O(n^{-L})$ can be shown similarly, 
the proof is completed. 
\qed\halflineskip
\end{en-text}
%

\section{Asymptotic mixed normality of the global realized volatilities 
with a moving threshold}\label{0211181815}
%
\subsection{{\sred The GRV with a moving threshold}}
In this section, we will consider a situation where the intensity of jumps 
is moderate. Then it is possible to keep the cut-off ratio of the data small, and 
to get a precise estimate for the integrated volatility. 
\begin{en-text}
{\xred 
The ordinary realized volatility is useless when 
the observations $(X_\tj)_{j\in I_n}$, $I_n=\{1,...,n\}$, are contaminated with jumps. 
To remove the influence of jumps, we will consider a realized volatility estimator with 
a global jump filter.
}
\end{en-text}
Let 
\bea\label{0301030324}
\delta_0\in\big(0,1/4\big)\quad\text{and}\quad\delta_1\in\big(0,1/2\big).
\eea
%
%
In the context of the global jump filtering, 
given a collection $(\mfS_{n,j-1})_{j\in I_n}$ 
of {\xred nonnegative} 
random variables, we {\xred consider} the index set {\xred$\calj_n$} given by 
\bea\label{0301030241}
{\xred\calj}_n &=& 
\big\{j\in I_n;\> V_j<V_{(s_n)}\big\}
\eea
where
\bea\label{0301030242}
V_j &=& \big|{\xred(\mfS_{n,j-1}^-)^{1/2}}
\>\Delta_jX\big|
\eea
and 
\bea\label{0302110546}
s_n \yeq n-\lfloor Bn^{\delta_1}\rfloor
\eea
for a positive constant $B$. 
{\xred 
Here $x^-=1_{\{x\not=0\}}x^{-1}$ for $x\in\bbR$. 
\begin{remark}\rm
It is natural to set a spot volatility estimator of $\sigma_\tjm^2$ in $\mfS_{n,j-1}$ 
though not definitively necessary (Remark \ref{0302110521}). 
In Section \ref{Sec3}, we discussed some constructions of $\mfS_{n,j-1}$. 
In the terminology of Section \ref{Sec2}, the cut-off rate by $\calj_n$ is 
$\alpha_n=\lfloor Bn^{\delta_1}\rfloor/n$, 
$\calj_n=\calj_n(\alpha_n)$ 
and $\alpha_n$ goes to $0$ as $n$ tends to $\infty$. 
We note that the definition of $V_j$ is different from that in (\ref{0211281531}). 
\end{remark}
}
%

%

%
For estimation of $\Theta$ of (\ref{0211281514}), we consider 
the global realized volatility {\fred(GRV)} with a moving threshold 
\bea\label{0302101836}
{\xred{\bf G}}_n 
&=& 
\sum_{j\in{\xred\calj}_n}q_n^{-1}\big|\Delta_jX\big|^2H_{n,j}
\eea
where $(q_n)_{n\in\bbN}$ is a sequence of positive numbers, 
and 
\bea\label{0301030351} 
H_{n,j}&=&1_{\{|\Delta_jX|<B_0n^{-\frac{1}{4}-\delta_0}\}}
\eea
for a positive constant $B_0$. 
%

{\xred 
Here, $\sigma=(\sigma_t)_{t\in[0,T]}$ and 
$b=(b_t)_{t\in[0,T]}$ are \cadlag adapted processes. 
We will assume (\ref{0301030324}) and the following conditions. 
\begin{en-text}
{\sred 
\bd
\im[{\bf [G1]}] 
{\bf (i)} For every $p>1$, $\sup_{t\in[0,T]}\|\sigma_t\|_p<\infty$ and 
\beas 
\big\|\wt{\sigma}_t-\wt{\sigma}_s\big\|_p &\leq& C(p)|t-s|^{1/2}
\quad(t,s\in[0,T])
\eeas
for some constant $C(p)$ for every $p>1$. 
\im[{\bf (ii)}] $\sup_{t\in[0,T]}\|b_t\|_p<\infty$ for every $p>1$. 
\im[{\bf (iii)}] $\sigma_t\not=0$ a.s. for every $t\in[0,T]$, an 
$\sup_{t\in[0,T]}\big\|\sigma_t^{-1}\big\|_p<\infty$ for every $p>1$. 
%
\ed
\end{en-text}
{\xred 
\bd
\im[{\bf [G1$^o$]}] 
For every $p>1$, $\sup_{t\in[0,T]}\|\sigma_t\|_p+\|b_t\|_p<\infty$. 
\ed
}

{\xred 
\bd
\im[{\bf [G2$^o$]}] 
$q_n>0$ $(n\in\bbN)$ and 
$q_n-1=o(n^{-1/2})$ as $n\to\infty$. 
\ed
}
\halflineskip

{\xred 
\begin{remark}\label{0302110521}\rm 
Theoretically, we may set $\mfS_{n,j-1}=1$. 
Condition ${\xred[G2^o]}$ is satisfied with $q_n=1$. 
Asymptotically the choice $(\mfS_{n,j-1},q_n)=(1,1)$ is sufficient and valid. 
However, in practice, a natural choice is a pair 
$S_{n,j-1}$ {\xred satisfying $[G2]$ in Section \ref{Sec2}} and $q_n=q(\alpha_n)$, 
{\xred where the function $q$ is defined by (\ref{0302110914}) in Section \ref{Sec2}.}
\end{remark}
}

{\xred 
\begin{en-text}
The jump part $J$ of $X$ is assumed to be finitely active, i.e., 
$N_T<\infty$ a.s., $N_t=\sum_{s\leq t}1_{\{\Delta J_s\not=0\}}$. 
We assume that the distribution of the variable $N_T$ 
depends on $n$ 
and $N_T$ 
may diverge as $n\to\infty$. 
More precisely, the following condition is assumed 
for $\Lambda_n=\#\{j\in I_n;\>$\koko
\end{en-text}
For the jump part $J$ of the semimartingale $X$, we only assume 
\beas 
\Lambda_n &:=&\#\{j\in I_n;\>\Delta_jJ\not=0\}<\infty\quad a.s. 
\eeas
for every $n\in\bbN$, 
and the following estimate:
\begin{en-text}
{\xred The following condition models a fairly high jump-intensity of $\ol{N}$, 
in particular, exogenous heavy contamination with jumps. }
\bd
\im[{\bf [G3]}]
There exists a constant $\xi\geq0$ 
such that 
$\|\ol{N}_T\|_p=O(n^{\xi})$ as $n\to\infty$ 
for every $p>1$. 
\ed
\end{en-text}
}
{\xred 
\bd
\im[{\bf [G3$^o$]}]
There exists a constant $\xi\geq0$ 
such that 
$\|{\xred \Lambda_n}\|_p=O(n^{\xi})$ 
as $n\to\infty$ 
for every $p>1$. 
\ed
\begin{remark}\rm
The diverging $\Lambda_n$ models high intensity of the jump part for a fixed $n$ 
in practice. 
Mathematically, we are assuming that the process $\sigma$ is independent of $n$. 
This makes sense naturally in particular when the jumps are exogenous. 
It is sufficient for the limit theorem by using the \cadlag property of $\sigma$. 
On the other hand, though details are omitted, we can treat $\sigma$ depending on $n$ 
if 
uniform $L^\inftym$-continuity of $\sigma$ 
and uniformity in $[G1^o]$ are satisfied. 
\end{remark}
}

{\xred 
Define $\Gamma$ by 
\beas
\Gamma &=& 2T\int_0^T \sigma_t^4dt.
\eeas
Extend $(\Omega,\calf,P)$ so that there is 
a standard normal random variable $\zeta$ independent of $\calf$ 
on the extension. 
The $\calf$-stable convergence is denoted by $\to^{d_s}$. 
We obtain asymptotic mixed normality of the global realized volatility ${\bf G}_n$ with 
a moving threshold. 
\begin{theorem}\label{0302110536}
Suppose that $[G1^o]$, $[G2^o]$ and $[G3^o]$ are satisfied. 
Suppose that $\xi<2\delta_0$. Then 
\beas 
n^{1/2}\big(
{\xred {\bf G}_n}
-\Theta\big)
&\to^{d_s}&
\Gamma^{1/2}\zeta
\eeas
as $n\to\infty$. 
\end{theorem}
{\xred 
Theorem \ref{0302110536} follows from Theorem \ref{0211291415}, 
that is presented in a slightly more general setting. 
}
\begin{en-text}
{\color{gray}
Just combine Lemmas 
\ref{0211281630}, \ref{0211281845} and \ref{0211291358}. 
}
\end{en-text}

\subsection{The WGRV with a moving threshold}
%
Suppose that a collection $(\mfS_{n,j-1})_{j\in I_n}$ ($n\in\bbN$) 
of positive random variables is given. 
Consider constants $\delta_0$ and $\delta_1$ satisfying (\ref{0301030324}), 
and 
$V_j$ and $s_n$ given by 
(\ref{0301030242}) 
and (\ref{0302110546}), respectively. 
We define 
{\fred 
the Winsorized global realized volatility (WGRV) with a moving threshold by
}
\bea\label{0302111147}
{\bf W}_n 
&=& 
\sum_{j \in I_n} q_n^{-1} 
\big\{ |\Delta_j X| \wedge \mathfrak{S}_{n,j-1}^{1/2} V_{(s_n)} \big) \big\}^2 H_{n,j}, 
\eea
where $(q_n)_{n \in \mathbb{N}}$ is a sequence of positive numbers. 
The error of the WGRV has the same limit as GRV ${\bf G}_n$.
}}


\begin{theorem}
	Suppose that {\xred$[G1^o]$}, {\xred$[G2^o]$} and {\xred$[G3^o]$} are satisfied. 
	Suppose that $\xi < 2\delta_0$. 
	Then 
	\beas 
	n^{1/2}\big(
	{\bf W}_n-\Theta\big)
	&\to^{d_s}&
	\Gamma^{1/2}\zeta
	\eeas
	as $n\to\infty$
where $\zeta$ is a standard Gaussian random variable independent of $\calf$.
\end{theorem}
\proof 
{\xred Let $\wt{X}=X-J$.}
It suffices to show that 
$n^{1/2} \| {\bf W}_n - {\xred{\bf G}}_n  \|_p \to 0$
as $n\to\infty$ for every $p > 1$. 
{\xred From (\ref{0302111147}), }
\beas
{\bf W}_n - {\xred{\bf G}}_n  
&=& 
\sum_{j \in {\xred\calj}_n^c} q_n^{-1} \mathfrak{S}_{n,j-1} V_{(s_n)}^2 H_{n,j}. 
\eeas
%
We have 
\begin{align*}
	\left\| \sum_{j \in {\xred\calj}_n^c} \mathfrak{S}_{n,j-1} V_{(s_n)}^2 H_{n,j} \right\|_p
	&\leq 
	\left\| \sum_{j \in {\xred\calj}_n^c} |\Delta_j X|^2 1_{\{ j \in \mathfrak{L}_n\}} H_{n,j} \right\|_p 
	+ \left\|  \sum_{j \in {\xred\calj}_n^c} |\Delta_j \tilde{X}|^2 1_{\{ j \in \mathfrak{L}_n^c\}} H_{n,j} \right\|_p  \\
	&\lesssim n^{-1/2-2\delta_0+\xi} + n^{-1+\ep+\delta_1}
\end{align*}
{\xred as $n\to\infty$ for any $\ep>0$.}
Since $\delta_1 < 1/2$, $\xi < 2\delta_0$ 
{\xred and $\ep$ is arbitrary}, we obtain the desired 
{\xred convergence from Theorem \ref{0302110536}.}
\qed\halflineskip

\subsection{Stability of the realized volatility under missing}

{\sred 
We are about establishing asymptotic mixed normality of 
the integrated volatility estimator having a moving threshold. 
We will solve this problem by showing a stability of estimation under elimination of a certain portion of the data. 
{\xred In other words, this is a question of stability under missing data. 
In what follows, we will consider the variable ${\bf V}_n$ defined by 
\bea\label{0301030318}
{\bf V}_n 
&=& 
\sum_{j\in\calm_n}q_n^{-1}\big|\Delta_jX\big|^2H_{n,j}
\eea
where $(q_n)_{n\in\bbN}$ is a sequence of positive numbers,
$H_{n,j}$ is given in (\ref{0301030351}), and 
$\calm_n$ is an abstract random index set in $I_n$.} 
It is not necessary to specify {\xred$\calm_n$ like $\calj_n$}  
by (\ref{0301030241}) and (\ref{0301030242}). 
\begin{en-text}
\bd
\im[{\bf [G2$'$]}]
\bd
\im[{\bf (i)}] 
The positive random variables $\mfS_{n,j-1}$ $(j\in I_n, n\in\bbN)$ satisfy 
\beas 
\sup_{n\in\bbN}\sup_{j\in I_n}\big(\|\mfS_{n,j-1}\|_p+\|\mfS_{n,j-1}^{-1}\|_p\big)
&<& 
\infty
\eeas
for every $p>1$. 

\im[{\bf (ii)}] 
$|q_n-1|=o(n^{-1/2})$ as $n\to\infty$. 
\ed
\ed
\end{en-text}
\begin{en-text}
{\colorg The following green part is not necessary. Delete it later. 
Define two events $\mfA_n$ and $\mfB_n$ by 
\beas 
\mfA_n 
&=& 
\bigcup_{j\in I_n}
\bigg[\big\{j\in\calm_n^c\big\}\cap\big\{\Delta_jN=0\big\}\bigg]
\eeas
and 
\beas
\mfB_n
&=& 
\bigcap_{j\in I_n}\bigg[\big\{V_j\geq V_{(s_n)}\big\}\cup\big\{|\Delta_jJ|\leq n^{-\frac{1}{4}-\delta_0}\big\}\bigg]
\eeas

\begin{lemma}\label{0211281550}
For any $L>0$, the following estimates hold as $n\to\infty$:
\bd
\im[{\bf (i)}] 
$P\big[\mfA_n\cap\mfB_n^{\>c}\big]\yeq O(n^{-L})$
\im[{\bf (ii)}] 
$P\big[\mfA_n^{\>c}\big]\yeq O(n^{-L})$
\im[{\bf (iii)}] 
$P\big[\mfB_n\big]\yeq 1-O(n^{-L})$.
\ed
\end{lemma}
\proof \koko
}
\end{en-text}

Let 
\beas 
{\bf V}_n^\dagger
&=& 
\sum_{j\in\calm_n}q_n^{-1}\big|\Delta_j\wt{X}\big|^2H_{n,j}. 
\eeas
Recall $\wt{X}=X-J$.

\bd
\im[{\bf [G2$'$]}]
\bd
\im[{\bf (i)}] 
For every $n\in\bbN$, $\calm_n$ is a random set in $I_n$ such that 
$\#\big(I_n\setminus\calm_n\big)\leq B_1n^{\delta_1}$ $(n\in\bbN)$ for some positive constant $B_1$. 
\im[{\bf (ii)}] $q_n>0$ $(n\in\bbN)$ and 
$q_n-1=o(n^{-1/2})$ as $n\to\infty$. 
\ed
\ed
\halflineskip

\begin{lemma}\label{0211281630}
Suppose that {\xred$[G1^o]$}, $[G2']$ and {\xred$[G3^o]$} are satisfied. 
Suppose that $\xi<2\delta_0$. Then 
\beas 
n^{1/2}\big\|{\bf V}_n-{\bf V}_n^\dagger\big\|_p
&\to& 
0
\eeas
as $n\to\infty$ for every $p>1$. 
\end{lemma}
\proof We have the estimate
\bea\label{0211281840}
n^{1/2}\big\|{\bf V}_n-{\bf V}_n^\dagger\big\|_p
&\leq&
2q_n^{-1}\Phi_n^{(\ref{0211281652})}+q_n^{-1}\Phi_n^{(\ref{0211281653})},
\eea
where
\bea\label{0211281652}
\Phi_n^{(\ref{0211281652})}
&=& 
n^{1/2}\bigg\|
\sum_{j\in\calm_n}\big|\Delta_j\wt{X}\Delta_jJ\big|
H_{n,j}
\bigg\|_p,
\eea
and
\bea\label{0211281653}
\Phi_n^{(\ref{0211281653})}
&=& 
n^{1/2}\bigg\|
\sum_{j\in\calm_n}\big|\Delta_jJ\big|^2
H_{n,j}
\bigg\|_p
\eea
for $p>1$. 
By using the inequality
\beas 
|\Delta_jJ|H_{n,j}&\leq&\big(|\Delta_j\wt{X}|
+{\sred B_0\>}
n^{-\frac{1}{4}-\delta_0}\big)1_{\{\Delta_jJ\not=0\}},
\eeas
we obtain
\begin{en-text}
\bea\label{0211281656}
\Phi_n^{(\ref{0211281652})}
&\leq&
n^{1/2}\bigg\|
\sum_{j\in\calm_n}\big|\Delta_j\wt{X}\Delta_jJ\big|
H_{n,j}1_{\mfB_n}
\bigg\|_p
+
n^{1/2}\bigg\|
\sum_{j\in\calm_n}\big|\Delta_j\wt{X}\Delta_jJ\big|
H_{n,j}1_{\mfB_n^c}
\bigg\|_p
\nn\\&\leq&
n^{1/2}\bigg\|
\max_{j\in I_n}\big\{\big|\Delta_j\wt{X}\big|(|\Delta_j\wt{X}|
+{\sred B_0\>}n^{-\frac{1}{4}-\delta_0}\big)\big\}
N_T
\bigg\|_p
+
n^{1/2}\bigg\|
\max_{j\in\calm_n}\big|\Delta_j\wt{X}\big|
\sum_{j\in\calm_n}\big|\Delta_jJ\big|
H_{n,j}1_{\mfB_n^c}
\bigg\|_p
\nn\\&\simleq&
n^{-2\delta_0}\big\|N_T\big\|_p
+
\eea
\end{en-text}
\beas
\Phi_n^{(\ref{0211281652})}
&\leq&
n^{1/2}\bigg\|
\max_{j\in I_n}\big\{\big|\Delta_j\wt{X}\big|(|\Delta_j\wt{X}|
+{\sred B_0\>}n^{-\frac{1}{4}-\delta_0}\big)\big\}
\bigg\|_{2p}
\big\|
{\xred\Lambda_n}
\big\|_{2p}
\nn\\&\simleq&
n^{-\frac{1}{4}-\delta_0+\xi+\ep}
\eeas
as $n\to\infty$ for any $\ep>0$ and $p>1$. 
Therefore, 
\bea\label{0211281656}
\Phi_n^{(\ref{0211281652})}
&\to&
0
\eea
for every $p>1$ 
since $\xi<2\delta_0<\frac{1}{4}+\delta_0$. 
Similarly, 
\beas
\Phi_n^{(\ref{0211281653})}
&\simleq&
n^{1/2}\bigg\|
\max_{j\in I_n}(|\Delta_j\wt{X}|^2+n^{-\frac{1}{2}-2\delta_0}\big)
\bigg\|_{2p}
\big\|
{\xred\Lambda_n}
\big\|_{2p}
\nn\\&\simleq&
n^{-2\delta_0+\xi+\ep}
\eeas
as $n\to\infty$ for any $\ep>0$ and $p>1$ since 
{\xred$\delta_0<1/4$}. 
In particular, 
\bea\label{0211281834}
\Phi_n^{(\ref{0211281653})}
&\to&
0
\eea
as $n\to\infty$ since $\xi<2\delta_0$. 
Now the proof is completed with 
(\ref{0211281840}), (\ref{0211281656}) and (\ref{0211281834}). 
\qed\halflineskip

Define $\wt{\bf V}_n$ by 
\beas 
\wt{\bf V}_n
&=& 
\sum_{j\in I_n}\big|\Delta_j\wt{X}\big|^2{\sred.}
\eeas
\begin{lemma}\label{0211281845}
Suppose that $\xi<1/2$. Then
\beas 
n^{1/2}\big\|{\bf V}_n^\dagger-\wt{\bf V}_n\big\|_p
&\to&
0
\eeas
as $n\to\infty$ for every $p>1$. 
\end{lemma}
\proof 
{\sred 
Recall that $\delta_0<1/4$ and $\delta_1<1/2$. 
Define ${\bf V}_n^\ddagger$ by 
\beas 
{\bf V}_n^\ddagger
&=& 
\sum_{j\in\calm_n}q_n^{-1}\big|\Delta_j\wt{X}\big|^2.
\eeas
Then 
\beas
n^{1/2}\big\|{\bf V}_n^\dagger-{\bf V}_n^\ddagger\big\|_p
&\simleq&
n^{1/2}\bigg\|\sum_{j\in\calm_n}\big|\Delta_j\wt{X}\big|^2|H_{n,j}-1|\bigg\|_p
\nn\\&\leq&
n^{1/2}\bigg\|\sum_{j\in\calm_n}\big|\Delta_j\wt{X}\big|^2|H_{n,j}-1|
1_{\{{\xred\Delta_jJ\not=0}\}}\bigg\|_p
\nn\\&&
+n^{1/2}\bigg\|\sum_{j\in\calm_n}\big|\Delta_j\wt{X}\big|^2
1_{\{|\Delta_j\wt{X}|>n^{-\frac{1}{4}-\delta_0}\}}1_{\{{\xred\Delta_jJ=0}\}}\bigg\|_p
\nn\\&\leq&
n^{1/2}\bigg\|\max_{j\in I_n}\big|\Delta_j\wt{X}\big|^2{\xred\Lambda_n}\bigg\|_p
+O(n^{-L})
\nn\\&\simleq&
n^{-\half+\ep+\xi}
\eeas
for any positive number $\ep>0$. 
Here $L$ is an arbitrary positive number greater than $1/2$, and 
we used the inequality $\delta_0<1/4$ to get $O(n^{-L})$. 
Since $\xi<1/2$, we obtain 
\bea\label{0301021857}
n^{1/2}\big\|{\bf V}_n^\dagger-{\bf V}_n^\ddagger\big\|_p
&=& 
o(1)
\eea
as $n\to\infty$ for every $p>1$. 
}

{\sred 
From the condition $q_n-1=o(n^{-1/2})$ of $[G2']$ (ii),} obviously, 
\bea\label{0211281849}
n^{1/2}\big\|{\sred{\bf V}_n^\ddagger}-\wt{\bf V}_n\big\|_p
&\leq& 
n^{1/2}\bigg\|\sum_{j\in I_n\setminus\calm_n}\big|\Delta_j\wt{X}\big|^2\bigg\|_p
+o(1)
\nn\\&\simleq&
n^{-\half{\sred + \ep}+\delta_1}+o(1)
\yeq
o(1)
\eea
as $n\to\infty$ for every $p>1$ since 
$\#(I_n\setminus\calm_n)\simleq n^{\delta_1}$ with $\delta_1<1/2$ 
{\sred 
thanks to $[G2']$ (i) and (\ref{0301030324}),}
and 
{\sred 
\beas 
\bigg\|\max_{j\in I_n}\big|\Delta_j\wt{X}\big|^2\bigg\|_p
&=&
O(n^{-1+\ep})
\eeas
for any $p>1$ and any positive number $\ep$. 
Proof ends with (\ref{0301021857}) and (\ref{0211281849}).
}
\qed\halflineskip

\begin{lemma}\label{0211291358}
Suppose that {\xred$[G1^o]$} is satisfied. Then 
\beas 
n^{1/2}\big(\wt{\bf V}_n-\Theta\big)
&\to^{d_s}&
\Gamma^{1/2}\zeta
\eeas
as $n\to\infty$. 
\end{lemma}
\proof
We have
\bea\label{0211291352}
\wt{\bf V}_n
&=&
\sum_{j\in I_n}
\left(\int_\tkm^\tk \sigma_tdw_t+\int_\tkm^\tk b_tdt\right)^2
\nn\\&=&
\Phi^{(\ref{0211280431})}_n+\Phi^{(\ref{0211280432})}_n+
2\Phi^{(\ref{0211280433})}_n
+\Phi^{(\ref{0211280435})}_n
\eea
where 
\bea\label{0211280431}
\Phi^{(\ref{0211280431})}_n
&=& 
\sum_{j\in I_n}
2\int_\tjm^\tj\int_\tjm^t\sigma_sdw_s\sigma_tdw_t{\colorb , }
\eea
\bea\label{0211280432}
\Phi^{(\ref{0211280432})}_n
&=& 
\sum_{j\in I_n}
\int_\tjm^\tj \sigma_t^2dt
{\xred\>\yeq\Theta},
\eea
\bea\label{0211280433}
\Phi^{(\ref{0211280433})}_n
&=& 
\sum_{j\in I_n}
\int_\tjm^\tj \sigma_tdw_t\int_\tjm^\tj b_tdt,
\eea
%
and 
\bea\label{0211280435}
\Phi^{(\ref{0211280435})}_n
&=& 
\sum_{j\in I_n}
\bigg(\int_\tjm^\tj b_tdt\bigg)^2.
\eea

Since $b$ is a \cadlag process, 
for any $\ep>0$, there exists a number $\delta>0$ such that 
$P\big[w'(b,\delta)\geq\ep\big]<\ep$. Here $w'(b,\delta)$ is 
a modulus of continuity defined by 
\beas 
w'(b,\delta)
&=&
\inf_{(s_i)\in\cals_\delta}\max_i\sup_{r_1,r_2\in[s_{i-1},s_i)}\big|b_{r_1}-b_{r_2}\big|,
\eeas
where $\cals_\delta$ is the set of sequences $(s_i)$ such that 
$0=s_0<s_1<\cdots<s_v=T$ and $\min_{i=1,...,v-1}(s_i-s_{i-1})>\delta$. 
Let 
\bea\label{0211290604}
\dot{\Phi}^{(\ref{0211280433})}_n
&=& 
\sum_{j\in I_n}
\int_\tkm^\tk \sigma_tdw_t\int_\tjm^\tj (b_t-b_\tjm)dt{\colorb .}
\eea
Write 
\beas 
{\sf E}_j &=& \int_\tjm^\tj \sigma_tdw_t,
\nn\\
{\sf V}_j &=& n^{1/2}\bigg|\int_\tjm^\tj \sigma_tdw_t \bigg| \int_\tjm^\tj\big(|b_t|+|b_\tjm|\big)dt.
\eeas
For $\omega\in\Omega$ such that 
$w'(b(\omega),\delta)<\ep$, 
there exists a $(s_i)$ (depending on $\omega$) such that 
\beas &&
\max_i\sup_{r_1,r_2\in[s_{i-1},s_i)}\big|b_{r_1}(\omega)-b_{r_2}(\omega)\big|<\ep,
\nn\\&&
\min_{i=1,...,v-1}(s_i-s_{i-1})>\delta. 
\eeas
For $n>T/\delta$, 
all intervals $[\tjm,\tj)$ ($j\in I_n$) includes at most one point among $(s_i)$, 
therefore the number of intervals $[\tjm,\tj)$ that include some one $s_i$ 
is at most $T/\delta$. 
The increment of $b(\omega)$ in $[\tjm,\tj)$ is less than $\ep$ 
if $[\tjm,\tj)\cap\{s_i\}{\xred=\emptyset}$. 
Thus, we have the inequality 
\beas 
\big\|n^{1/2}\dot{\Phi}^{(\ref{0211280433})}_n\big\|_p
&\leq&
\bigg\|\sum_{j\in I_n}n^{1/2}|{\sf E}_j|\bigg\|_p\ep h
+\bigg\|\max_{j\in I_n}{\sf V}_j\bigg\|_p\frac{T}{\delta}
+\bigg\|\sum_{j\in I_n}{\sf V}_j\bigg\|_{2p}P\big[w'(b,\delta)
{\xred\geq}\ep\big]^{\frac{1}{2p}}
\eeas
for every $p>1$. 
Therefore, 
\beas 
\big\|n^{1/2}\dot{\Phi}^{(\ref{0211280433})}_n\big\|_p
&\leq&
C\bigg[
\ep
+\bigg(n^{-1/2}+\sum_{j\in I_n}\big\|{\sf V}_j1_{\{{\sf V}_j>n^{-1/2}\}}\big\|_p\bigg)\frac{T}{\delta}
+\ep^{\frac{1}{2p}}\bigg]
\nn\\&\leq&
C'\big(\ep+n^{-1/2}+\ep^{\frac{1}{2p}}\big)
\eeas
for all $n>T/\delta$, where $C$ and $C'$ are some constants independent of $n$. 
Consequently, 
\bea\label{0211291342}
\lim_{n\to\infty}\big\|n^{1/2}\dot{\Phi}^{(\ref{0211280433})}_n\big\|_p
&=&
0
\eea
for every $p>1$. 
Moreover, for 
\beas
\ddot{\Phi}^{(\ref{0211280433})}_n
&=& 
\sum_{j\in I_n}
\int_\tjm^\tj \sigma_tdw_t\int_\tkm^\tk b_\tjm dt
\yeq 
\sum_{j\in I_n}hb_\tjm 
\int_\tjm^\tj \sigma_tdw_t, 
\eeas
we have 
\bea\label{0211291346}
\lim_{n\to\infty}\big\|n^{1/2}\ddot{\Phi}^{(\ref{0211280433})}_n\big\|_p
&=&
0
\eea
for every $p>1$, by orthogonality. 
From (\ref{0211291342}) and (\ref{0211291346}), 
\bea\label{0211291347}
\lim_{n\to\infty}\big\|n^{1/2}\Phi^{(\ref{0211280433})}_n\big\|_p
&=&
0
\eea
for every $p>1$. 

Obviously, 
\bea\label{0211291355}
\big\|n^{1/2}\Phi^{(\ref{0211280435})}_n\big\|_p
&=&
0
\eea
for every $p>1$. 
Now, we can show the claim of the lemma by using 
(\ref{0211291352}), (\ref{0211291347}) and (\ref{0211291355}) 
together with the mixture type of martingale central limit theorem 
applied to ${\xred n^{1/2}}\Phi^{(\ref{0211280431})}_n$, 
{\xred with the aid of the \cadlag property of $\sigma$.}
\qed\halflineskip

\begin{theorem}\label{0211291415}
Suppose that {\xred$[G1^o]$}, $[G2']$ and {\xred$[G3^o]$} are satisfied. 
Suppose that $\xi<2\delta_0$. Then 
\beas 
n^{1/2}\big(
{\colorb {\bf V}_n}
-\Theta\big)
&\to^{d_s}&
\Gamma^{1/2}\zeta
\eeas
as $n\to\infty$. 
\end{theorem}
\proof 
Just combine Lemmas 
\ref{0211281630}, \ref{0211281845} and \ref{0211291358}. 
\qed\halflineskip


\section{Constant volatility}\label{0211281401}
The case of constant $\sigma$ is specific and theoretical treatments can be slightly different from those of the previous sections. 
In this situation, we do not need to pre-estimate the local spot volatility, 
and hence, we can take $S_{n,j}=1$ constantly and 
no approximation error is caused. 
$\sigma_t=\theta \tilde{\sigma}_t$ is also the case if $\tilde{\sigma}_\tjm$ are  observable. 
For example, the GRV with a fixed cut-off rate $\alpha$ is redefined as 
$$
\bbV_n^0 (\alpha) = \sum_{j \in \calj_n^0(\alpha)} q(\alpha)^{-1} |\Delta_j X|^2 K_{n,j}, 
$$
where  
$$
\calj_n^0(\alpha) = \big\{ j\in I_n;\> |\Delta_j X| < |\Delta X|_{(s_n(\alpha))} \big\}. 
$$
Then we have the following theorem. 
Note that we do not need the condition $[G2]$, 
and $\gamma_0$ in $[G2]$(ii) can be arbitrarily close to $1/2$. 

\begin{theorem}
	Suppose that $[G1]$ and $[G3]$ are fulfilled. 
	Suppose that $\xi < \frac{1}{2}$. Let $\alpha \in (0, 1)$ and 
	$\beta_0 < \frac{1}{2} - \xi$. 
	Then 
	$$
	  \| \bbV_n^0(\alpha) - \Theta\|_p = O(n^{-\beta_0})
	$$
	as $n \to \infty$ for every $p>1$. 
\end{theorem}

The other global-threshold estimators are discussed similarly. 


\section{Simulation studies}\label{Sec7}
In this section, we conduct several numerical simulations to see that 
our global realized volatility estimators outperform those proposed in previous studies. 

\subsection{The case of compound Poisson jumps}\label{SubSec7-1}

Here we consider {\tred a process $X=(X_t)_{t\in[0,1]}$ satisfying 
	the stochastic differential equation}
\begin{equation}
	dX_t = \theta X_t dt + 
	(\sigma + \eta X_t^2 )^{\frac{1}{4}} dw_t + dJ_t, \quad t \in [0, 1]{\colorb ,} 
\end{equation} 
{\tred with $X_0=1$}, 
where $J_t$ is the jump part of $X$. 
{\tred In this section,}
we assume that $J$ is a compound Poisson process of the form 
$
J_t = \sum_{i =1}^{N_t}  \xi_i,
$
where $(N_t)_t$ is a Poisson process with intensity $\lambda > 0$ and 
$(\xi_i)_i$ are independently and normally distributed random variable with 
mean $\mu$ and variance $\nu^2$.
For the intensity parameter, we consider both cases where 
$\lambda$ is high and low. 
Our aim is to estimate the integrated volatility $\Theta = \int_0^1 (\sigma + \eta {\tred X}_t^2)^{\frac{1}{2}} dt$. 

{\tred By simulation, we will}
compare {\tred the performance of} 
the {\colorb threshold realized volatility (TRV),} bipower variation (BV), minimum realized volatility (minRV), the GRV, and the WGRV, 
{\tred where TRV, BV and minRV are given by
\begin{align*}
	\mathrm{TRV}_n &= \sum_{j=1}^n |\Delta_j X|^2 1_{\{ |\Delta_j X| \leq n^{-\rho} \}}, \quad \rho \in (0, 1/2), \\
	\mathrm{BV}_n &= \frac{\pi}{2} \sum_{j=1}^{n-1} |\Delta_j X| |\Delta_{j+1} X|, \\
	\mathrm{minRV}_n &= \frac{\pi}{\pi-2} \sum_{j=1}^{n-1} |\Delta_j X|^2 \wedge |\Delta_{j+1} X|^2, 
\end{align*}
respectively. 
The package YUIMA (cf. \cite{Yuima2014}, \cite{iacus2017simulation}) 
was used for the simulation studies below. 
}

Note that, although TRV is based on threshold method, it is completely 
different from our GRV, since TRV employs a deterministic threshold 
and never uses information of other increments. In this sense, 
TRV is based on a ``local" approach. 

The set-up of simulation is as follows. The number of samples is $n = 2000$. 
We repeat calculating the estimators 500 times to obtain their average and quantile. 
The true parameters are 
$\theta = 0.2, \ \sigma = 1, \ \eta = 3, \ \mu = 0.3, \ \nu = 0.2$. 
Throughout this subsection, we set the cut-off ratio $\alpha = 0.2$ for 
{\colorb GRV and WRGV} with {\tred a local volatility estimator $S_{n,j-1}$}. 
That is, we trim the upper 20\% of absolute increments. 
While it may seem that we eliminate too many observations and 
the estimator suffers from downside bias, GRV and WGRV estimate the integrated volatility 
well thanks to the adjustment coefficient by $q(\alpha)$ and ${\sf w}(\alpha)$.
{\colorb 
In calculating the TRV, we set $\rho = 0.45, 0.2, 0.1$ to see 
the effect of the choice of this parameter on the accuracy of estimation. }
Note that {\tred $\sigma_s$ in (\ref{0301191155})} 
is not directly observable and depends on {\tred $X_t$}. 
Hence,  {we need \tred $S_{n,j-1}$} 
to normalize the increment 
$\Delta_i {\tred X}$ when constructing the GRV. In this simulation, 
we use the LGRV {\tred$\bbL_{n,j}(\alpha_0)$ of} (\ref{0211260215}) 
{\tred with $\alpha_0=0.2$,} 
and {\tred the local minRV $\bbM_{n, j}$ of} (\ref{loc_minRV}) 
{\tred for $S_{n,j-1}$. 
	We adopt $\kappa_n = \lfloor 10 \times n^{0.45} \rfloor$
	for the length of a subinterval to calculate these local volatilities.} 
Moreover, we calculate GRV without normalization (defined in Section \ref{0211281401}) 
for comparison. 
{\colorb
Note that $\kappa_n$ depends on two tuning parameters, the choice of 
which can affect the precision of estimation. 
We argue this point in the final Section \ref{SubSec7-3}. 
}

We use the following labels as in Table \ref{Table01} 
to describe the estimators. 

\renewcommand{\arraystretch}{1.5}
\begin{table}[H]
	\label{Table01}
	\centering
		\caption{Definitions of estimators}
	\begin{tabular}[t]{lllll}
		\hline
		Label & Method & Spot volatility & Cut-off ratio $\alpha$ & Exponent $\rho$ for truncation  \\
		\hline \hline
		{\tt trv[$\rho$]} & TRV & -- & -- & 0.45, 0.2, 0.1 \\		
		{\tt bv} & BV & -- & -- & --  \\
		{\tt mrv}  & minRV & -- & -- & -- \\
		{\tt grv.lgrv[$\alpha$]} & GRV & GRV & 0.2 & --\\
		{\tt grv.mrv[$\alpha$]} & GRV & minRV  & 0.2 & -- \\
		{\tt wgrv.lgrv[$\alpha$]} & WGRV & GRV & 0.2 & -- \\
		{\tt wgrv.mrv[$\alpha$]} & WGRV & minRV & 0.2 & -- \\
		{\tt grv[$\alpha$]} & GRV & -- & 0.2, 0.1, 0.05 & -- \\
		{\tt grv.lgrv.mov} &GRV & GRV & depends on $n$& -- \\		
		{\tt wgrv.lgrv.mov} & WGRV & GRV & depends on $n$ & --  \\				
		\hline
	\end{tabular}
\end{table}
\renewcommand{\arraystretch}{1.0}

\subsubsection{The case of high intensity: GRV with fixed cut-off ratio}
First, we deal with the case of high intensity. 
Here we set $\lambda = 30$ so that the data includes many jumps. 
The example of a sample path and its increments are shown in Figure \ref{fig1}. 
Obviously, there are many large spikes in the data, suggesting the 
existence of jumps. 

Note that the volatility is non-constant here. 
In fact, in Panel (b) of Figure \ref{fig1}, the size of increments tend to increase 
as time passes. 
Hence, to estimate the volatility, we have to use estimated spot volatilities 
to normalize the increments. 

{\tred
In this example, we show the error ratios of GRV and WGRV with shrinking cut-off ratio (tuning parameters that determine the cut-off ratio $\alpha_n = \lfloor B n^{\delta_1} \rfloor$ for these estimators are $B=10$ and $\delta_1=0.45$, 
the same as those used in the next subsection). 
Theoretically, they are available in the case of moderate intensity of jumps. 
We show their results just for reference. 
We will discuss the case of moderate intensity in more detail in the next subsection. 
}

\begin{figure}[H]
	\begin{center}
		\begin{tabular}{c}
			
			\begin{minipage}{0.5\hsize}
				\begin{center}
					\includegraphics[keepaspectratio, scale=0.55]{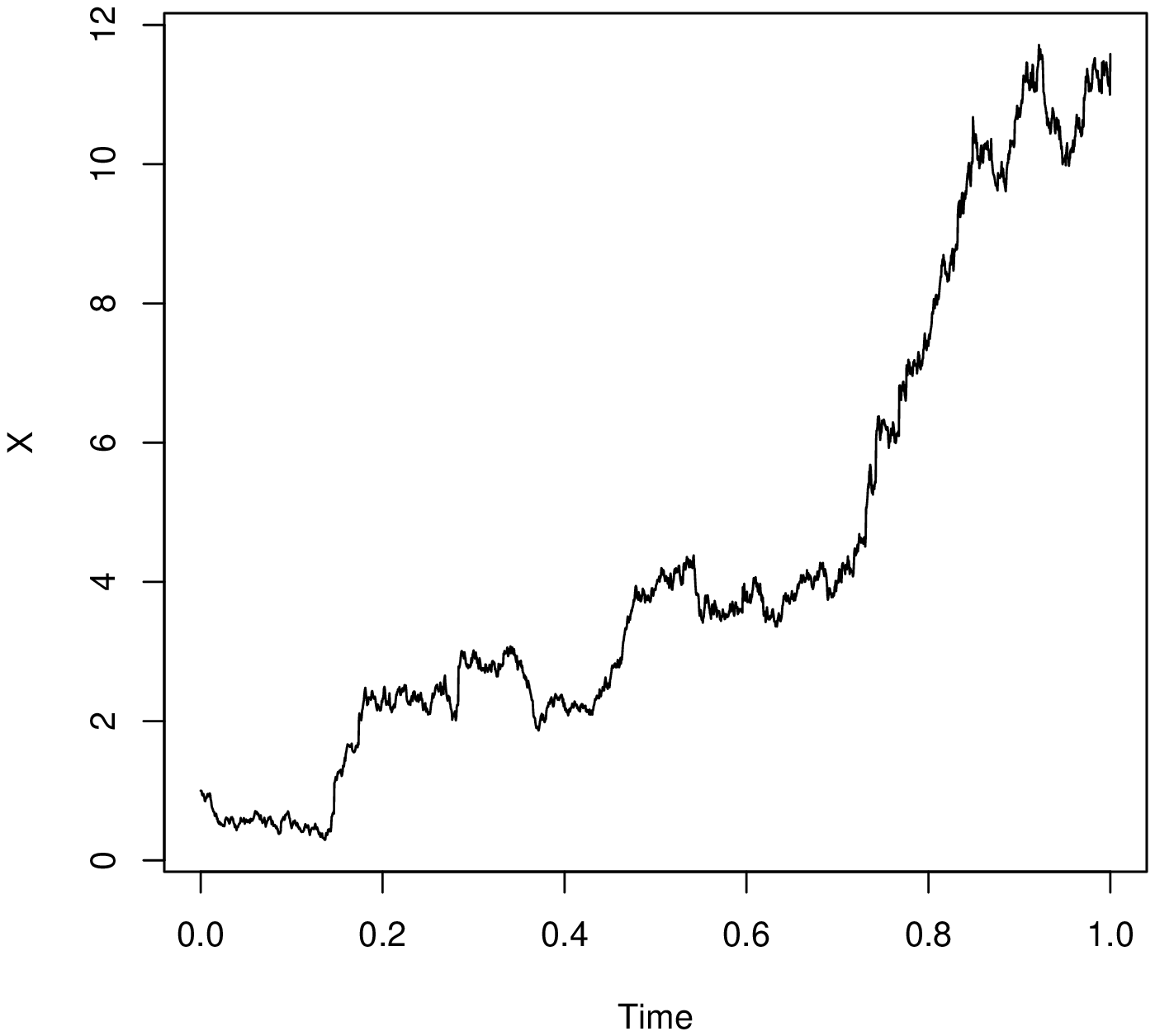}
					{(a) Sample path of $X$}
				\end{center}
			\end{minipage}
			
			\begin{minipage}{0.5\hsize}
				\begin{center}
					\includegraphics[keepaspectratio, scale=0.55]{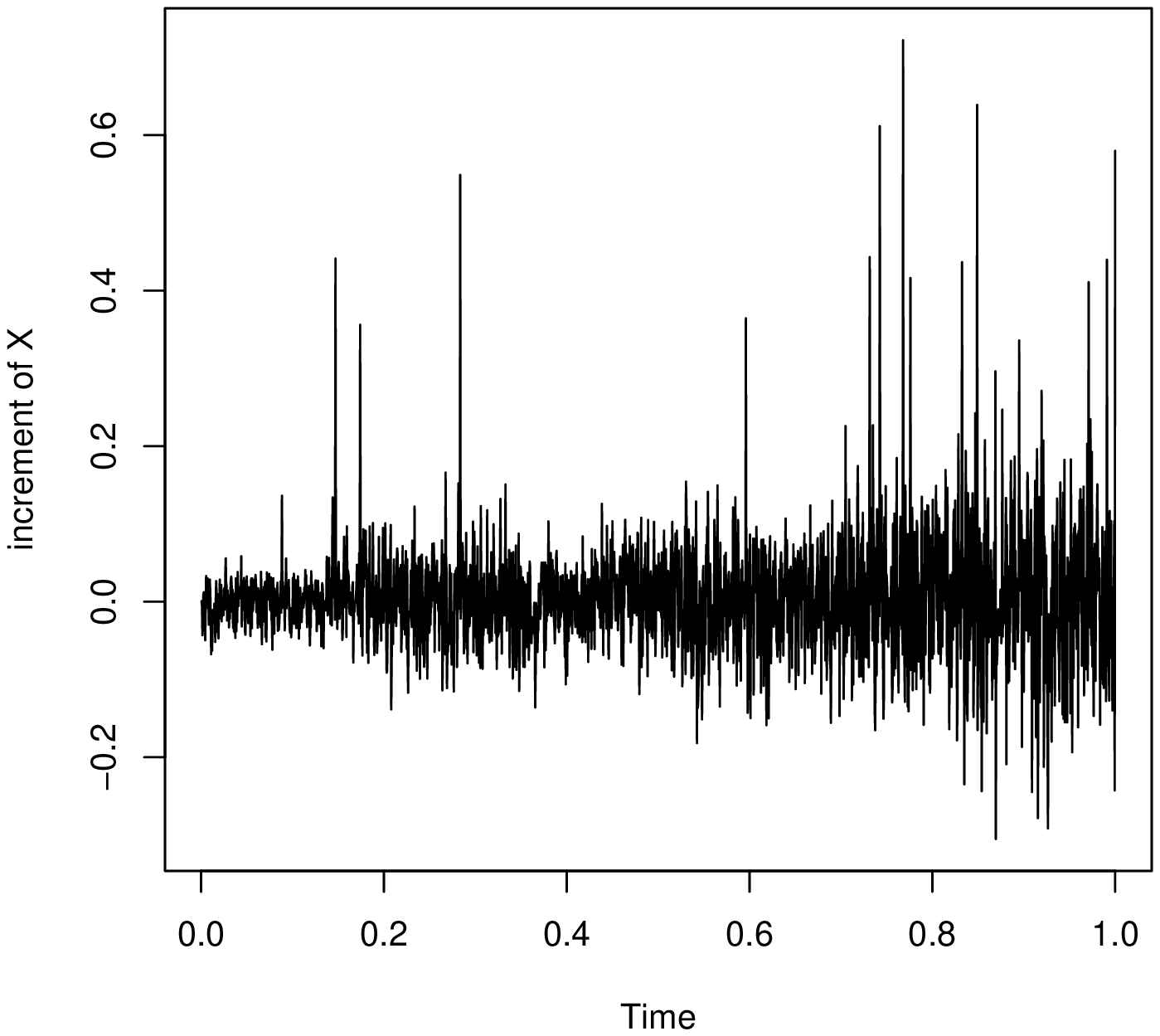}
					{(b) Increment of $X$}
				\end{center}
			\end{minipage}
			
		\end{tabular}
		\caption{Sample path of $X$ and its increments ($\lambda = 30$)}
		
		\label{fig1}		
	\end{center}
	
\end{figure}

{\colorg Table \ref{Table05} shows the summary of error ratios (percentage deviation of estimated values from the true value
for each estimator)}, and 
Figure \ref{fig2} 
{\tred gives their box plots.} 
In this case, both BV ({\tt bv}) and minRV ({\tt mrv}) seem to suffer from upward bias due to jumps. In particular, the BV deviates from the true value considerably. 
On the other hand, GRV with normalization perform well with errors concentrating around zero ({\tt grv.lgrv, grv.mrv}). Note that, although WGRV performs relatively well, it seems to have a small upward bias ({\tt wgrv.lgrv, wgrv.mrv}). 
This suggests that, if there are many large jumps, 
using an upper quantile ($V_{(s_n(\alpha))}$) may sometimes lead to biases 
rather than obtaining a robust estimate. 

The three right {\tred box plots} 
in this figure ({\tt grv[0.20], grv[0.10], grv[0.05]}) are the results of GRV without 
normalizing increments by local-global filters, with the cut-off ratio 
$\alpha = 0.2, \ 0.1, \ 0.05$, respectively.
We see that they seem to be less precise (especially when $\alpha$ is large {\tred extremely small}) than 
GRV or WGRV with local volatility.  
This result suggests that, if we do not normalize increments by spot volatilities in the case of non-constant volatility, we end up obtaining 
inappropriate estimates. 

Intuitively, when we ignore normalization, we tend to eliminate increments 
where volatility is high (because they are typically large), even if they come from the Brownian motion, while keeping relatively small jumps which we should actually remove. 
In addition, theoretically, the adjusting constant $q(\alpha)$ in the definition of GRV 
(\ref{0211260834}) comes from the standard normal distribution. 
Therefore, when the volatility is non-constant, we should normalize the increments 
$|\Delta_i {\tred X}|$ by local volatility to make them 
approximately {\tred standard} normally distributed.

\renewcommand{\arraystretch}{1.5}
\begin{table}[H] \centering 
	\caption{Summary of error ratios: $\lambda = 30$}
	\vspace{2mm}
	\begin{tabular}{@{\extracolsep{5pt}} lccccc} 
		
		\hline
		& Min. & 1st Qu. & Median & 3rd Qu. & Max. \\ 
		\hline \hline
		{\tt trv[0.45]} & -29.88 & -28.16 & -27.42 & -26.45 & -21.64 \\ 
		{\tt trv[0.20]} & -9.40 & -3.57 & -1.95 & -0.67 & 3.54 \\ 
		{\tt trv[0.10]} & 0.32 & 3.38 & 4.55 & 6.11 & 21.61 \\ 
		{\tt bv} & -0.25 & 2.66 & 3.87 & 5.11 & 8.98 \\ 
		{\tt mrv} & -2.61 & 0.36 & 1.43 & 2.48 & 5.75 \\ 
		{\tt grv.lgrv[0.20]} & -6.65 & -1.03 & -0.10 & 0.88 & 3.64 \\ 
		{\tt grv.mrv[0.20]} & -6.65 & -1.33 & -0.40 & 0.64 & 3.65 \\ 
		{\tt wgrv.lgrv[0.20]} & -3.50 & -0.39 & 0.60 & 1.46 & 3.68 \\ 
		{\tt wgrv.mrv[0.20]} & -3.53 & -0.50 & 0.44 & 1.32 & 3.72 \\ 
		{\tt grv[0.20]} & -9.29 & -4.39 & -3.21 & -1.91 & 2.15 \\ 
		{\tt grv[0.10]} & -10.26 & -2.84 & -1.76 & -0.76 & 2.71 \\ 
		{\tt grv[0.05]} & -13.85 & -4.41 & -1.77 & -0.42 & 3.38 \\ 
		{\tt grv.lgrv.mov} & -9.13 & -1.18 & -0.00 & 0.95 & 3.55 \\ 
		{\tt wgrv.lgrv.mov} & -3.19 & -0.14 & 0.80 & 1.62 & 4.61 \\ 
		\hline \\[-1.8ex] 
	\end{tabular} 
	\label{Table05}
\end{table} 
\renewcommand{\arraystretch}{1.0}

\begin{figure}[H]
\begin{center}
		\begin{tabular}{c}
	\begin{minipage}{0.4\hsize}
		\includegraphics[keepaspectratio, scale=0.7]{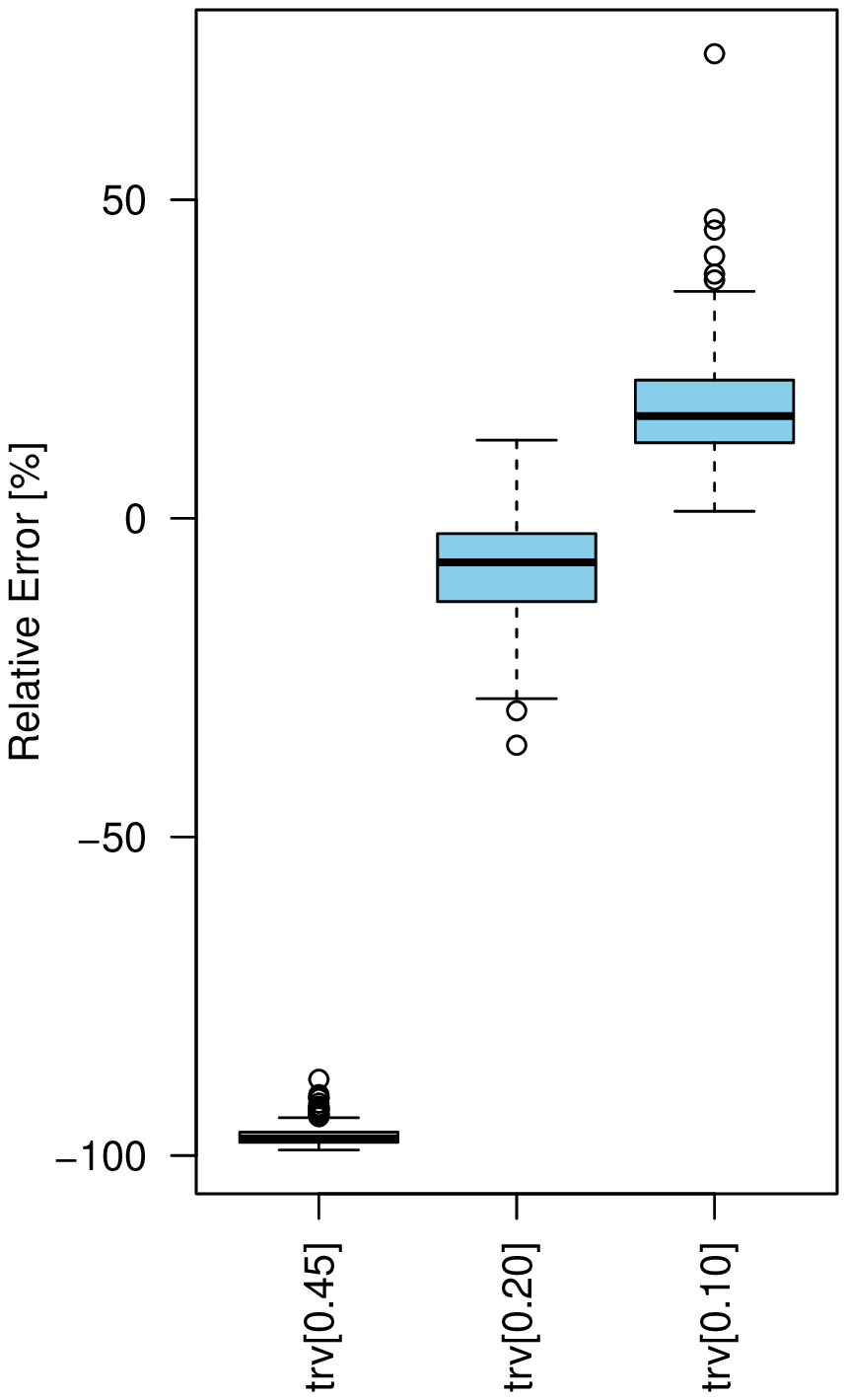}
	\end{minipage}
	
	\begin{minipage}{0.6\hsize}
		\includegraphics[keepaspectratio, scale=0.7]{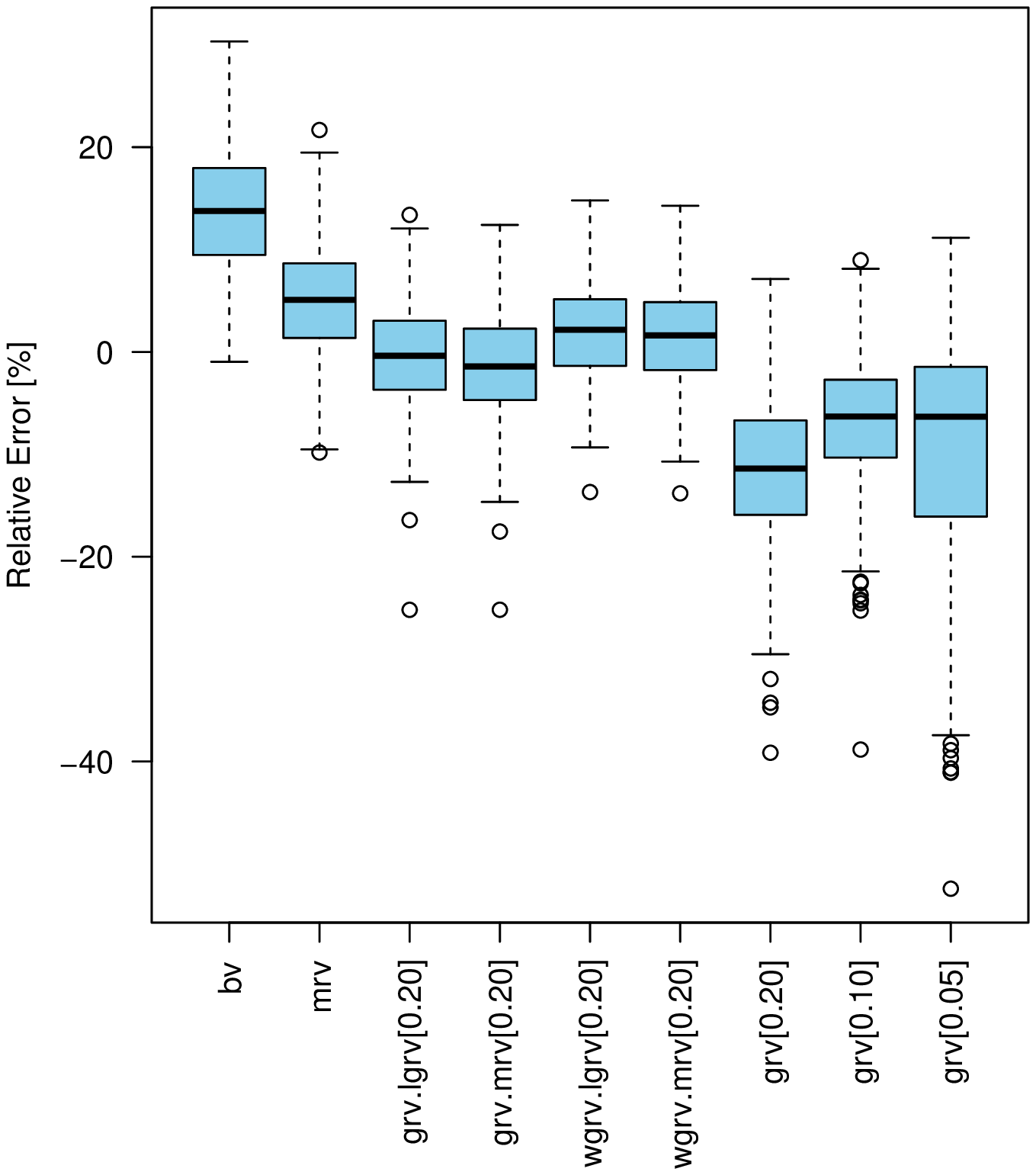}
	\end{minipage}
\end{tabular}
	\caption{Error ratios [\%] for the case of high intensity: $\lambda = 30$}
	\label{fig2}
\end{center}
\end{figure}

{\colorg 
The good news is that they also perform well even in the case of 
extremely high intensity. 
We consider $\lambda = 50$ here to see their accuracy. 
Figure \ref{fig2EX} shows a sample path and its increments. 
It is obvious there are numerous jumps and one can easily imagine 
that the standard realized volatility estimator can never 
estimate the true volatility. 
{\colorg 
Table \ref{Table08} and Figure \ref{fig9} show error ratios of each estimator for $\lambda=50$. 
It shows that GRV and WGRV with cut-off ratio $\alpha=0.2$ perform well even 
in the case of high intensity of jumps. 
}

\begin{figure}[H]
	\begin{center}
		
		\begin{tabular}{c}

			\begin{minipage}{0.5\hsize}
				\begin{center}
					\includegraphics[keepaspectratio, scale=0.55]{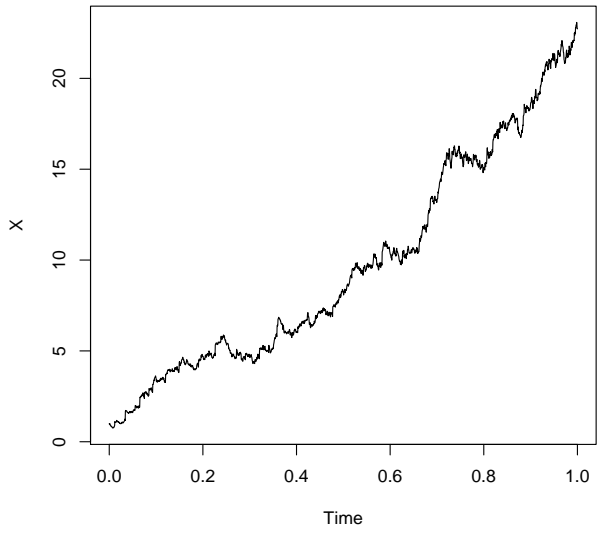}
					{(a) Sample path of $X$}
				\end{center}
			\end{minipage}
			
			\begin{minipage}{0.5\hsize}
				\begin{center}
					\includegraphics[keepaspectratio, scale=0.55]{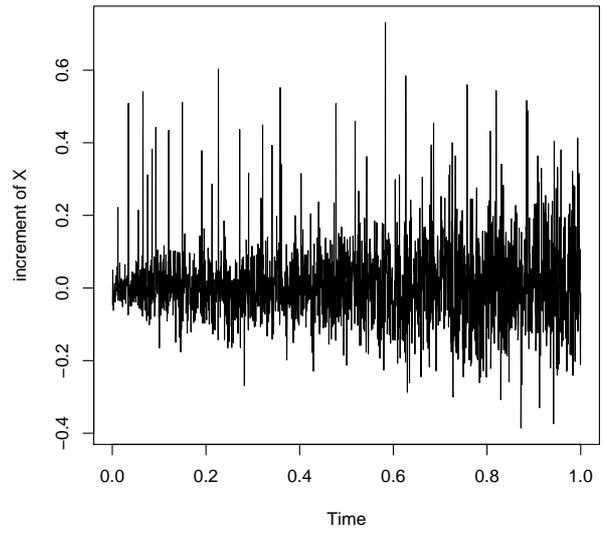}
					{(b) Increment of $X$}
				\end{center}
			\end{minipage}
			
		\end{tabular}
		\caption{Sample path of $X$ and its increments ($\lambda = 50$)}
		\label{fig2EX}		
	\end{center}
	
\end{figure}

}

\renewcommand{\arraystretch}{1.5}
\begin{table}[H] \centering
	\caption{Summary of error ratios: $\lambda = 50,  \ \alpha = 0.2$} 
	\vspace{2mm}
	\begin{tabular}{lccccc}
		\hline 
		& Min. & 1st Qu. & Median & 3rd Qu. & Max. \\ 
		\hline \hline 
		{\tt trv[0.45]} & -99.52 & -98.90 & -98.52 & -98.09 & -94.42 \\ 
		{\tt trv[0.20]} & -52.18 & -27.61 & -19.86 & -12.98 & 2.35 \\ 
		{\tt trv[0.10]} & 0.78 & 11.76 & 15.51 & 20.33 & 37.69 \\ 
		{\tt bv} & 2.71 & 13.60 & 17.47 & 20.89 & 37.19 \\ 
		{\tt mrv} & -6.36 & 3.40 & 7.72 & 11.58 & 27.20 \\ 
		{\tt grv.lgrv[0.20]} & -38.58 & -7.51 & -2.38 & 1.82 & 18.41 \\ 
		{\tt grv.mrv[0.20]} & -39.65 & -8.54 & -3.51 & 0.81 & 16.58 \\ 
		{\tt wgrv.lgrv[0.20]} & -10.64 & 0.87 & 3.99 & 7.44 & 25.22 \\ 
		{\tt wgrv.mrv[0.20]} & -12.60 & -0.08 & 3.18 & 6.74 & 23.23 \\ 
		{\tt grv[0.20]} & -33.92 & -15.51 & -11.40 & -7.01 & 5.42 \\ 
		{\tt grv[0.10]} & -56.46 & -23.57 & -11.63 & -4.84 & 9.21 \\ 
		{\tt grv[0.05]} & -66.14 & -40.55 & -29.94 & -17.97 & 5.88 \\ 
		{\tt grv.lgrv.mov} & -48.17 & -16.03 & -6.55 & -0.14 & 16.16 \\ 
		{\tt wgrv.lgrv.mov} & -10.31 & 2.16 & 5.85 & 9.61 & 31.21 \\ 
		\hline
	\end{tabular}
	\label{Table08}
\end{table} 
\renewcommand{\arraystretch}{1.0}

\begin{figure}[H]
	\begin{center}
		\begin{tabular}{c}
			\begin{minipage}{0.4\hsize}
				\includegraphics[keepaspectratio, scale=0.7]{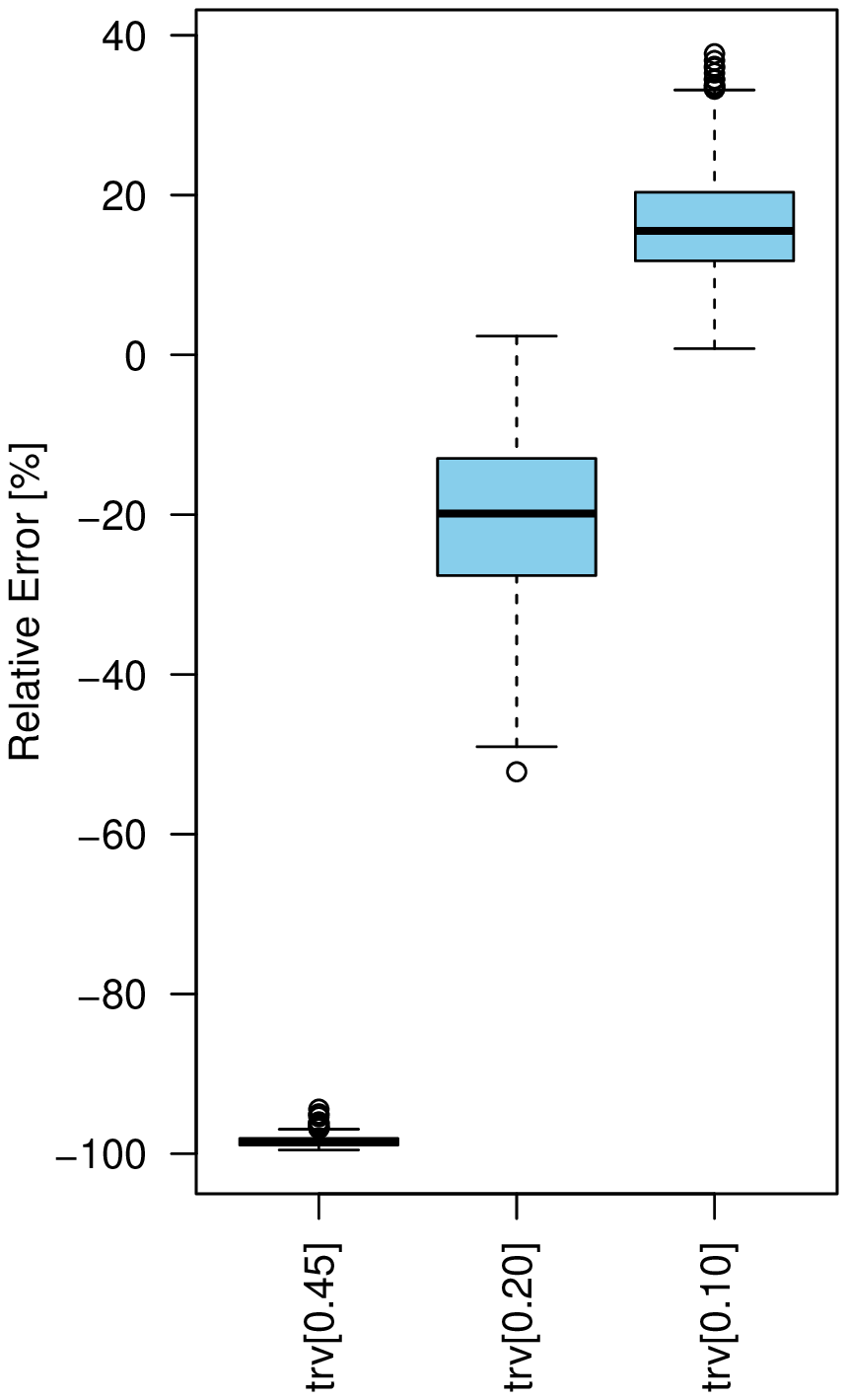}
			\end{minipage}
			
			\begin{minipage}{0.6\hsize}
				\includegraphics[keepaspectratio, scale=0.7]{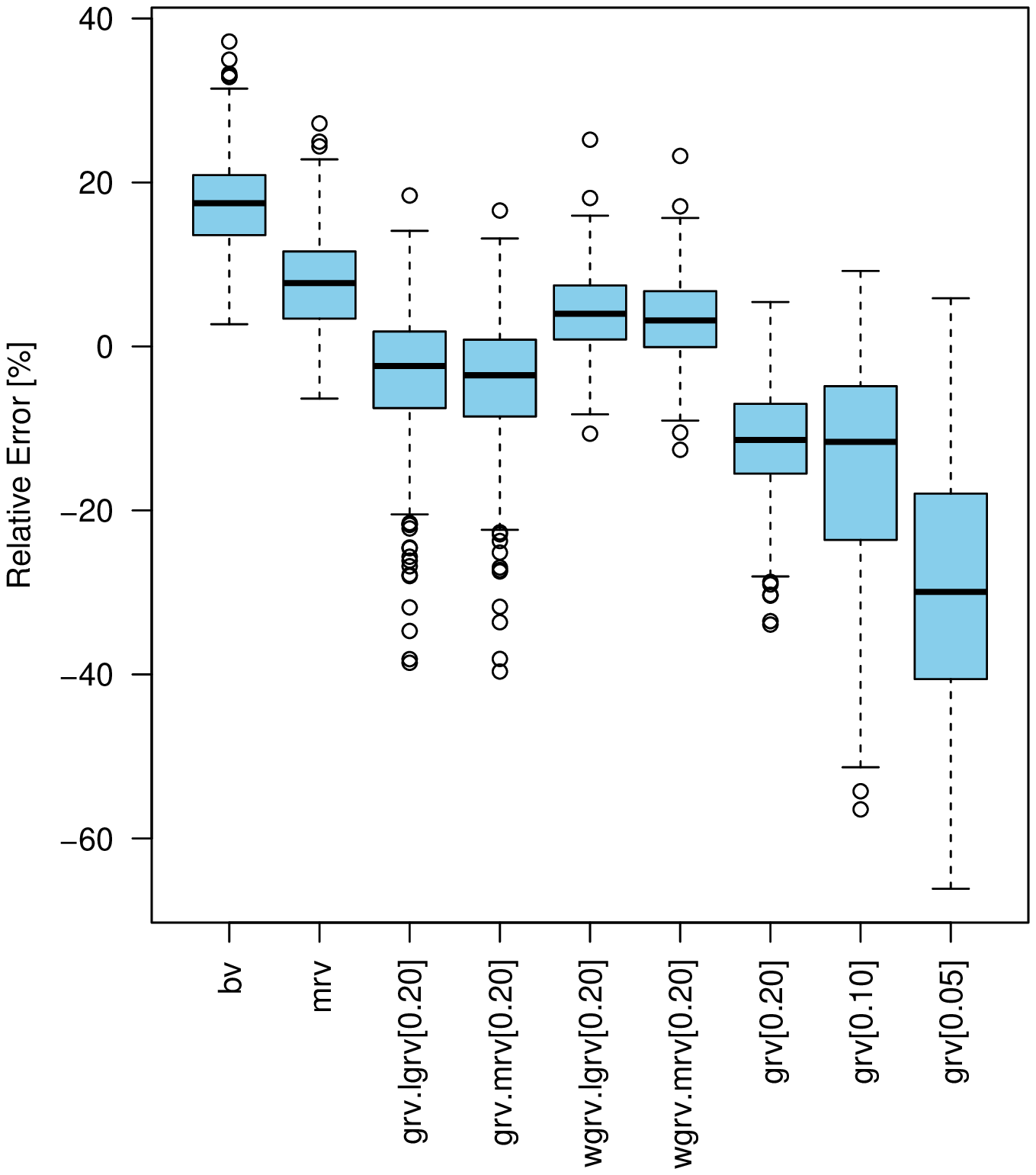}
			\end{minipage}
		\end{tabular}
		\caption{Error ratios [\%] for the case of high intensity: $\lambda = 50$}
		\label{fig9}
	\end{center}
\end{figure}

{\colorg
By taking $\alpha$ larger, accuracy improves. 
Table \ref{Table04} shows the error ratios of each estimator in the case of $\lambda = 50$, 
with $\alpha = \alpha_0$ ranging from 0.1 to 0.5. 
We can see that, for large cut-off ratio ($\alpha=\alpha_0 = 0.4, 0.5$), 
GRV and WGRV with spot volatilities 
({\tt grv.lgrv}, {\tt grv.mrv}, {\tt wgrv.lgrv}, {\tt wgrv.mrv}) still 
perform well. 
Looking in more detail, we see that the local GRV outperforms 
local minRV for both GRV and WGRV. 
This example imply that we should take a cut-off ratio 
quite large in order to obtain a precise estimate.  
	
}

\vspace{5mm}
\renewcommand{\arraystretch}{1.5}
\begin{table}[H] 
\begin{center}
	\caption{Error ratios [\%] for the case of extremely high intensity: $\lambda = 50$}
	\vspace{2mm}
\begin{tabular}{lccccc}
	\hline
	& \multicolumn{5}{c}{Cut-off ratio ($\alpha = \alpha_0$)}  \\ \cline{2-6}
	& 0.1    & 0.2    & 0.3    & 0.4    & 0.5    \\ \hline
	{\tt trv[0.45]} & -98.40 & -98.40 & -98.40 & -98.40 & -98.40 \\
	{\tt trv[0.20]} & -20.35 & -20.35 & -20.35 & -20.35 & -20.35 \\
	{\tt trv[0.10]} & 16.19  & 16.19  & 16.19  & 16.19  & 16.19  \\
	{\tt bv}            & 17.60  & 17.60  & 17.60  & 17.60  & 17.60  \\
	{\tt mrv}           & 7.71   & 7.71   & 7.71   & 7.71   & 7.71   \\
	{\tt grv.lgrv}[$\alpha$]      & -18.04 & -3.49  & -0.41  & -0.96  & -1.47  \\
	{\tt grv.mrv}[$\alpha$]       & -18.05 & -4.35  & -2.62  & -3.90  & -4.89  \\
	{\tt wgrv.lgrv}[$\alpha$]     & 10.28  & 4.05   & 2.27   & 1.35   & 0.59   \\
	{\tt  wgrv.mrv} [$\alpha$]     & 10.87  & 3.21   & 1.06   & -0.30  & -1.40  \\
	{\tt  grv[0.20]} & -11.50 & -11.50 & -11.50 & -11.50 & -11.50 \\
	{\tt  grv[0.10]} & -14.60 & -14.60 & -14.60 & -14.60 & -14.60 \\
	{\tt  grv[0.05]} & -28.85 & -28.85 & -28.85 & -28.85 & -28.85 \\
	{\tt  grv.lgrv.mov}  & -8.31  & -8.76  & -8.78  & -8.80  & -8.83  \\
	{\tt wgrv.lgrv.mov} & 4.99   & 6.26   & 6.48   & 6.63   & 6.73  \\ 
	\hline
\end{tabular}
\label{Table04}
\end{center}
\end{table}
\renewcommand{\arraystretch}{1.0}

\subsubsection{The case of moderate intensity: GRV with a shrinking cut-off ratio}
Next, we consider the case of low intensity. 
In this case, we can use shrinking cut-off rate. 
Recall that the shrinking cut-off rate is defined by  
$\alpha_n = \lfloor B n^{\delta_1} \rfloor / n$. 
In this simulation, we set $B = 10$ and $\delta_1 = 0.45$, so the 
cut-off rate is then $\alpha_n = 0.1525$. 

The {\colorg error ratios} are shown in Table \ref{Table07} Figure \ref{fig3}. 
All global-filtering estimators perform well (for GRVs with fixed cut-off ratio, 
we set $\alpha = 0.2$ as before). 
These results suggest that if there are not so many jumps in the data, 
it would be advisable to use as many data as possible by making 
the cut-off ratio small. 
{\colorb
Note that TRV still has bias, especially for $\rho=0.45$. 
This implies that the accuracy of estimation is still highly 
vulnerable to the choice of $\rho$ for TRV, 
even in the case of moderate intensity of jumps. 	 
}

\renewcommand{\arraystretch}{1.5}
\begin{table}[H] \centering 
	\caption{Summary of error ratios: $\lambda = 5$} 
	\vspace{2mm}
	\begin{tabular}{@{\extracolsep{5pt}} lccccc} 
		\hline 
		& Min. & 1st Qu. & Median & 3rd Qu. & Max. \\ 
		\hline \hline
		{\tt trv[0.45]} & -29.61 & -27.40 & -26.09 & -23.89 & -18.27 \\ 
		{\tt trv[0.20]} & -2.79 & -0.58 & 0.14 & 0.85 & 3.86 \\ 
		{\tt trv[0.10]} & -2.39 & 1.16 & 2.52 & 4.25 & 13.51 \\ 
		{\tt bv} & -2.26 & 0.31 & 1.29 & 2.24 & 6.34 \\ 
		{\tt mrv} & -3.63 & -0.69 & 0.29 & 1.24 & 4.27 \\ 
		{\tt grv.lgrv[0.20]} & -3.62 & -1.44 & -0.55 & 0.21 & 4.79 \\ 
		{\tt grv.mrv[0.20]} & -3.59 & -1.40 & -0.55 & 0.30 & 5.00 \\ 
		{\tt wgrv.lgrv[0.20]} & -3.52 & -1.11 & -0.39 & 0.44 & 4.55 \\ 
		{\tt wgrv.mrv[0.20]} & -3.53 & -1.06 & -0.40 & 0.46 & 4.62 \\ 
		{\tt grv[0.20]} & -7.14 & -2.82 & -1.75 & -0.75 & 2.90 \\ 
		{\tt grv[0.10]} & -5.46 & -2.15 & -1.32 & -0.46 & 3.17 \\ 
		{\tt grv[0.05]} & -4.15 & -1.63 & -0.86 & -0.10 & 3.23 \\ 
		{\tt grv.lgrv.mov} & -3.54 & -1.32 & -0.50 & 0.27 & 4.87 \\ 
		{\tt wgrv.lgrv.mov} & -3.15 & -0.96 & -0.34 & 0.47 & 4.17 \\ 
		\hline 
	\end{tabular} 
	\label{Table07} 
\end{table} 
\renewcommand{\arraystretch}{1.0}

\begin{figure}[H]
		\begin{tabular}{c}
			\begin{minipage}{0.4\hsize}
					\includegraphics[keepaspectratio, scale=0.7]{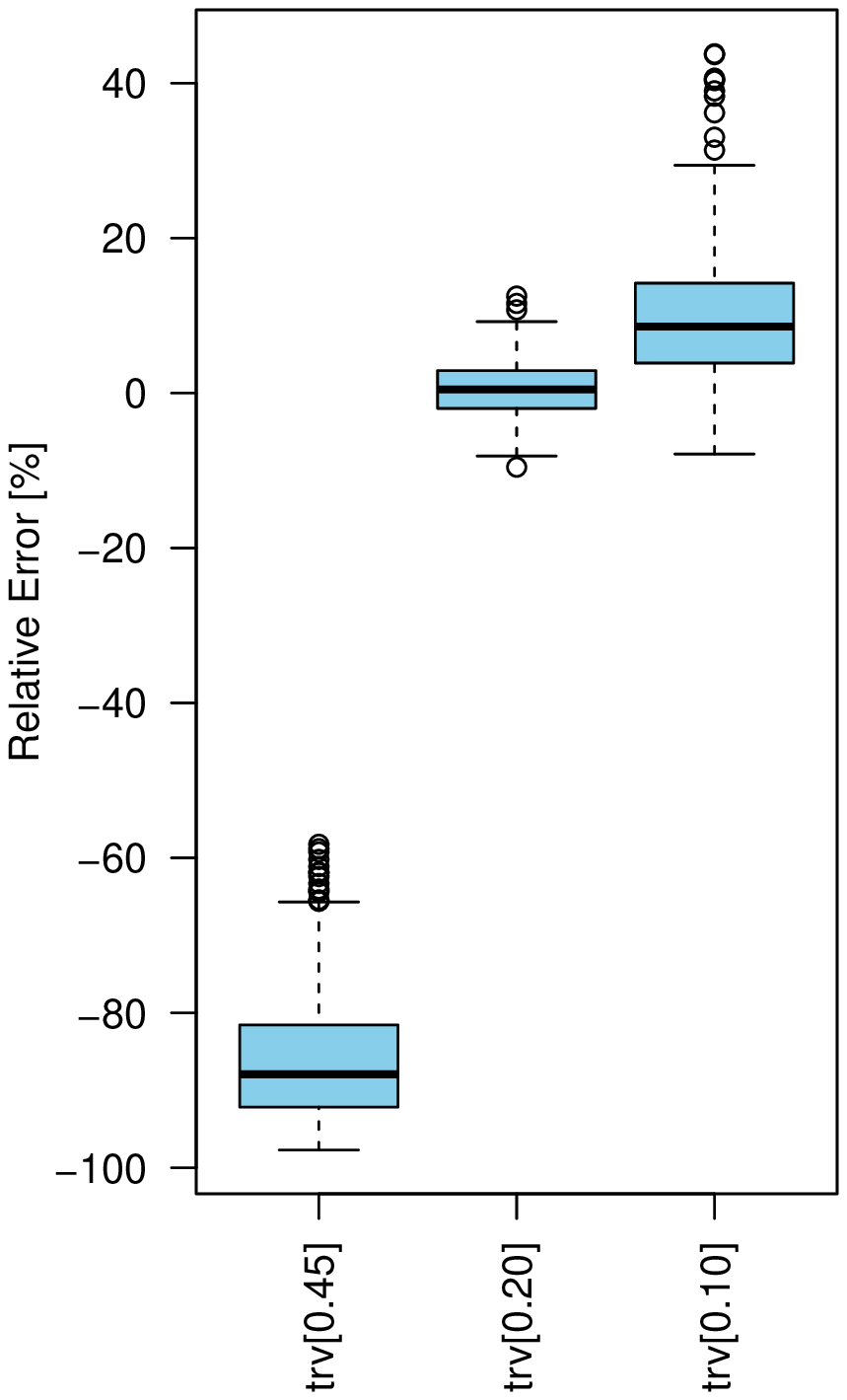}
			\end{minipage}
			
			\begin{minipage}{0.6\hsize}
					\includegraphics[keepaspectratio, scale=0.7]{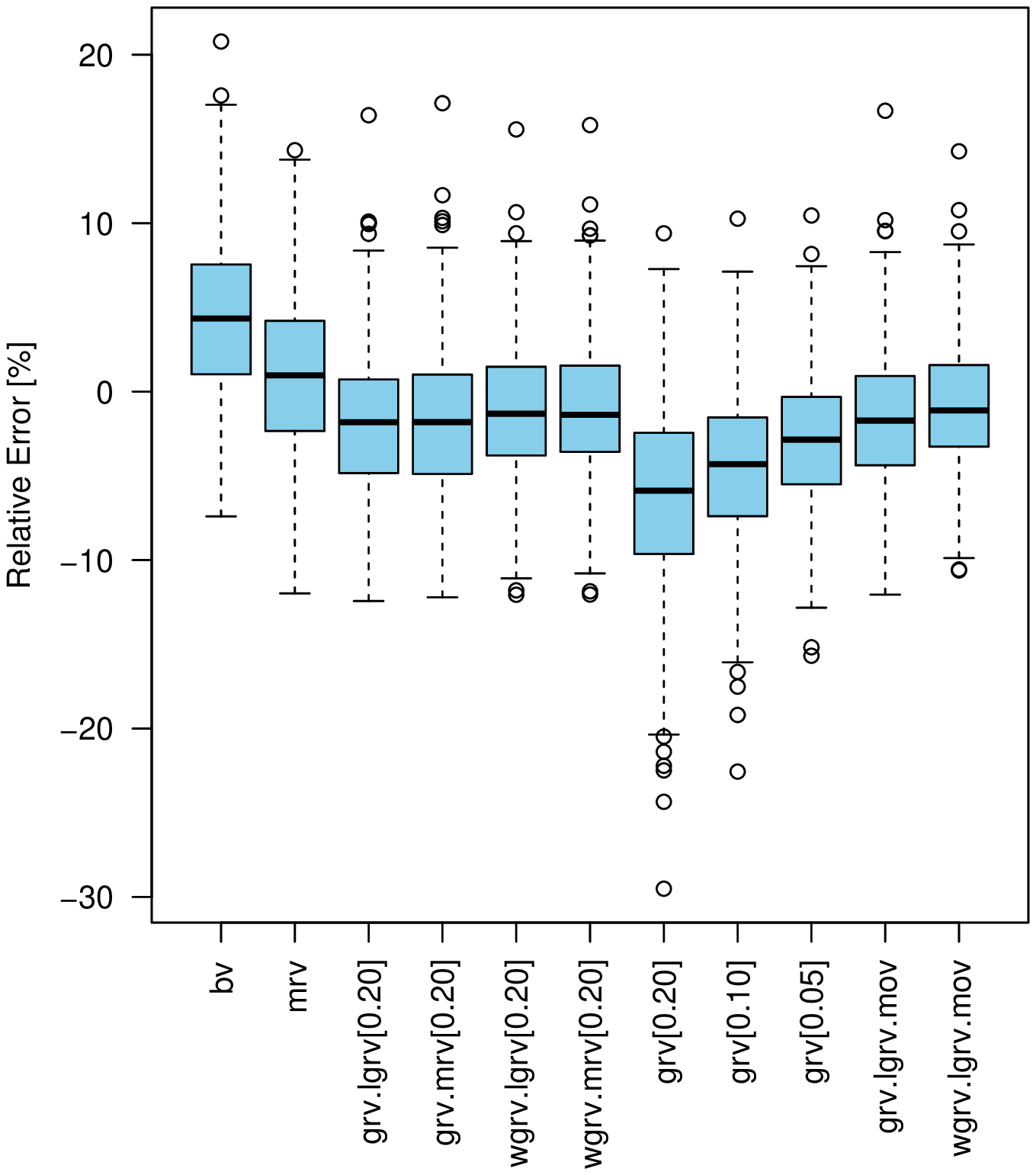}
			\end{minipage}
		\end{tabular}
		\caption{Error ratios [\%] for the case of low intensity: $\lambda = 5$}
		\label{fig3}		
\end{figure}

{\colorg
For GRV and WGRV with shrinking cut-off ratio, 
we proved the asymptotic mixed normality. 
Hence, the distribution of the Studentized errors 
$\Gamma^{-1/2} \sqrt{n} ({\bf V}_n - \Theta)$ and 
$\Gamma^{-1/2} \sqrt{n} ({\bf W}_n - \Theta)$ are expected to follow the standard normally distribution. 

Figure \ref{fig7} shows QQ plots comparing theoretical quantiles of the standard normal distribution and the Studentized errors of GRV and WGRV estimators with shrinking threshold, BV and minRV. 
In this example, {\tt wgrv.lgrv.mov} outperforms the others. 
It is close to the standard normal distribution. 
On the other hand, {\tt grv.lgrv.mov} seems to deviate {\tred from} $N(0,1)$.  
We can also see that {\tt bv} are far from $N(0,1)$, implying 
that it is not appropriate even in the case of low intensity.

\begin{figure}[H]
	\begin{center}
		\begin{tabular}{c}
			\begin{minipage}{0.5\hsize}
				\includegraphics[keepaspectratio, scale=0.7]{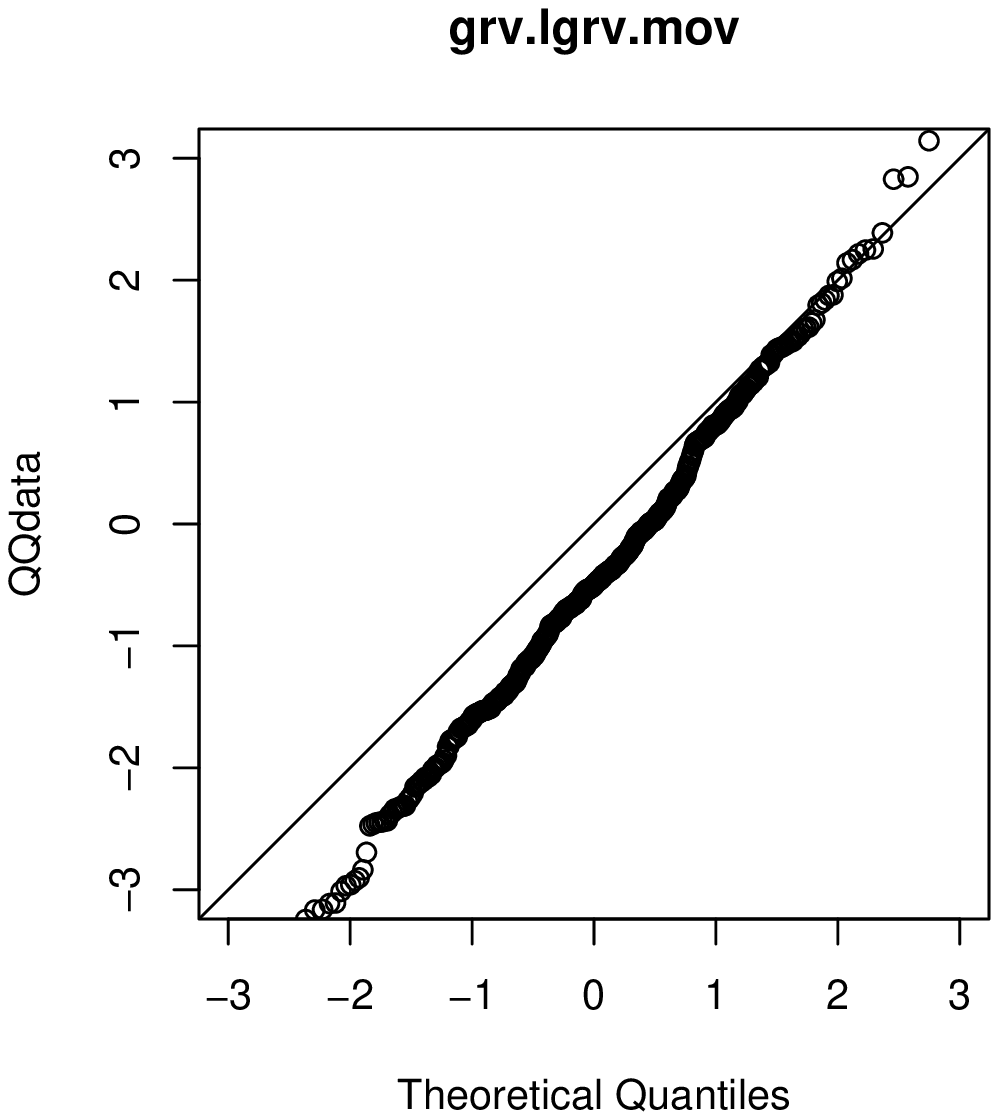}
			\end{minipage}
			
			\begin{minipage}{0.5\hsize}
				\includegraphics[keepaspectratio, scale=0.7]{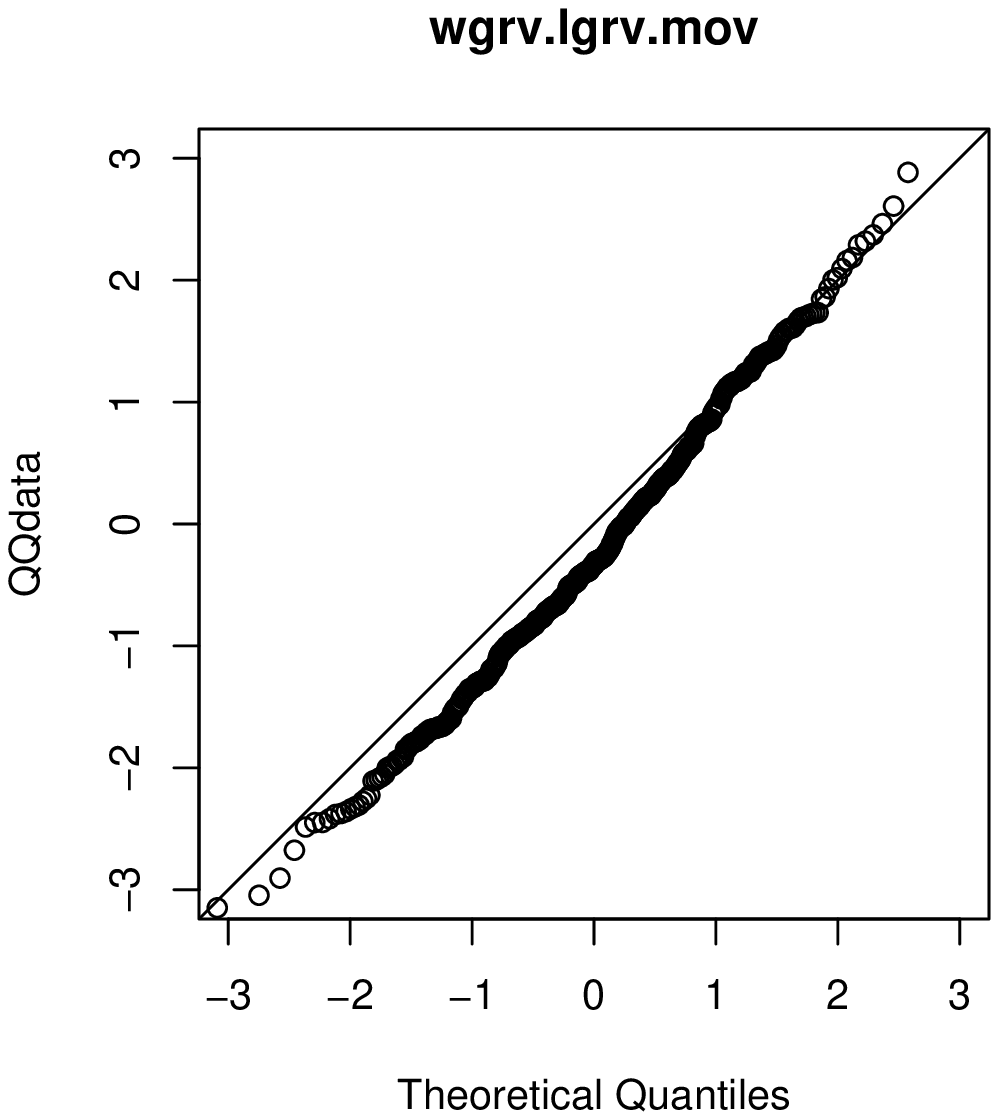}
			\end{minipage} \\

			\begin{minipage}{0.5\hsize}
				\includegraphics[keepaspectratio, scale=0.7]{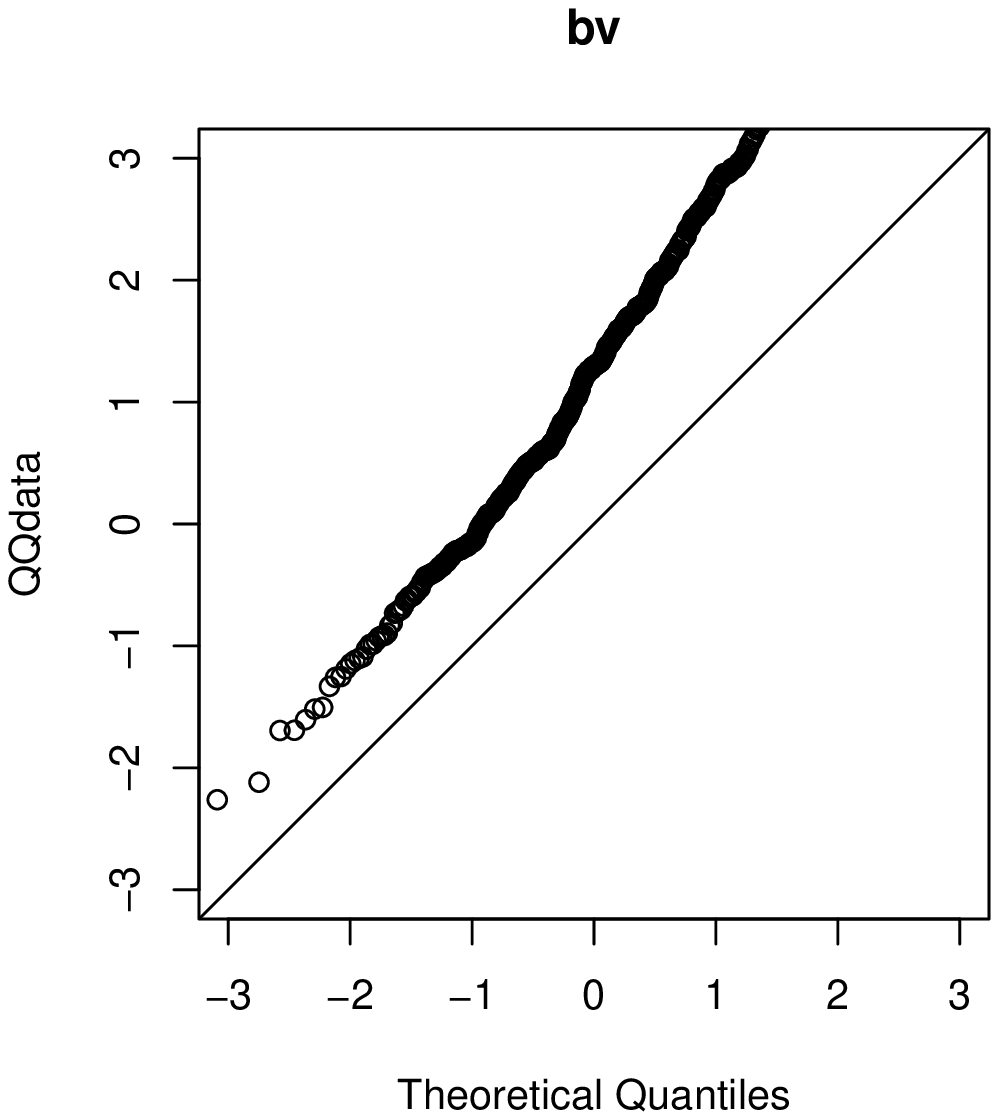}
			\end{minipage}
			
			\begin{minipage}{0.5\hsize}
				\includegraphics[keepaspectratio, scale=0.7]{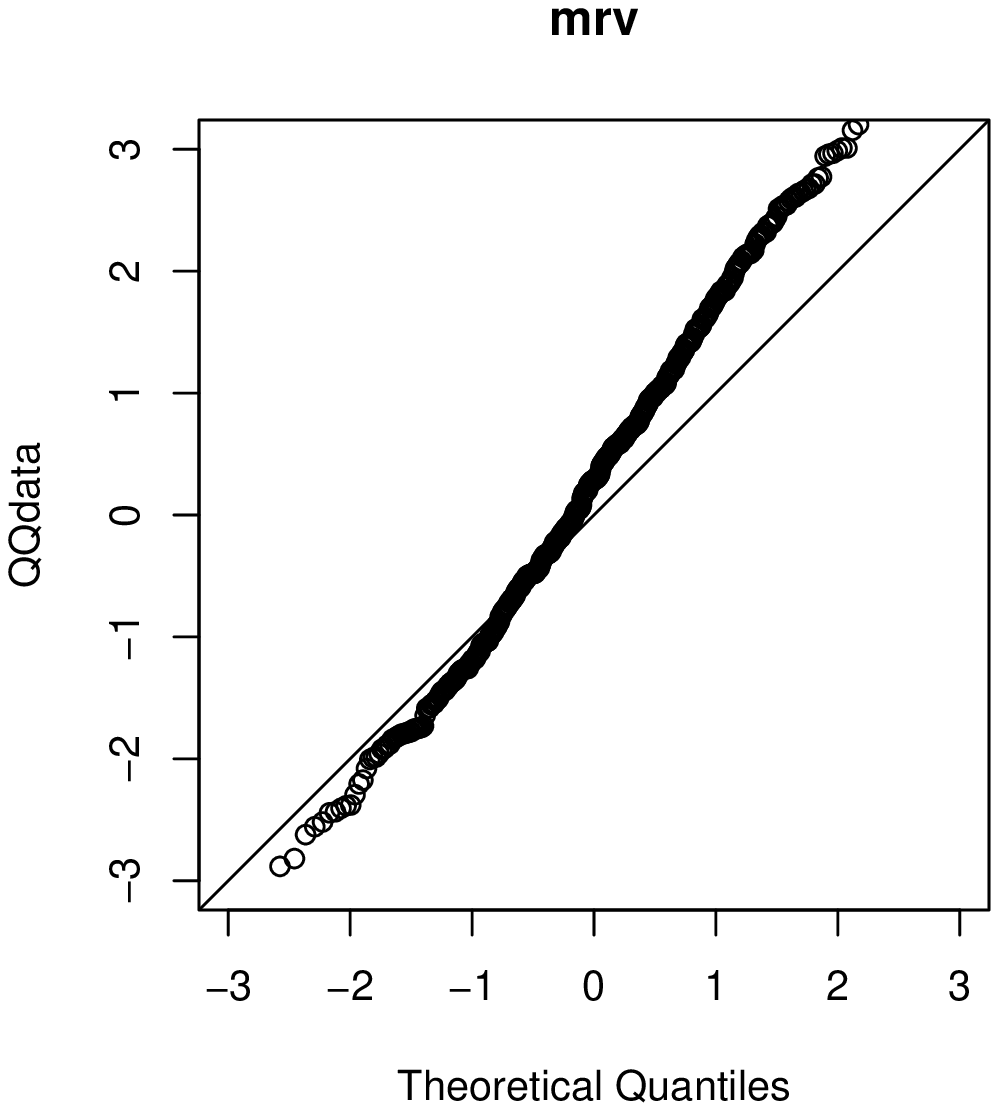}
			\end{minipage}		

		\end{tabular}
		\caption{QQ plot for Studentized errors: $\lambda = 5$}
		\label{fig7}
	\end{center}
\end{figure}

}

{\colorb 
{\colorg The important tuning parameter for the shrinking threshold GRV} is the exponent $\delta_1$, 
an appropriate choice of which may strongly depend on the intensity of jumps. Recall that small $\delta_1$ means that we keep almost all the samples untrimmed. 
Table \ref{Table02} shows average {\colorg error ratios}  
of GRV and WGRV with shrinking 
cut-off ratio for several values of intensity $\lambda$ and the parameter $\delta_1$. 
For moderate intensity ($\lambda = 5, 10$), 
{\colorg the average ratios are not so large for small $\delta_1$. }
On the other hand, for high intensity ($\lambda = 30, 50$), this is not the case. 
Indeed, as for GRV, estimation errors are quite large downward for small $\delta_1$. 
This can be interpreted that its multiplication of $q(\alpha_n)^{-1}$ for GRV is  insufficient to compensate its 
elimination of jumps (small $\delta_1$ implies small $\alpha_n$, making $q(\alpha_n)^{-1}$ close to 1). 
Moreover, as for WGRV, there occur large upward biases for small $\delta_1$, 
since it keeps almost large increments and uses an extremely large increment 
for winsorization. 

It is worth noting that large $\delta_1$ makes both GRV and WGRV accurate 
{\colorg to a certain extent}, 
even in the case of high intensity of jumps. 
Thus, in practice, one may use shrinking cut-off GRV and WGRV 
by setting the tuning parameter $\delta_1$ sufficiently close to $1/2$. 

{\colorg 
However, as Figure \ref{fig8} implies, the errors are not normally distributed as 
theory predicts when the intensity of jumps is extremely high. 
We should be aware that GRV and WGRV with shrinking cut-off may suffer from 
some biases in the case of extremely intensive jumps. 
We may consider using a large fixed cut-off ratio (as discussed in 
the previous subsection) in such a situation. 
}

\begin{figure}[H]
	\begin{center}
		\begin{tabular}{c}
			\begin{minipage}{0.5\hsize}
				\includegraphics[keepaspectratio, scale=0.7]{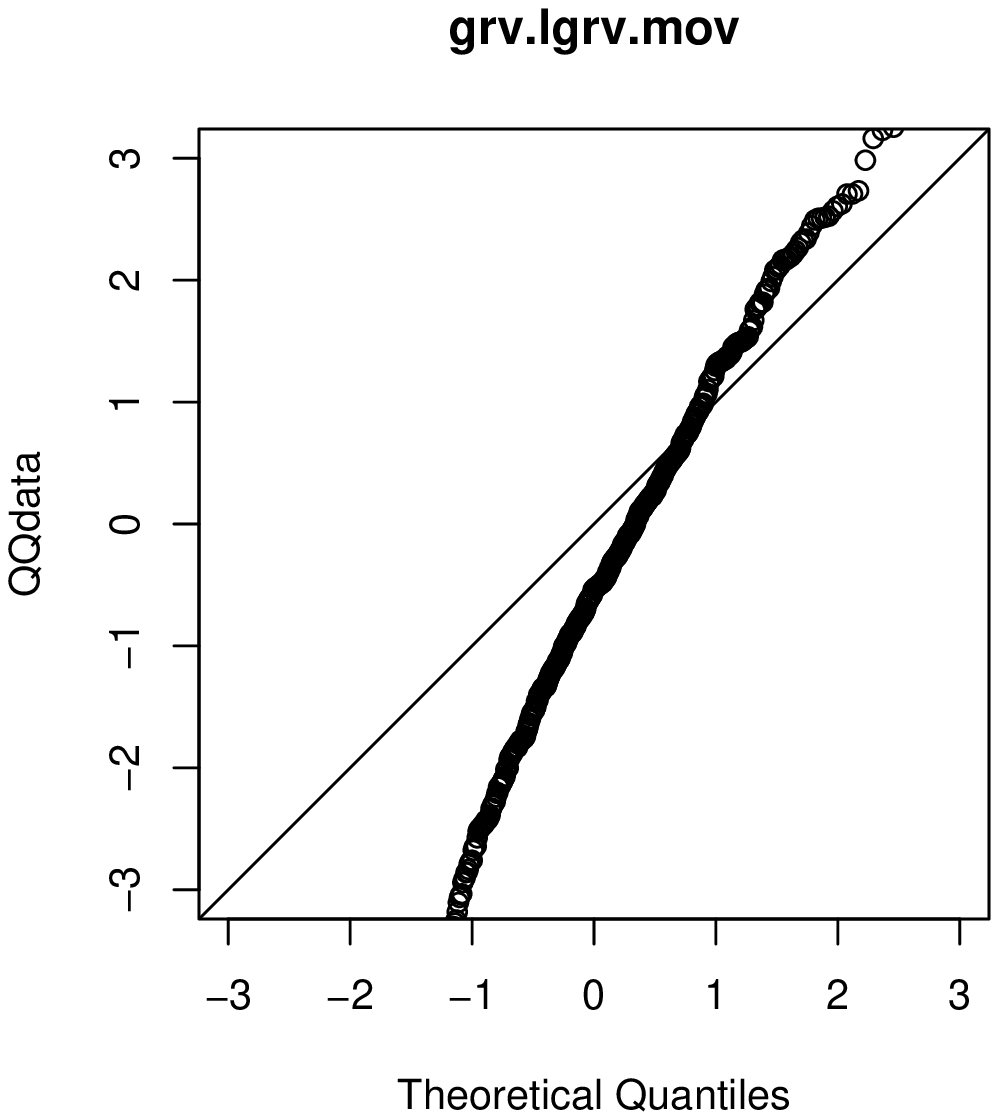}
			\end{minipage}
			
			\begin{minipage}{0.5\hsize}
				\includegraphics[keepaspectratio, scale=0.7]{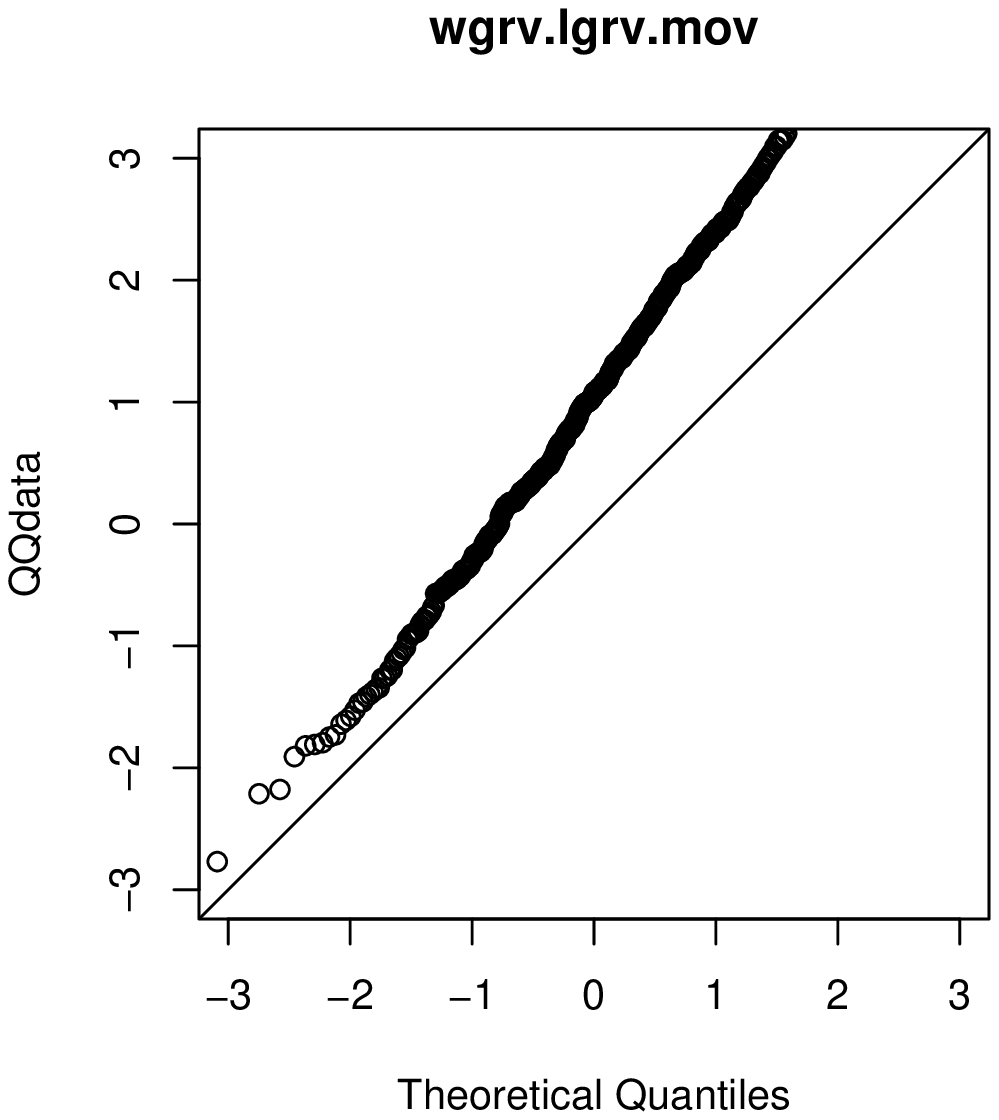}
			\end{minipage} 
		\end{tabular}
		\caption{QQ plot for Studentized errors: $\lambda = 50, \ \delta_1 = 0.49$}
		\label{fig8}
	\end{center}
\end{figure}

}

\renewcommand{\arraystretch}{1.5}
\begin{table}[H]
	\caption{Average error ratios [$\%$] of GRV and WGRV with shrinking cut-off ratio}
	
	\begin{minipage}{0.5\hsize}
		\begin{center}
		\vspace{3mm}
		{(a) GRV}\\
		\vspace{3mm}
		\begin{tabular}{ccccc}\hline
			\multicolumn{1}{c}{} & \multicolumn{4}{c}{Intensity of jumps ($\lambda$) }\\
			\cline{2-5}
			$\delta_1$ & 5 & 10 & 30 & 50 \\ \hline
			0.10    & -0.41 & -3.07 & -23.89 & -43.57 \\
			0.20    & -0.28 & -1.56 & -18.95 & -39.12 \\
			0.30    & -0.58 & -0.63 & -11.40 & -31.14 \\
			0.40    & -1.22 & -0.65 & -2.91  & -17.39 \\
			0.45   & -1.57 & -1.09 & -0.57  & -8.76  \\
			0.49   & -1.82 & -1.50 & -0.43  & -3.01 \\ \hline
		\end{tabular}
		\end{center}
	\end{minipage}
	\hfill
	\begin{minipage}{0.5\hsize}
		\begin{center}
		\vspace{3mm}
		{(b) WGRV}\\
		\vspace{3mm}
			\begin{tabular}{ccccc}\hline
				\multicolumn{1}{c}{} & \multicolumn{4}{c}{Intensity of jumps ($\lambda$) }\\
				\cline{2-5}
				$\delta_1$ & 5 & 10 & 30 & 50 \\ \hline
				0.10    & 2.03  & 5.88  & 60.71 & 208.18 \\
				0.20    & 1.18  & 3.43  & 22.96 & 66.24  \\
				0.30    & 0.43  & 1.76  & 10.27 & 29.11  \\
				0.40    & -0.32 & 0.60  & 4.30  & 11.13  \\
				0.45   & -0.80 & 0.00  & 2.74  & 6.26   \\
				0.49   & -1.11 & -0.46 & 1.85  & 3.83\\  \hline
			\end{tabular}
		\end{center}
	\end{minipage}
\label{Table02}
\end{table}
\renewcommand{\arraystretch}{1.0}

\subsubsection{The case of constant volatility}
Since we assumed that the volatility is location-dependent in the previous 
sections, the normalization by estimated spot volatilities is needed to 
obtain an accurate estimator. 
However, if the true volatility of data is constant, 
we may ignore normalization. 

Here we set $\eta = 0$ so that the data is driven by a constant-volatility diffusion process. The intensity is $\lambda = 30$. 
{
\colorg
The summary table of estimated values are shown in Table \ref{Table06}. 
Obviously, all types of GRV and WGRV outperform other estimators. 
}

Figure \ref{fig4} shows the {\colorg error ratios} of this case. 
The GRVs without normalization ({\tt grv[0.20]}, {\tt grv[0.10]} and {\tt grv[0.05]}) perform as well as those with normalization. 
This suggests that, if the true process can be thought as constant-volatility, 
we may skip normalization (calculation of spot volatilities) procedure. 

However, it would be more typical that the volatility is non-constant. 
Thus, basically, it would be advisable to use normalization.

\renewcommand{\arraystretch}{1.5}
\begin{table}[H] \centering 
	\caption{Summary table of estimated values: $\lambda = 30$}
	\vspace{2mm}
	\begin{tabular}{@{\extracolsep{5pt}} lccccc}  
		\hline 
		& Min. & 1st Qu. & Median & 3rd Qu. & Max. \\ 
		\hline \hline
		{\tt trv[0.45]} & 0.41 & 0.44 & 0.45 & 0.46 & 0.49 \\ 
		{\tt trv[0.20]} & 0.97 & 1.11 & 1.16 & 1.22 & 1.45 \\ 
		{\tt trv[0.10]} & 1.67 & 2.49 & 2.76 & 3.08 & 4.81 \\ 
		{\tt bv} & 1.16 & 1.43 & 1.53 & 1.63 & 2.52 \\ 
		{\tt mrv} & 0.92 & 1.03 & 1.07 & 1.13 & 2.52 \\ 
		{\tt grv.lgrv[0.20]} & 0.94 & 1.01 & 1.05 & 1.07 & 1.17 \\ 
		{\tt grv.mrv[0.20]} & 0.94 & 1.02 & 1.05 & 1.09 & 1.32 \\ 
		{\tt wgrv.lgrv[0.20]} & 0.95 & 1.02 & 1.05 & 1.08 & 1.16 \\ 
		{\tt wgrv.mrv[0.20]} & 0.95 & 1.02 & 1.05 & 1.08 & 1.17 \\ 
		{\tt grv[0.20]} & 0.94 & 1.01 & 1.05 & 1.07 & 1.17 \\ 
		{\tt grv[0.10]} & 0.95 & 1.02 & 1.06 & 1.08 & 1.16 \\ 
		{\tt grv[0.05]} & 0.96 & 1.04 & 1.07 & 1.09 & 1.17 \\ 
		{\tt grv.lgrv.mov} & 0.95 & 1.02 & 1.05 & 1.08 & 1.16 \\ 
		{\tt wgrv.lgrv.mov} & 0.95 & 1.02 & 1.05 & 1.08 & 1.15 \\ 
		\hline \hline
		{\tt True Value} & 1.00 & 1.00 & 1.00 & 1.00 & 1.00 \\ 
		\hline 
	\end{tabular} 
\label{Table06}
\end{table} 
\renewcommand{\arraystretch}{1.0}

\begin{figure}[H]
	\begin{tabular}{c}
		\begin{minipage}{0.4\hsize}
			\includegraphics[keepaspectratio, scale=0.7]{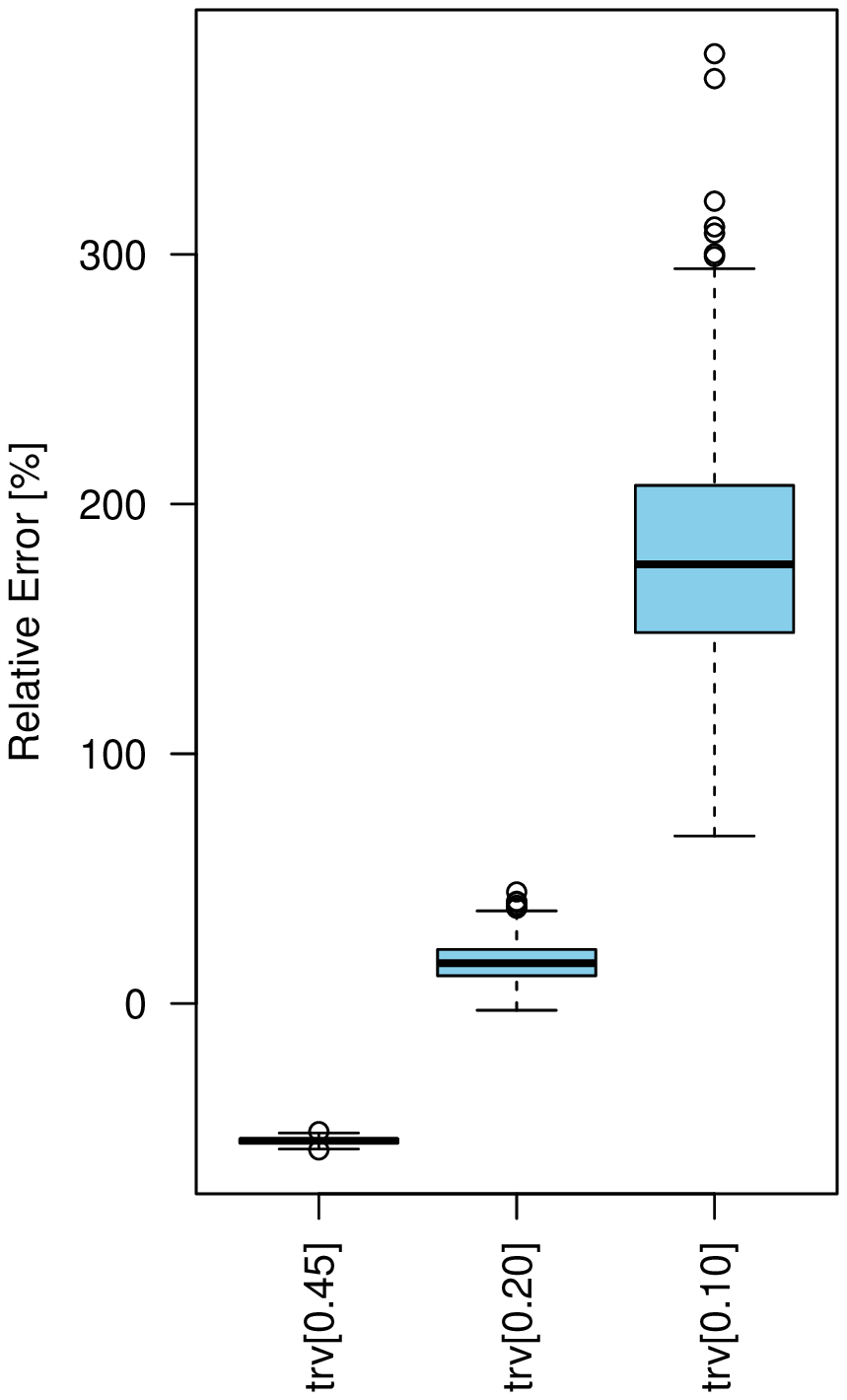}
		\end{minipage}
		
		\begin{minipage}{0.6\hsize}
			\includegraphics[keepaspectratio, scale=0.7]{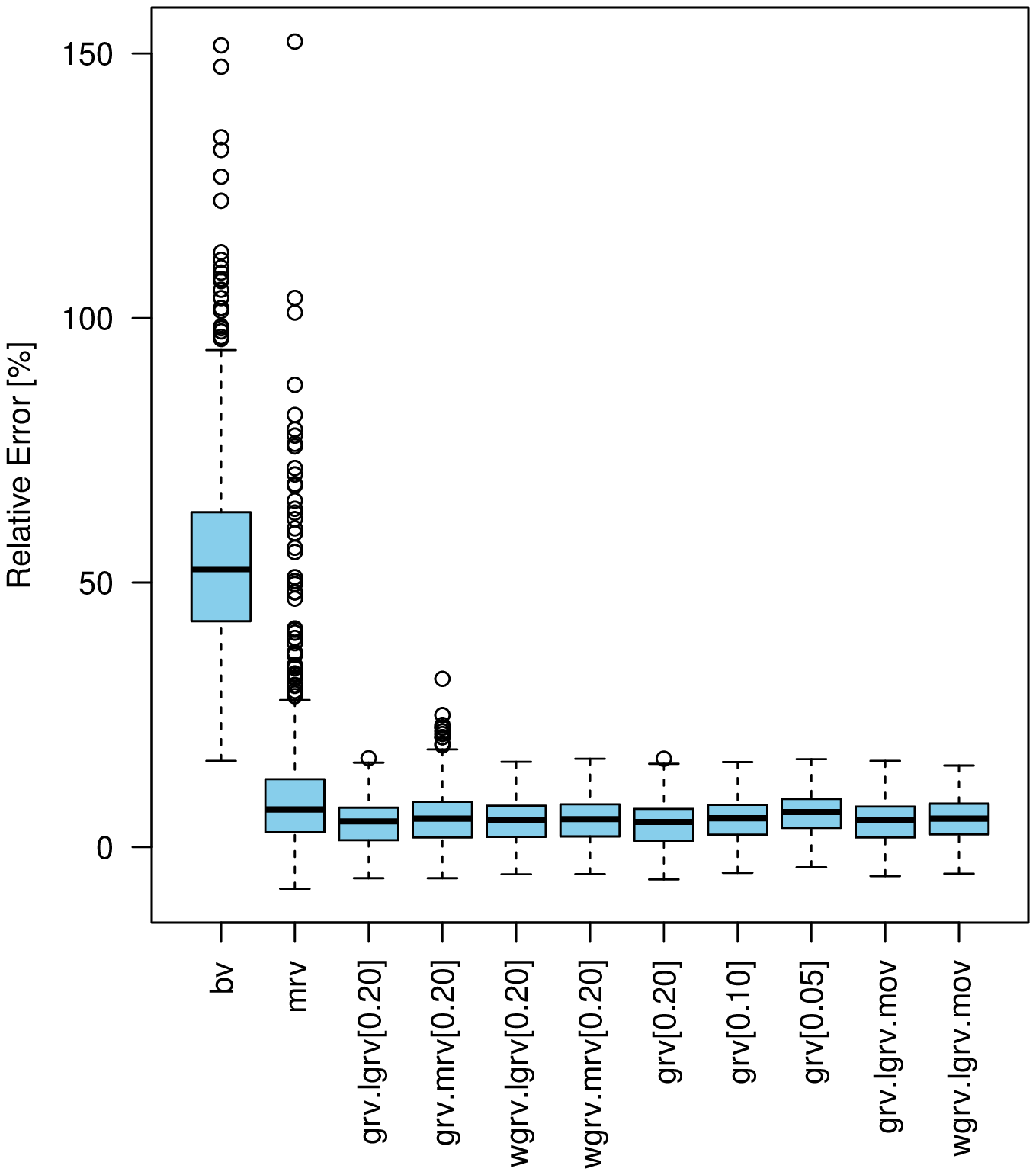}
		\end{minipage}
	\end{tabular}
	\caption{Error ratios [\%] results for the constant volatility: $\lambda = 30$}
	\label{fig4}		
\end{figure}

\subsection{The case of Neyman-Scott type clustering jumps}
As the previous examples show, the minRV performs relatively well in the case of 
compound Poisson type jumps. 
However, even if the intensity of jumps is small, the minRV 
may suffer from an upward bias depending on the structure of jumps.  
In particular, if there are consecutive jumps (which is quite rare 
for compound Poisson processes), the minRV loses it advantage. 
Here we show an example of such a situation. 

We consider the case that the data-generating process is given by 
{\tred $X=U+J$,} 
where {\tred $U$} is the continuous part and $J$ is the jump part. 
Here we assume that $J$ {\tred is a {\it marked Neyman-Scott clustering process}} 
({\tred simply denoted by} NS hereafter), instead of a compound Poisson process. 

{\tred The} NS process is a typical point process representing consecutive jumps. 
That is, there may be jumps within some consecutive intervals. 
This leads to upward bias of BV and minRV because the both of two adjacent increments can consist of large jumps. 
The NS process is constructed as follows. 

{\colorb 
	\bd
	
	\im[(1)] 
	Set ``centers" on the time interval $[0, {\tred 1}]$ by a Poisson process $(N_t^0)$ with intensity $\lambda_0$. A center is defined as the point $t \in [0, {\tred 1}]$ which 
	satisfies $\Delta N_t = 1$. 
	
	\im[(2)] 
	For each center $c \in [0, {\tred 1}]$, choose the number $N_c$ of ``children," 
	assuming $N_c$ is Poisson-distributed with mean $\lambda_c$. 
	
	\im[(3)]
	For each center $c \in [0, {\tred 1}]$, generate independently and exponentially distributed random variables $\left( v_i^{(c)} \right)_{1 \leq i  \leq N_c}$ 
	with mean $h$. 
	Then the location of child $i$ derived from center $c$ is defined as $c - v_i^{(c)}$. 
	This defines the location of a jump.
	
	\im[(4)]
	For each child $i$, generate an independently and normally distributed random variable $\xi_i \sim N(0, \nu_J^2)$. This determines the size and direction of a jump $\Delta J_s$.
	
	\im[(5)]
	The NS process is defined as $J_t = \sum_{s \in [0, t]} \Delta J_s$. 
	
	\ed
}

We 
{\tred generate $X=U+J$,} 
where ${\tred U}$ is the Brownian semimartingale {\tred independent of $J$}, 
satisfying the {\tred stochastic differential equation
	\begin{equation}
		dU_t = \theta\> U_t dt + 
		(\sigma + \eta\> U_t^2 )^{\frac{1}{4}} dw_t
	\end{equation} 
	with $U_0=1$. 
}
We set $\lambda_0 = \lambda_c = 5$ and $\nu_J = 0.5$. 
For the continuous part {\colorg $U$}, we use $\theta = 0.2, \ \sigma = 1, \ \eta = 3$. 
{\tred As before, the number $n$ of samples is $n = 2000$, and the number 
of trials is 500. }

\renewcommand{\arraystretch}{1.5}
\begin{table}[H] \centering 
	\caption{Summary table of error ratios: Neyman-Scott clustering jumps} 
	\vspace{2mm}
	\begin{tabular}{@{\extracolsep{5pt}} lccccc} 
		\hline  
		& Min. & 1st Qu. & Median & 3rd Qu. & Max. \\ 
		\hline \hline
		{\tt trv[0.45]} & -97.04 & -88.67 & -82.29 & -75.74 & -59.83 \\ 
		{\tt trv[0.20]} & -74.35 & -29.07 & -11.09 & 6.46 & 138.67 \\ 
		{\tt trv[0.10]} & -68.02 & -14.24 & 3.17 & 24.43 & 157.66 \\ 
		{\tt bv} & -54.60 & 5.52 & 27.19 & 69.17 & 369.71 \\ 
		{\tt mrv} & -67.31 & -1.40 & 19.59 & 61.35 & 300.83 \\ 
		{\tt grv.lgrv[0.20]} & -74.39 & -31.45 & -14.02 & 3.64 & 136.14 \\ 
		{\tt grv.mrv[0.20]} & -70.53 & -26.44 & -9.19 & 8.39 & 139.11 \\ 
		{\tt wgrv.lgrv[0.20]} & -74.49 & -31.31 & -13.35 & 4.32 & 137.83 \\ 
		{\tt wgrv.mrv[0.20]} & -74.48 & -30.70 & -12.73 & 4.23 & 136.96 \\ 
		{\tt grv[0.20]} & -74.89 & -33.38 & -16.67 & 0.73 & 134.64 \\ 
		{\tt grv[0.10]} & -74.70 & -32.63 & -15.03 & 1.88 & 136.44 \\ 
		{\tt grv[0.05]} & -74.38 & -31.68 & -14.04 & 3.40 & 136.67 \\ 
		{\tt grv.lgrv.mov} & -74.50 & -31.63 & -13.79 & 3.84 & 138.23 \\ 
		{\tt wgrv.lgrv.mov} & -74.32 & -30.82 & -12.89 & 4.34 & 139.78 \\ 
		\hline
	\end{tabular} 
	\label{Table10} 
\end{table} 
\renewcommand{\arraystretch}{1.0}

{\colorg Table \ref{Table10} and Figure \ref{fig5} show the error ratios} in the case of NS jumps. 
Because of the possible consecutive jumps, both bipower variation and 
minRV have upward bias, whereas {\colorb GRV and WRGV} are all robust to such clustering jumps. 
This suggests that the  {\colorb GRV and WRGV perform} very well for various 
structures of jumps. 

\begin{figure}[H]
	\includegraphics[keepaspectratio, scale=0.7]{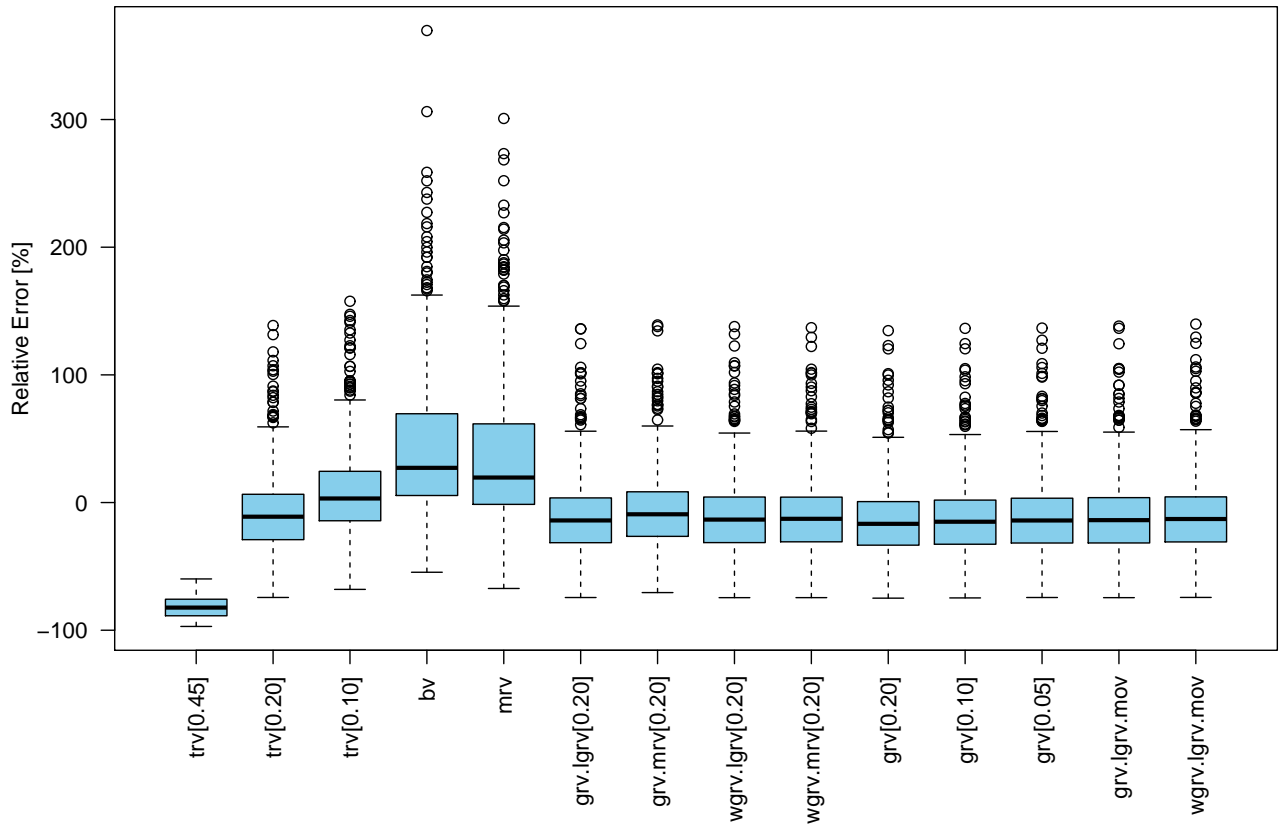}
	\caption{Error ratios [\%] for the case of Neyman-Scott clustering jumps}
	\label{fig5}
\end{figure}

{\colorb

\subsection{A remark on estimation of spot volatilities} \label{SubSec7-3}
Finally, we argue how estimation of spot volatilities affect the 
accuracy of GRV and WGRV. 

We have used $\kappa_n = \lfloor B n^c \rfloor = \lfloor 10n^{0.45} \rfloor = 305$ for local GRV and local minRV and seen that GRV and WGRV with these spot volatilities perform highly well. 
However, the choice of $\kappa_n$ may affect the accuracy of 
GRV and WGRV. 
In fact, if the true volatility varies greatly, a wide subinterval (a large $\kappa_n$) 
leads to imprecise estimation of spot volatilities and causes 
misdetection of jumps by using such information. 
Therefore, it ends up obtaining biases of GRV and WGRV. 

To see this, consider the following SDE: 
\begin{align*}
	dX_t = \theta X_t dt + (\sigma + \eta \sin^2 X_t) dw_t + dJ_t, 
\end{align*}   
where $J_t = \sum_{j=1}^{N_t} \xi_j$ is the same compound Poisson process with intensity $\lambda$ as in Section \ref{SubSec7-1}. 
We set $\sigma = 1, \eta = 5, \lambda = 10, \mu = 0.3,  \nu = 0.2$. 
{\tred Again, the number $n$ of samples is $n = 2000$, and the number 
	of trials is 500. }
{\tred In this example}, the volatility $(\sigma + \eta \sin^2 X_t)^2$ swings in the range $[1, 36]$. 
A sample path of this model is shown in Figure \ref{fig6}. 
The volatility alternates between low and high in short time intervals, so the 
estimation of spot volatility requires an appropriate choice of $\kappa_n$. 
\begin{figure}[H]
	\begin{center}
		\begin{tabular}{c}
			
			\begin{minipage}{0.5\hsize}
				\begin{center}
					\includegraphics[keepaspectratio, scale=0.55]{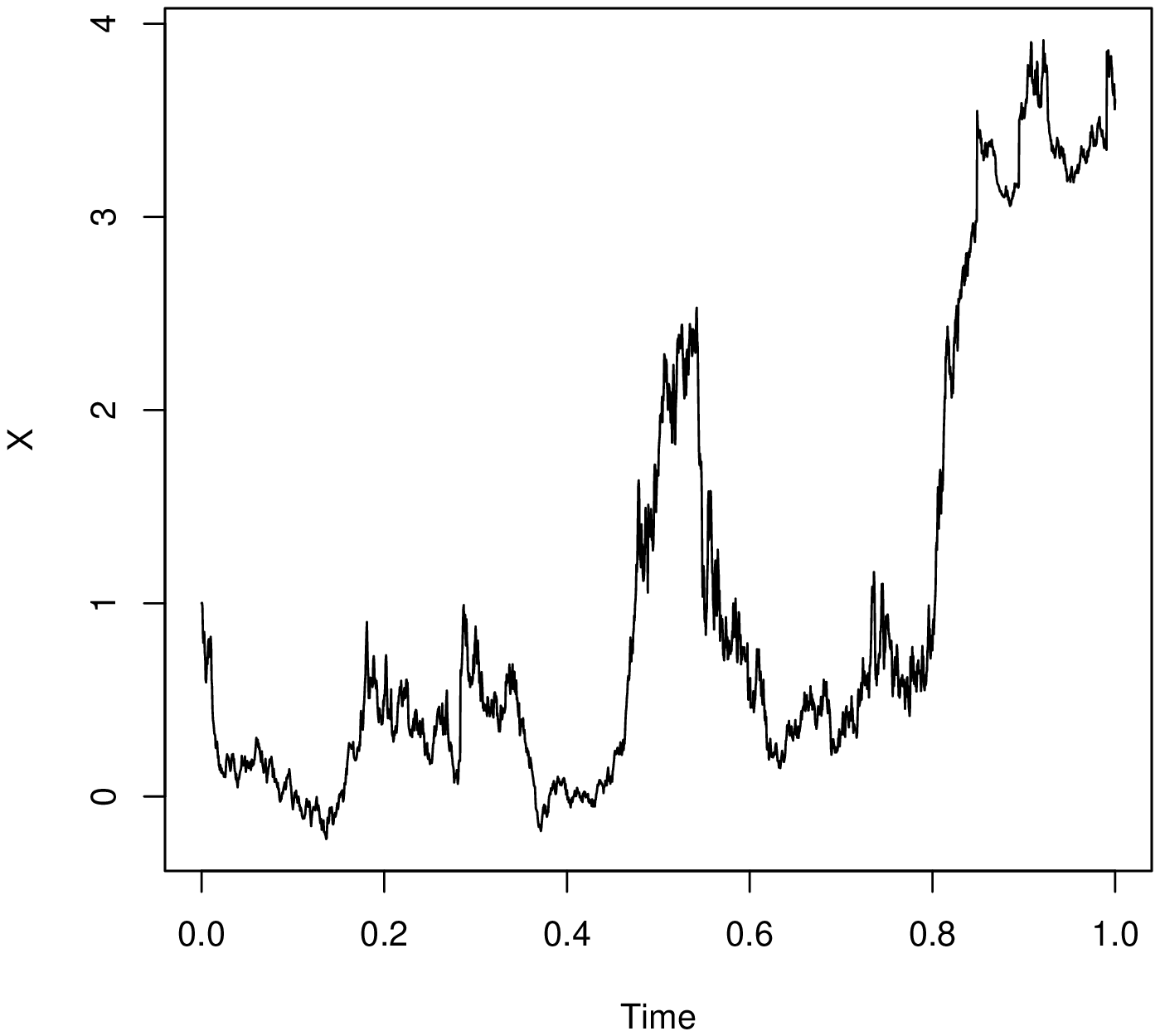}
					{(a) Sample path of $X$}
				\end{center}
			\end{minipage}
			
			\begin{minipage}{0.5\hsize}
				\begin{center}
					\includegraphics[keepaspectratio, scale=0.55]{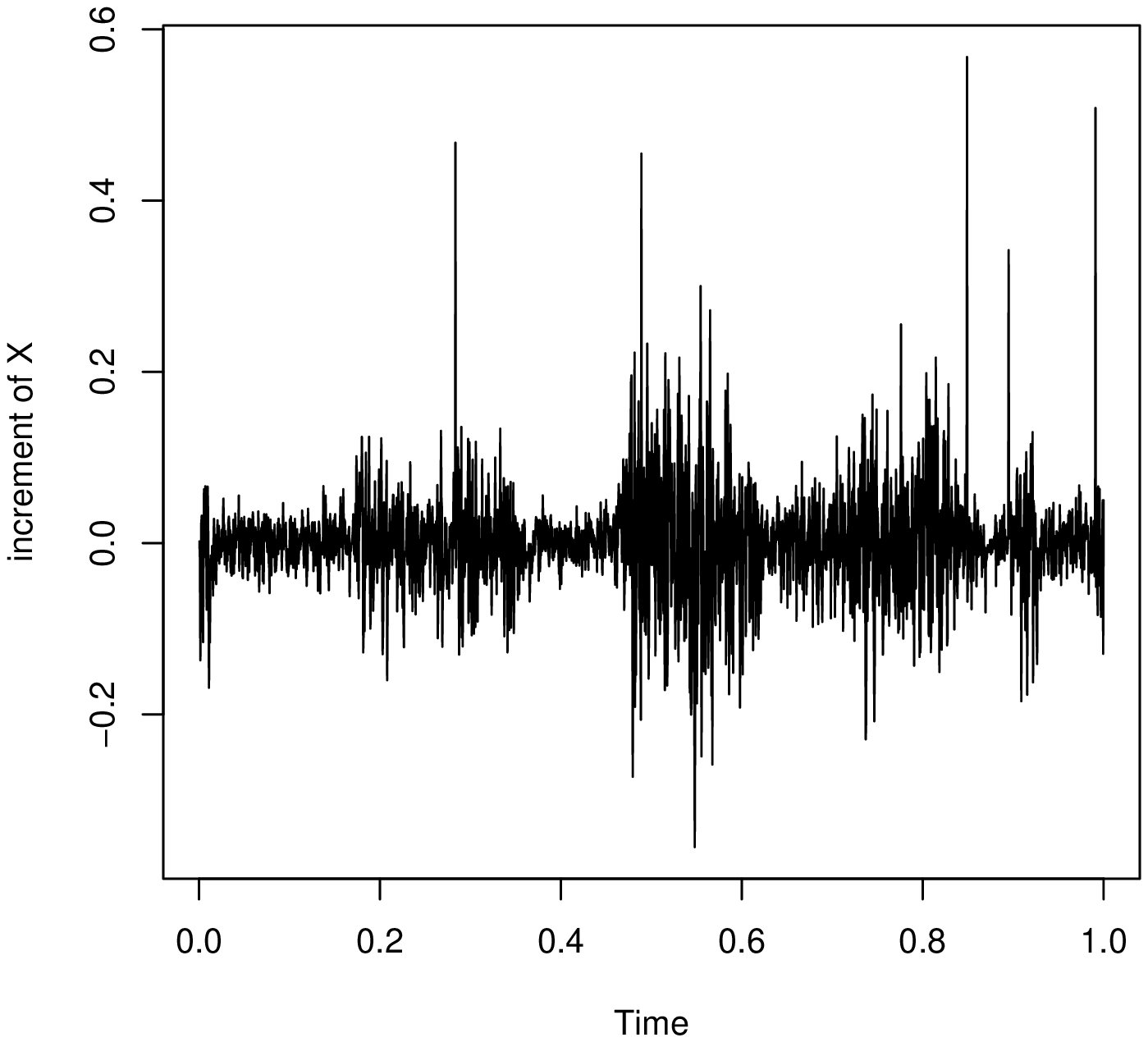}
					{(b) Increment of $X$}
				\end{center}
			\end{minipage}
			
		\end{tabular}
		\caption{Sample path of $X$ and its increments}
		
		\label{fig6}		
	\end{center}
	
\end{figure}

{\colorg Table \ref{Table03} shows the summary and average error ratios of GRV and WGRV, respectively.} 
for several values of $c$ and $B$ that determine the width $\ol{\kappa}_n = 2\kappa_n + 1$ of subintervals for 
spot volatility estimation. 
This indicates that large $B$ and $c$ (wide subinterval) tend to give imprecise estimates. 
Since the volatility varies in a wide subinterval as Figure \ref{fig6} shows, 
the estimated spot volatility is prone to deviate the true value. 
This leads to misdetection of jumps, and thus distorts the estimate of GRV and WGRV. For instance, an underestimated 
spot volatility makes normalized increments too large, so the increments 
are likely to be regarded as jumps and eliminated from calculation of the estimates. As a result, GRV and WGRV are underestimated.   
In this example, it seems that small values such as $c=0.1,  0.2$ and $B = 1, 5$ are 
preferable. 

This example suggests that we should choose the tuning parameters $B$ and $c$ carefully, 
especially when volatility switches between high and low states frequently. 
After all, the proper choice of tuning parameters, such as $B$ and $c$, while observing the data in detail, is needed to obtain precise estimates by GRV and WGRV.

\renewcommand{\arraystretch}{1.5}
\begin{table}[H]
	\caption{Average error ratios of GRV and WGRV for different $c$ and $B$ determining the width $\ol{\kappa}_n$}
	
		\vspace{3mm}
	\begin{minipage}{0.5\hsize}
		\begin{center}
			\vspace{3mm}
			{(a) GRV with local GRV ({\tt grv.lgrv})}\\
			\vspace{3mm}
			\begin{tabular}{ccccc}\hline
				\multicolumn{1}{c}{} & \multicolumn{4}{c}{$B$}\\
				\cline{2-5}
				$c$ & 1 & 5 & 10 & 20 \\ \hline
				0.10 & 3.41   & -5.80  & -9.31  & -15.17 \\
				0.20 & -1.19  & -9.31  & -15.17 & -24.28 \\
				0.30 & -5.26  & -16.48 & -26.09 & -36.26 \\
				0.40 & -9.31  & -27.63 & -37.64 & -44.13 \\
				0.45 & -12.36 & -33.71 & -42.02 & -46.39 \\
				0.49 & -15.46 & -37.92 & -44.30 & -47.37 \\ \hline
			\end{tabular}
		\end{center}
	\end{minipage}
	\hfill
	\begin{minipage}{0.5\hsize}
		\begin{center}
			\vspace{3mm}
			{(b) GRV with local minRV ({\tt grv.mrv})}\\
			\vspace{3mm}
			\begin{tabular}{ccccc}\hline
				\multicolumn{1}{c}{} & \multicolumn{4}{c}{$B$}\\
				\cline{2-5}
				$c$ & 1 & 5 & 10 & 20 \\ \hline
				0.10 & 6.21   & -4.21  & -6.92  & -11.85 \\
				0.20 & 0.69   & -6.92  & -11.85 & -20.91 \\
				0.30 & -3.85  & -13.08 & -22.75 & -33.82 \\
				0.40 & -6.92  & -24.44 & -35.37 & -43.14 \\
				0.45 & -9.33  & -31.01 & -40.38 & -46.01 \\
				0.49 & -12.11 & -35.72 & -43.37 & -47.22 \\ \hline
			\end{tabular}
		\end{center}
	\end{minipage} \\ 

	\vspace{3mm}
	\begin{minipage}{0.5\hsize}
	\begin{center}
	\vspace{3mm}
	{(c) WGRV with local GRV ({\tt wgrv.lgrv})}\\
	\vspace{3mm}
	\begin{tabular}{ccccc}\hline
		\multicolumn{1}{c}{} & \multicolumn{4}{c}{$B$}\\
		\cline{2-5}
		$c$ & 1 & 5 & 10 & 20 \\ \hline
		0.10 & 15.27 & 1.33   & -3.69  & -9.64  \\
		0.20 & 8.81  & -3.69  & -9.64  & -17.78 \\
		0.30 & 1.88  & -10.83 & -19.33 & -28.89 \\
		0.40 & -3.69 & -20.74 & -30.20 & -37.17 \\
		0.45 & -6.97 & -26.38 & -34.80 & -39.68 \\
		0.49 & -9.94 & -30.50 & -37.36 & -40.95 \\ \hline
	\end{tabular}
	\end{center}
	\end{minipage}
	\hfill
	\begin{minipage}{0.5\hsize}
	\begin{center}
	\vspace{3mm}
	{(d) WGRV with local minRV ({\tt wgrv.mrv})}\\
	\vspace{3mm}
	\begin{tabular}{ccccc}\hline
		\multicolumn{1}{c}{} & \multicolumn{4}{c}{$B$}\\
		\cline{2-5}
		$c$ & 1 & 5 & 10 & 20 \\ \hline
		0.10 & 22.42 & 3.91   & -1.64  & -7.54  \\
		0.20 & 13.26 & -1.64  & -7.54  & -15.38 \\
		0.30 & 4.77  & -8.71  & -16.88 & -26.58 \\
		0.40 & -1.64 & -18.29 & -28.05 & -35.98 \\
		0.45 & -4.96 & -23.98 & -33.10 & -39.12 \\
		0.49 & -7.78 & -28.37 & -36.22 & -40.71 \\ \hline
	\end{tabular}
	\end{center}
	\end{minipage}

	\label{Table03}
\end{table}
\renewcommand{\arraystretch}{1.0}
}

\section{Concluding remarks}
In this paper, we construct the global realized volatility estimator 
in the nonparametric context. 
We proved the consistency and the asymptotic normality 
of  
GRV and WGRV, and, by numerical simulations, we show that these new approaches outperform 
previous studies which use increments within a single or two intervals.  

Our new approach for eliminating jumps is highly versatile. 
For example, by normalization, it works well when the volatility of data 
is driven by a nonconstant-volatility process. 
Moreover, both GRV and WGRV 
are accurate enough in the case of not only compound-Poisson sporadic jumps but also Neyman-Scott consecutive jumps.  

The global-filtering method could be extended to the covariance estimation 
even under the nonsynchronous sampling scheme. 
Furthermore, this approach could also be applied to construct 
a test statistic for jump.  
Also, it is valuable to apply our approach to empirical research of 
high-frequency time series data. 
These are important topics for future research. 


\bibliographystyle{spmpsci}      
\bibliography{bibtex-20180615-20191212-20200312-20200418-r+}  

\end{document}

%% file: nakamacro291103+++.tex
\def\infm{{\infty\text{--}}}
\def\inftym{\infm}
\def\koko{{\coloroy{koko}}}
\def\bd{\begin{description}}
\def\ed{\end{description}}

\def\D2{\bbD_{2,\infty-}}

\def\tj{{t_j}}
\def\tjm{{t_{j-1}}}

\def\D{{\bf D}}

\def\cald{{\cal D}}

\def\calf{{\cal F}}

\def\calj{{\cal J}}
\def\calk{{\cal K}}
\def\call{{\cal L}}
\def\calm{{\cal M}}

\def\cals{{\cal S}}

\def\calv{{\cal V}}

\def\ds{\displaystyle}
\def\yeq{\>=\>}
\def\yleq{\>\leq\>}

\def\simleq{\ \raisebox{-.7ex}{$\stackrel{{\textstyle <}}{\sim}$}\ }

\def\ep{\epsilon}
\def\half{\frac{1}{2}}

\def\cadlag{c\`adl\`ag\ }
\def\up{\uparrow}
\def\down{\downarrow}

\def\y{\vspace*{3mm}\\}
\def\halflineskip{\vspace*{3mm}}
\def\nn{\nonumber}
\def\be{\begin{equation}}
\def\ee{\end{equation}}
\def\bea{\begin{eqnarray}}
\def\eea{\end{eqnarray}}
\def\beas{\begin{eqnarray*}}
\def\eeas{\end{eqnarray*}}
\def\bi{\begin{itemize}}
\def\ei{\end{itemize}}
\def\im{\item}
\def\bd{\begin{description}}
\def\ed{\end{description}}
\newcommand{\bbD}{{\mathbb D}}

\newcommand{\bbF}{{\mathbb F}}

\newcommand{\bbL}{{\mathbb L}}
\newcommand{\bbM}{{\mathbb M}}
\newcommand{\bbN}{{\mathbb N}}
\newcommand{\bbO}{{\mathbb O}}

\newcommand{\bbR}{{\mathbb R}}

\newcommand{\bbV}{{\mathbb V}}
\newcommand{\bbW}{{\mathbb W}}

\newcommand{\bbZ}{{\mathbb Z}}